\documentclass[english,11pt]{elsarticle} 

\usepackage[T1]{fontenc}
\usepackage{float}
\usepackage{amsmath}
\usepackage{amsthm}
\usepackage{amssymb}
\usepackage{amsfonts}
\usepackage{bbm}
\usepackage{graphicx}
\usepackage{epstopdf}
\usepackage[pdftex,colorlinks,linkcolor={blue},citecolor={blue}]{hyperref}

\usepackage{caption}

\usepackage{nomencl}

\usepackage{footnote}
\makesavenoteenv{tabular}

\usepackage{algorithm}
\usepackage{algorithmic}

 \usepackage{pgfplots}
  \pgfplotsset{compat=newest}
  \usetikzlibrary{plotmarks}
  \usetikzlibrary{arrows.meta}
  \usepgfplotslibrary{patchplots}
  \usepackage{grffile}
  \usepackage{amsmath}

\usepackage{epsfig}
\usepackage{subfigure}
\usepackage{amssymb,amsmath,amsfonts,layout,graphicx}
\usepackage{makeidx}
\usepackage{babel}
\usepackage{tikz}
\usepackage{sublabel}
\usepackage{cases}
\usepackage{algorithm}
\usepackage{algorithmic}

\usepackage
[
        letterpaper,
        left=1.91cm,
        right=1.91cm,
        top=1.91cm,
        bottom=2cm,
]
{geometry}

\newtheorem{assumption}{Assumption}
\newtheorem{proposition}{Proposition}

\newtheorem{remark}{Remark}

\DeclareMathOperator*{\argmax}{arg\,max}

\makenomenclature 

\begin{document}
\begin{frontmatter}
 \title{\textbf{Spatial Pricing in Ride-Sourcing Markets under a Congestion Charge}}	

\author[1staddress]{Sen Li}
\ead{cesli@ust.hk}
\author[1staddress]{Hai Yang}
\ead{cehyang@ust.hk}
\author[2ndaddress]{Kameshwar Poolla}
\ead{poolla@berkeley.edu}
\author[2ndaddress]{Pravin Varaiya}
\ead{varaiya@berkeley.edu}
\address[1staddress]{Department of Civil and Environmental Engineering,  Hong Kong University of Science and Technology}
\address[2ndaddress]{Department of Electrical Engineering and Computer Science,  University of California, Berkeley}

\begin{abstract}
This paper studies the optimal spatial pricing for a ride-sourcing platform subject to a congestion charge. The platform determines the ride prices over the transportation network to maximize its profit, while the regulatory agency imposes the congestion charge to reduce traffic congestion in the urban core.  A network economic equilibrium model is proposed to capture the intimate interactions among passenger demand, driver supply, passenger and driver waiting times, platform pricing, vehicle repositioning and  flow balance over the transportation network. The overall optimal pricing problem is cast as a non-convex program. An algorithm  is proposed to approximately compute its  optimal solution, and a tight upper bound is established to evaluate its performance loss with respect to the globally optimal solution. Using the proposed model, we compare the impacts of three forms of congestion charge: (a) a one-directional cordon charge on ride-sourcing vehicles that enter the congestion area; (b) a bi-directional cordon charge on ride-sourcing vehicles that enter or exit the congestion area; (c) a trip-based congestion charge on all ride-sourcing trips. We show that the one-directional congestion charge not only reduces the ride-sourcing traffic in the congestion area, but also reduces the travel cost outside the congestion zone and benefits passengers in  these underserved areas. We establish that in all  congestion charge schemes the largest share of the tax burden is carried by the platforms, as opposed to passengers and drivers. We further show that compared to other congestion charges, the one-directional cordon charge is more effective in congestion mitigation: to achieve the same congestion-mitigation target, it imposes  a smaller cost on passengers, drivers, and the platform. On the other hand, compared with the other  charges, the trip-based congestion charge is more effective in revenue-raising:  to raise the same tax revenue, it leads to a smaller loss to passengers, drivers, and the platform.  We validate these results through realistic numerical studies for San Francisco. 
\end{abstract}


\nomenclature[01]{$\alpha$}{Trade-off between money and waiting time for ride-sourcing passengers}
\nomenclature[02]{$\beta$}{Trade-off between money and in-vehicle travel time for ride-sourcing passengers}
\nomenclature[03]{$\epsilon$}{Parameter of the passenger demand logit model}
\nomenclature[04]{$\sigma$}{Parameter of the driver supply logit model}
\nomenclature[05]{$\eta$}{Parameter of the driver repositioning logit model}
\nomenclature[06]{$\lambda_{ij}$}{Arrival rate (per minute) of ride-sourcing passengers traveling from zone $i$ to zone $j$}
\nomenclature[07]{$\lambda_{ij}^0$}{Arrival rate (per minute) of potential passengers traveling from zone $i$ to zone $j$}
\nomenclature[08]{$\lambda_{11}$}{Arrival rate (per minute) of  passengers traveling from $\mathcal{C}$ to $\mathcal{C}$}
\nomenclature[09]{$\lambda_{12}$}{Arrival rate (per minute) of  passengers traveling from $\mathcal{C}$ to $\mathcal{R}$}
\nomenclature[10]{$\lambda_{21}$}{Arrival rate (per minute) of  passengers traveling from $\mathcal{R}$ to $\mathcal{C}$}
\nomenclature[11]{$\lambda_{22}$}{Arrival rate (per minute) of  passengers traveling from $\mathcal{R}$ to $\mathcal{R}$}
\nomenclature[12]{$\mathbb{P}_{ij}$}{Probability of vehicle repositioning from zone $i$ to zone $j$ before congestion charge}
\nomenclature[13]{$\tilde{\mathbb{P}}_{ij}$}{Probability of vehicle repositioning from zone $i$ to zone $j$ after congestion charge}
\nomenclature[14]{$\Pi(k|i\rightarrow j)$}{Probability that a repositioning vehicle from $i$ to $j$ is intercepted at zone $k$}
\nomenclature[15]{$\sigma_i$}{Probability of getting matched to a passenger while a driver traverses zone $i$}
\nomenclature[16]{$\mathcal{C}$}{Set of indexes of the congested zones}
\nomenclature[17]{$\mathcal{R}$}{Set of indexes of the remote zones}
\nomenclature[18]{$\mathcal{S}$}{Set of indexes of the southwest corner of San Francisco}
\nomenclature[19]{$\mathcal{V}$}{Set of indexes of all the $M$ zones}

\nomenclature[20]{$c_{ij}$}{Generalized travel cost ($\$$) of ride-sourcing passengers from zone $i$ to zone $j$}
\nomenclature[20]{$d_{ij}$}{Average trip distance (mile) from zone $i$ to zone $j$}
\nomenclature[20]{$d^{\mathcal{C}}_{ij}$}{Average trip distance (mile) from zone $i$ to zone $j$ within the congestion area $\mathcal{C}$ }
\nomenclature[20]{$d^{\mathcal{R}}_{ij}$}{Average trip distance (mile) from zone $i$ to zone $j$  within the remote area $\mathcal{R}$}
\nomenclature[21]{$f_{ij}$}{Intended driver repositioning flow (per minute) from zone $i$ to zone $j$}
\nomenclature[22]{$\tilde{f}_{ij}$}{Actual driver repositioning flow (per minute) from zone $i$ to zone $j$}
\nomenclature[23]{$L$}{Parameter of the pickup time model}
\nomenclature[24]{$M$}{Number of zones in the city} 
\nomenclature[25]{$N$}{Vehicle hours supplied by the ride-sourcing drivers}
\nomenclature[26]{$N_0$}{Potential vehicle hours supplied by the drivers}
\nomenclature[27]{$N_i^I$}{Average number of idle drivers (idle vehicle hours)  in zone $i$}
\nomenclature[28]{$N_\mathcal{C}$}{Total number of ride-sourcing vehicles in the congestion area.}
\nomenclature[29]{$N_\mathcal{R}$}{Total number of ride-sourcing vehicles in the remote area.}
\nomenclature[30]{$N_\mathcal{S}$}{Total number of ride-sourcing vehicles in the southwest corner of San Francisco.}
\nomenclature[31]{$q$}{Average driver wage per unit time}
\nomenclature[32]{$r_i$}{Ride fare (\$/min) for ride-sourcing trips starting from zone $i$}
\nomenclature[33]{$t_{ij}$}{Travel time (minute) from zone $i$ to zone $j$}
\nomenclature[34]{$w_i^p$}{Average waiting time (minute) for passengers starting from zone $i$}
\nomenclature[35]{$w_i^d$}{Average waiting time (minute) for drivers in zone $i$ to pick up the next passenger}
\nomenclature[36]{$v_c$}{Average speed (mph) in the congestion area.}
\nomenclature[37]{$v_r$}{Average speed (mph) in the remote area.}

\end{frontmatter}

\printnomenclature

\section{Introduction}
Ride-sourcing platforms like Uber, Lyft and Didi, have greatly changed the way many people commute in the city. They offer on-demand mobility services to hundreds of millions of passengers from anywhere at anytime. They generate substantial travel demand that would not have existed. They replace many trips that otherwise would have been taken by public transit, taxis and walking.   They also raise public concerns about congestion, emissions, and social equity.

As ride-sourcing platforms keep disrupting  urban transportation, some city governments have responded by imposing various congestion-mitigation policies. In Feb 2019,  New York City (NYC) passed a Tax Law  imposing a \$2.75/trip congestion surcharge on all ride-sourcing trips that begin in, end in, or pass through the congestion area of NYC in Manhattan, south of and excluding 96th Street \cite{NYC2019surcharge}. In Apr 2019, the New York State Assembly passed the Traffic Mobility Act (Bill No. 09633)  authorizing the Metropolitan Transportation Authority (MTA) to charge a cordon price between \$10 and \$15 on vehicles entering  south of 61st Street in Manhattan by 2021 \cite{NYC_cordon}. It is expected that this cordon price will raise \$3.5B/year for MTA to upgrade its aging transit infrastructure.  In Jan 2020, the City of Chicago introduced a tiered congestion surcharge on ride-sourcing trips \cite{chicago_surcharge}: a \$3/trip charge on solo rides that start or end in the designated downtown area, a \$1.25/trip charge on other solo rides, and a \$0.65/trip charge on other shared rides. The revenue collected from this tax will be used to subsidize public transit and cab drivers. In addition to NYC and Chicago,  San Francisco imposed a Traffic Congestion Mitigation excise tax of 1.5\% to 3.25\%  (effective January 1, 2020) on fares for rides originating in San Francisco that are facilitated by commercial ride-share companies or are provided by an autonomous vehicle or private transit services vehicle \cite{SFmayor}. Massachusetts \cite{Massa_cordon} is proposing similar congestion-mitigating measures  for legislative approval. Bottom line: the regulatory landscape is changing. 

Despite  congestion-limiting regulations that target  ride-sourcing traffic, there is still a debate  whether ride-sourcing  platforms actually  contribute to  increased traffic congestion. On the one hand, various studies have demonstrated the positive correlation between ride-sourcing platforms and traffic congestion. For instance,  the San Francisco County Transportation Authority \cite{castiglione2018tncs} estimated that ride-sourcing platforms contributed 50\% of the increase in traffic congestion in San Francisco between 2010 and 2016. Schaller \cite{schaller2017empty} reported that the number of ride-sourcing vehicles  increased by 59\% in NYC between 2013 and 2017, while  average traffic speed declined by 15\% and  vehicle miles traveled  (VMT) increased by 36\%. A recent study \cite{qian2020impact} collected  ride-sourcing data in NYC and concluded that the growth of ride-sourcing platforms is the major contributing factor making urban traffic congestion worse. On the other hand, another study  concluded that ride-sourcing vehicles do not account for a significant portion of the overall traffic. Fehr \& Peers examined the combined VMT by Uber and Lyft in six metropolitan cities in the US \cite{balding2019}. They showed that the VMT of Uber and Lyft is vastly outstripped by personal and commercial vehicles,   accounting only for 13.4 percent of VMT in San Francisco County, 8 percent in Boston and 7.2 percent in Washington, DC. 
However, we emphasize that a city's traffic congestion exhibits a strong spatial and temporal pattern, with most of the congestion concentrated in a small area of the city for short periods of time. Moreover, the spatial and temporal pattern of the city's traffic is highly correlated with the demand and supply pattern of the ride-sourcing market. Therefore using  average VMT throughout the city  as the indicator of congestion will lead to a significant underestimate of the congestion caused by  ride-sourcing vehicles. Hence the spatial and temporal aspects are crucial in estimating and addressing these congestion externalities.

This paper focuses on the spatial aspect of the ride-sourcing market and proposes a  model to predict the impact of congestion charges on ride-sourcing services over a transportation network. An economic equilibrium model is formulated to capture the opposing incentives of passengers, drivers, and the platform. The model is built upon a traffic network to characterize the spatial distribution of passengers and drivers  in the ride-sourcing market. By incorporating the spatial aspect in the economic model, our framework can capture the intimate interactions among  passenger spatial distribution, driver spatial distribution, passenger waiting time, driver waiting time, idle driver repositioning, { traffic congestion}, network flow balance and location-differentiated pricing. The major contributions of this paper are summarized below:

\begin{itemize}

\item We develop a networked market equilibrium model that simultaneously captures  the essential factors in a ride-sourcing network, including passenger waiting time, driver waiting time, passenger spatial distribution, driver spatial distribution, driver repositioning strategies, { traffic congestion}, network flow balance and platform pricing. 

\item We cast the optimal spatial pricing problem as a non-convex program and  propose an efficient algorithm to approximately compute its optimal solutions. A tight upper bound is established to evaluate the performance of the proposed algorithm. 

\item We evaluate the impact of three forms of congestion charge: (a) a one-directional cordon price that penalizes all vehicles (including idle ones)  entering the congestion area, (b) a bi-directional cordon price that penalizes all vehicles (including idle ones) for entering or exiting the congestion area, and (c) a trip-based congestion charge on all ride-sourcing trips.\footnote{The trip-based congestion charge is imposed on a per-trip basis, and exempts idle vehicles cruising for passengers.} We show that the one-directional congestion charge not only reduces ride-sourcing traffic in the congestion area, but also reduces the travel cost outside of the congestion area and benefits passengers in these underserved zones. We verify that for all  three congestion charges the highest share of the tax burden is borne by the platform as opposed to passengers and drivers. We also show that compared with other congestion charges, the one-directional cordon price is most effective in congestion-mitigation, while the trip-based congestion charge is most effective in revenue-raising: to achieve the same traffic mitigating target (or revenue-raising target), the one-directional charge (or trip-based congestion charge) imposes a smaller cost on passengers, drivers and the platform.

\item The proposed framework is validated through realistic simulations  based on synthetic ride-sourcing data for San Francisco. Through sensitivity analysis, we show that the results and insights derived from our model are robust with respect to the variation of model parameters.

\end{itemize}

\section{Related work}
A growing literature  addresses the supply-demand equilibrium of mobility services over a transportation network.  In a pioneering work \cite{yang1998network}, a network model was developed to describe the movement of idle and occupied taxi vehicles on the road network to look for passengers and provide transportation services. It offers interesting insights on the interactions among the number of taxis, the average taxi utilization and the passenger waiting time. This model was later extended to incorporate various other  features, such as traffic congestion and demand elasticity \cite{wong2001modeling}, competition and regulation \cite{yang2002demand}, bilateral taxi-customer searching and meeting \cite{yang2010equilibria}, and multiple user classes and vehicle modes \cite{wong2008modeling}. { While these works primarily focus on street-hailing of taxi services, they were later adapted to the study of e-hailing services}. For instance, \cite{xu2019equilibrium} considered a traffic network model with (e-hailing) ride-sourcing services, which captured both intra-zone matching and inter-zone matching to account for cruising and deadheading vacant trips at the same time. By numerical studies, they showed that neglecting inter-zone matching may lead to significant bias in evaluating the vacancy/empty miles generated by the ride-sourcing platforms. \cite{ban2019general} developed a general economic equilibrium model to describe the equilibrium state of a transportation system with solo drivers and the ride-sourcing platform, which provides insights to understand the relationship between ride-sourcing service, deadhead miles, and its impact on congestion. \cite{ghili2020spatial} proposed a theoretical model to identify the spatial distortion of supply from demand, and showed that a smaller platform tends to distort the supply of drivers towards more densely populated areas due to network effects.

However,  all of the aforementioned works primarily focus on equilibrium analysis for a given price, so the price-formation mechanism is missing. In practice, ride-sourcing platforms use pricing as a powerful tool to balance supply and demand for more efficient matching, which defines the fundamental difference between the emerging ride-sourcing service and traditional taxi service. This has been widely studied in the network setting as the optimal spatial pricing problem, where the platform determines the location-differentiated prices to maximize its profit and passengers and drivers respond to these prices by making mode choices or repositioning decisions. For instance, \cite{bimpikis2019spatial} studied  spatial price discrimination in a ride-sharing platforms that serves a networks of locations, and showed that platform profit and consumer surplus at the equilibrium are maximized when the demand pattern is balanced across the transportation network. \cite{zha2018geometric} developed a discrete-time geometric matching framework over the transportation network to explore the effects of spatial pricing on the ride-sourcing market, and showed that the platform may increase the ride price to avoid inefficient supply if spatial price differentiation is not allowed. \cite{chen2019optimal} built a bi-level model to capture the decision-making of both platforms and ride-sourcing customers and revealed that the platform profit at a particular zone can be influenced by ride requests from other zones. \cite{guda2019your} considered surge pricing in a ride-sourcing market with independent drivers that strategically move between zones and showed that surge price can be useful even in zones where supply exceeds demand. \cite{ma2018spatio} proposed a spatial-temporal pricing mechanism for a multi-period, multi-location model and showed that the mechanism is incentive-aligned, individually rational, budget balanced, and  welfare-optimal. \cite{afifah2020spatial} studied the optimal spatial pricing problem with a Stackelberg framework considering the ride-sourcing platform's congestion externality. \cite{he2015modeling} studies the impact of the ride-sourcing platform on the taxi system by developing a spatial equilibrium model that not only balances the demand and supply of the taxi market but also captures the possible adoption of emerging e-hailing apps by the taxi drivers. \cite{qinpiggyback} considered a fleet of electrified vehicles providing electricity services and mobility services at the same time, and demonstrated its synergetic value in reducing the transportation spatial imbalance.   
{\em We emphasize that our paper differs from these works in that we simultaneously capture locational differentiated prices, passenger waiting time,  driver waiting time, idle driver repositioning, traffic congestion, and network flow balance.  These crucial elements are combined in a spatial pricing model to investigate the impact of congestion charges on the ride-sourcing networks.}

Road pricing has been a  research topic for decades. The idea was initially proposed in \cite{pigou2017economics}, which subsequently inspired several important works including \cite{walters1961theory}, \cite{vickrey1963pricing}, and \cite{beckmann1967optimal}. Various pricing schemes were explored \cite{yang2005mathematical}, \cite{lindsey2000traffic} following these seminal works. For instance, in \cite{may2000effects} the congested assignment network model was used to test four road pricing systems, with charges based on cordon crossed, distance traveled, time spent in traveling and time spent in congestion. They showed that when rerouting effects are considered, the benefit of road pricing is significantly smaller than expected. \cite{zhang2004optimal} considered the joint optimization of toll levels and toll locations on a road network using a mathematical program with mixed variables. \cite{yang2010road} proposed a convergent trial-and-error implementation method for a form of road pricing congestion control under the condition that  both the link travel time and the travel demand are unknown.  \cite{de2005congestion} considered the dynamic road pricing problem by modeling departure time decision and mode and route choices as endogenous variables. Various link-tolling schemes are analyzed using a dynamic network simulator. 
\cite{wu2011pareto}  determined a Pareto-improving pricing scheme for relieving  traffic congestion in a multimodal transportation network that maximizes  social benefit without increasing the travel expenses of the stakeholders. \cite{wu2012design} explored the effect of income on traveler choice and proposed a theoretical model to design efficient and equitable road pricing schemes. 
More recently, road pricing has been investigated for autonomous vehicles \cite{simoni2019congestion, mehr2019pricing, salazar2019intermodal} and ride-sourcing platforms \cite{li2019regulating}, \cite{li2020impact}.  It is important to note that road pricing for ride-sourcing platforms is substantially different from other road pricing schemes as it interferes with the complicated interaction among platform pricing, passenger demand and driver supply, which is a unique feature of the ride-sourcing market.

\section{Optimal Spatial Pricing: Market Equilibrium Model} 
\label{lowerlevel}

This section formulates a mathematical model of the optimal spatial pricing problem for a ride-sourcing platform over a transportation network. The model captures the complicated interactions among various endogenous variables, including ride fare, driver payment, passenger and driver distribution,  passenger and driver waiting time, vehicle repositioning, traffic congestion, and network flow balance.  The details of the model are presented below.

\subsection{Problem Setup}
Consider a city divided into $M$ zones. These zones are connected  by a road network expressed by a graph $\mathcal{G}(\mathcal{V}, \mathcal{E})$, where $\mathcal{V}$ denotes the set of vertices (or zones) and $\mathcal{E}$ denotes the set of edges (or roads). We investigate the ride-sourcing market at the zonal granularity, and assign an origin zone $i\in \mathcal{V}$ and a destination zone $j\in \mathcal{V}$ to each ride-sourcing trip. { Each trip is randomly initiated by a passenger who requests a pickup from zone $i$ through the user app. Upon receiving the request, the platform matches the passenger to the closest idle vehicle in the same zone\footnote{We assume each passenger is matched to the closest idle vehicle. When vehicles are densely distributed and the zones are reasonably large, there is a high probability that passengers are picked up by idle vehicles in the same zone. This paper exclusively considers this case for simplicity.}, if available.  After an idle driver is matched to the passenger, he/she travels to the passenger's location, picks up the passenger, and then chooses the shortest path from $i$ to $j$ to deliver the passenger to the destination. When a trip is completed, the driver can choose either to remain in the same zone or to cruise to a different zone to look for the next passenger. The idle driver remains cruising until he is matched to the next passenger.    
}



\subsection{Passenger Model}
Passengers make mode choices by comparing the costs of different commute options, such as ride-sourcing, public transit, taxi, and walking. For a passenger  traveling from zone $i$ to zone $j$, { the generalized cost of the ride-sourcing trip is defined as the weighted sum of the waiting time, the in-vehicle travel time, and the trip fare:
\begin{equation}
\label{generalizedcost}
	c_{ij} = \alpha {w^p_i}  + \beta t_{ij}+ r_i t_{ij},
	\end{equation}
where $w^p_i$ is the average passenger waiting time in zone $i$, $r_i$ is the per-time ride fare for trips starting from zone $i$, $t_{ij}$ is the average trip time from zone $i$ to zone $j$, $\alpha$ represents the passenger value-of-time when waiting for the ride, and $\beta$ represents the passenger's value of time when traveling on the road. Since travelers typically place a higher value on waiting time \cite{mohring1987values}, we have $\alpha>\beta$\footnote{{The value of passenger waiting time for ridesourcing services may be smaller than that of buses or subway. However, we emphasize that the proposed formulation and solution methodology do not rely on the assumption that  $\alpha>\beta$.}}.  It is worth emphasizing that $t_{ij}$ is an endogenous variable that depends on the traffic flow of the ride-sourcing market (we will delineate this relation in Section \ref{section_congestion}). It is also important to note that  $c_{ij}$ represents the {\em average} value of the generalized cost, while the travel cost for different passengers can be different due to {user} heterogeneity and randomness.  
}

The passenger waiting time $w^p_i$ is an endogenous variable that depends on the  supply and demand of the ride-sourcing market. Since we have assumed that passengers from zone $i$ are matched to drivers in the same zone, the passenger waiting time $w^p_i$ is a monotone function of the average number of idle vehicles (or equivalently, idle vehicle hours  $N^I_i$) in zone $i$.\footnote{In this paper, we regard idle vehicle hours as equivalent to the number of idle drivers. It depends on both supply and demand: it increases with respect to driver supply and decreases with respect to  passenger demand.}  With slight abuse of notation, we use $w^p_i(N^I_i)$ to denote this relation. The following assumption is  imposed:
\begin{assumption}
\({w^p_i}({N_i^I})\) is positive, strictly decreasing with respect to \({N_i^I}\), and $\lim_{N_i^i\rightarrow 0} w_i^p(N_i^I)=\infty$.
\label{assumption1}
\end{assumption}

\begin{remark}
The structure of the trip fare $r_it_{ij}$ closely approximates the industry practice. For instance, both Uber and Lyft have a fixed per-time fare and per-distance fare. They first calculate the total trip fare as the sum of a base fare, a time-based charge, and a distance-based charge, then this total trip fare is multiplied by a surge multiplier that reflects the real-time imbalance between supply and demand in each zone $i$. {The third term $r_it_{ij}$ in (\ref{generalizedcost}) can be viewed as an approximation of the total trip fare multiplied by the surge multiplier. This approximation introduces  error when the traffic speed is not uniform. However, we believe this is a reasonable simplification. Although adding an extra decision variable to represent the distance-based charge can be addressed by the proposed solution methodology, it would complicate the notation without providing extra insights.}
\end{remark}

We assume that passengers choose their transport mode based on the average travel cost of the ride-sourcing trip. The arrival rate of ride-sourcing passengers from zone $i$ to zone $j$ is determined by 
\begin{equation}
\label{demand}
	\lambda_{ij}  = {\lambda _{ij}^0}{F_p}(c_{ij}),
	\end{equation}
where \({\lambda^0_{ij}}\) is the arrival rate of potential passengers from zone $i$ to zone $j$ (total travel demand for all transport modes), and \({F_p}( \cdot )\) is the proportion of potential passengers who choose to take ride-sourcing. We assume that \({F_p}( \cdot )\) is  strictly decreasing with respect to the travel cost $c_{ij}$. Note that (\ref{demand}) includes the logit model as a special case.

{
\begin{remark}
We acknowledge that our model only produces long-term average outcomes for the ride-sourcing market. It does not capture the temporal dynamics (i.e., departure time, arrival time, etc). To address this concern, we can either incorporate the temporal aspect by considering a quasi-static model for each hour, or formulating a fully dynamic model with fine temporal granularity \cite{nourinejad2020ride}. However, when combined with the ride-sourcing network, both methods way render a much more complex temporal-spatial problem that is difficult to address. For this reason, we leave this as future work.  
\end{remark}
}

\subsection{ Driver Model}
In the long-run, drivers are sensitive to earnings and respond to the platform wage by subscribing to or unsubscribing from  the platform. In the short-run drivers decide whether to remain in the current zone or cruise to a different zone to look for the next passenger. 

The long-term driver decisions determine the total driver or vehicle hours in the overall ride-sourcing network as an increasing function of the average wage offered by the platform. The ride-sourcing vehicle hour is given by
\begin{equation}
\label{driversupply}
	N = {N_0}{F_d}(q),
\end{equation}
where $N$ is the total vehicle hour, $q$ is the driver's average hourly wage, $N_0$ is the supply of potential driver or vehicle hours, and $F_d(q)$ is a strictly increasing function representing the proportion of drivers willing to subscribe to the ride-sourcing platform at wage $q$. Note that (\ref{driversupply}) includes the logit model as a special case.

To model short-term driver repositioning decisions, we define $w_i^d$ as  the average driver waiting time in zone $i$. By Little's Law, the vehicle idle hour $N^I_i$ relates to the driver waiting time $w_i^d$ as
\begin{equation}
N^I_i=w_i^d\sum_{j=1}^N \lambda_{ij}.
\end{equation}
Drivers make repositioning decisions by comparing his expected earning (per unit time) in each zone. For each trip that originates from zone $i$, the average total trip fare  depends on the average trip time and the per-time trip fare:
\begin{equation}
\bar{e}_i=r_i\bar{t}_i=r_i\dfrac{\sum_{j=1}^N \lambda_{ij}t_{ij}}{\sum_{j=1}^N \lambda_{ij}}
\end{equation}
where $\bar{e}_i$ is the average trip fare for passengers from zone $i$, and $\bar{t}_i$ is the average trip-time for trips starting from zone $i$. Each idle driver can either remain in zone $i$ or cruise to zone $j$ to search for  the next passenger. If he remains in zone $i$, by the end of his next trip, he will experience an average waiting time of $w_i^d$,  an average trip time of $\bar{t}_i$, and earns a proportion\footnote{The proportion is determined by the commission rate of the platform, which is a uniform number that does not depend on the origin or destination of the trip. } of the average total trip fare $\bar{e}_i$. In this case, his expected earning (per-unit time) is proportional to $\dfrac{\bar{e}_i}{w_i^d+\bar{t}_i}$.  On the other hand, if he cruises to zone $j$ , by the end of his next trip, he will experience an average waiting time of $t_{ij}+w_j^d$, take an average trip time of $\bar{t}_j$, and earns a proportion of $\bar{e}_j$, which leads to a per-unit time earning that is proportional to $\dfrac{\bar{e}_j}{t_{ij}+w_j^d+\bar{t}_j}$. Under a logit model, the probability  of repositioning from zone $i$ to zone $j$ is\footnote{{Note that $\mathbb{P}_{ii}$ differs from $\mathbb{P}_{ij}$ in that it does not include $t_{ii}$. This is because after dropping off the passenger, if the idle driver decides to stay in the same zone,  she/he can directly start cruising for the next passenger, thus $t_{ii}$ should be excluded from his/her cost.} }
\begin{align}
\begin{cases}
&\mathbb{P}_{ij}=\dfrac{e^{\eta \bar{e}_j/ (t_{ij}+w_j^d+\bar{t}_j)}}{\sum_{k\neq i} e^{\eta\bar{e}_k/ (t_{ik}+w_k^d+\bar{t}_k)}+e^{\eta \bar{e}_i/ (w_i^d+\bar{t}_i)}}, \quad  j\neq i, \\
&\mathbb{P}_{ii}=\dfrac{e^{\eta \bar{e}_i/ (w_i^d+\bar{t}_i)}}{\sum_{k\neq i} e^{\eta\bar{e}_k/ (t_{ik}+w_k^d+\bar{t}_k)}+e^{\eta \bar{e}_i/ (w_i^d+\bar{t}_i)}}.
\end{cases}
\label{Discrete_choice}
\end{align}
For each zone $i$, the arrival rates of drivers ending their trip is $\sum_{k=1}^M \lambda_{ki}$. Since all incoming drivers need to make repositioning decisions, the {\em intended} vehicle rebalance flow from zone $i$ to zone $j$ (or itself) is determined by\footnote{Based on (\ref{Discrete_choice}),   when $t_{ij}$ is large, $\mathbb{P}_{ii}$ can be significantly greater than $\mathbb{P}_{ij}$. In this case, the majority of the drivers will stay in the same zone to seek the next passenger.}
\begin{equation}
\label{flow_intended}
f_{ij}=\mathbb{P}_{ij} \sum_{k=1}^M \lambda_{ki} .
\end{equation}

{
It is important to note that $f_{ij}$ is only the rebalance flow intended by the driver before he/she actually starts repositioning. During  repositioning each vehicle is still available to match  a nearby passenger. Therefore, the vehicle can be intercepted by other zones during the transitional period before it reaches the repositioning destination. To capture this, we denote $\mathcal{P}_{ij}\subset \mathcal{V}$ as the set of zones traversed by the shortest path between zone $i$ and zone $j$ (note that $i,j\in \mathcal{P}_{ij}$), and define $\Pi(k|i\rightarrow j)$ as the probability that a driver is intercepted by zone $k$ when he is on his way to reposition from zone $i$ to zone $j$.  In this case, the resulting rebalancing flow, $\tilde{f}_{ik}$, can be derived as
\begin{equation}
\label{flow_intercepted}
\tilde{f}_{ik}= \sum_{j: k\in \mathcal{P}_{ij}}    \Pi(k|i\rightarrow j) f_{ij},
\end{equation}
The intercepting probability $\Pi(k|i\rightarrow j)$ depends on the demand and supply in each zone along the shortest path $\mathcal{P}_{ij}$, the sequence of traversing these zones, and the duration that a traversing vehicle stays in these zones. To model this intercepting probability, let us first consider a single zone $i$ and denote $d_i$ as the average time spent to traverse zone $i$. In this case, an idle driver that traverses zone $i$ stays in zone $i$ for an average duration of $d_i$. During this period, there is a probability $\sigma_i$ that the idle vehicle is matched to a passenger before it exits zone $i$. To model $\sigma_i$, we propose an M/M/1 queue that captures the stochasticity of passengers and drivers in zone $i$, where passengers are modeled as servers, and drivers are viewed as jobs. This queuing model has been used in other works modeling the ride-sourcing platforms  \cite{banerjee2015pricing}. Based on the  M/M/1 queue, we derive that the waiting time $\tau$ of the drivers in zone $i$ is subject to an exponential distribution \cite{harrison1993response}, therefore we have
\begin{equation}
\label{CDF_exp}
\mathbb{P}(\tau\leq T) =1-e^{-T/w_i^d}.
\end{equation}
Based on (\ref{CDF_exp}), if an idle driver stays for a duration of $d_i$ in zone $i$, then the probability that it is matched to a passenger before it exits zone $i$ is
\begin{equation}
\label{zone_intercept}
\sigma_i=\mathbb{P}(\tau\leq d_i)=1-e^{-d_i/w_i^d}
\end{equation}
The intercepting probability $\Pi(k|i\rightarrow j)$ can be derived based on $\sigma_i$. In particular, if a repositioning vehicle from zone $i$ to zone $j$ is intercepted by zone $k$,  it indicates that this vehicle is not intercepted by any zones on the shortest path $\mathcal{P}_{ij}$ prior to zone $k$. To make this precise, we order the elements of $\mathcal{P}_{ij}$ based on the distance to zone $i$ in an increasing sequence, i.e., the first element $I_1\in\mathcal{P}_{ij}$ is zone $i$ itself, the second element $I_2\in\mathcal{P}_{ij}$ is the closest zone to $i$ on the shortest path $\mathcal{P}_{ij}$, the third element  $I_3\in\mathcal{P}_{ij}$ is the second closest zone to $i$, etc. Let $k$ be the $s$th element, which satisfies $I_s=k$, then for $k\neq j$, the intercepting probability can be derived as:
 \begin{equation}
 \label{overall_intercept}
 \Pi(k|i\rightarrow j)= (1-\sigma_{I_1})\cdot  (1-\sigma_{I_2}) \cdots (1-\sigma_{I_{s-1}}) \sigma_{I_s}.
\end{equation}  
Equation (\ref{overall_intercept}) indicates that if a repositioning vehicle is intercepted by zone $I_s=k$, then it must have passed all the zone $I_1, I_2, \ldots, I_{s-1}$ before it enters $k$. 

\begin{remark}
The above discussion is based on two underlying assumptions. First, we have assumed that the calculation of $\sigma_i$ and $\sigma_j$ are independent. This is a reasonable assumption if we only consider inter-zone matching: the matching radius is smaller than the diameter of each zone, and that each passenger is only matched to a vehicle within the same zone. The case of intra-zone matching \cite{xu2019equilibrium} is outside the scope of this paper. Second, we assume that when drivers make repositioning decisions based on the logit model (\ref{Discrete_choice}), they do not anticipate the possibility of being intercepted by other zones. We believe that this is a realistic assumption since the repositioning of human drivers is often guided by intuition instead of computer algorithms. While they can decide the repositioning strategies based on a rough estimate of prospective earnings,  they lack the computational resources and the real-time information to carry out more delicate calculations.   
\end{remark}

}
Finally, since each vehicle has three operating modes: (a) carrying a passenger, (b) on the way to pick up the passenger, (c) cruising with empty seats to look for the next passenger,  the total number of vehicle hours should be partitioned as
\begin{equation}
\label{conservation}
N=\sum_{i=1}^M \sum_{j=1}^M \lambda_{ij}t_{ij}+\sum_{i=1}^M \sum_{j=1}^M w_i^p \lambda_{ij}  +  \sum_{i=1}^M \sum_{j=1}^M   w_i^d \lambda_{ij},
\end{equation}
where the first term accounts for the operating hours of vehicles with passengers, the second term accounts for the operating hours of vehicles that are on their way to pick up passengers, and the third term accounts for the operating hours of idle vehicles.\footnote{Note that the vehicle hours associated with the repositioning vehicles are incorporated in the third term of  (\ref{conservation}) because all repositioning vehicles are available for matching during the repositioning process.}

\subsection{Traffic Congestion}
\label{section_congestion}
{The ride-sourcing vehicles contribute to the traffic congestion of the city. This affects the traffic speed, which further affects the time it takes to travel from zone $i$ to zone $j$. For this reason, the travel time $t_{ij}$ is an endogenous variable that depends on the spatial distribution of ride-sourcing vehicles. Here we present a traffic congestion model to capture this dependence. 
}

{
First, we note that the congestion externality of ride-sourcing vehicles is spatially asymmetric: in the core area of the city where streets are already heavily congested, introducing more ride-sourcing vehicles will significantly affect the traffic speed, whereas in remote areas where congestion rarely occurs, ride-sourcing vehicles only have a negligible impact on traffic speed. To model this, we consider a ``congestion area" $\mathcal{C}\subset \mathcal{V}$ that consists of all the congested zones in the urban core, and a ``remote area" ($\mathcal{R}=\mathcal{V}\backslash\mathcal{C}$)  that consists of all other zones in the suburb. Let $v_c$ and  $v_r$ be the average traffic speed within the congestion area and remote area, respectively.  Denote $d_{ij}$ as the average shortest-path travel distance from zone $i$ to zone $j$. Since a trip from zone $i$ to zone $j$ may traverse both $\mathcal{C}$ and $\mathcal{R}$, we denote 
 $d_{ij}^{\mathcal{C}}$ and $d_{ij}^{\mathcal{R}}$ as the overlap of $d_{ij}$ with respect to $\mathcal{C}$ and $\mathcal{R}$, respectively.\footnote{If the trip from $i$ to $j$ does not traverse any zone in the congestion area (or remote area), then $d_{ij}^{\mathcal{C}}=0$ (or $d_{ij}^{\mathcal{R}}=0$). In addition, we clearly have $d_{ij}=d_{ij}^{\mathcal{C}} +d_{ij}^{\mathcal{R}}$.} The travel time $t_{ij}$ can be then derived as\footnote{We assume that each side of the equation is properly scaled to have the same unit, wherever necessary. }
\begin{equation}
\label{travel_time}
t_{ij}=d_{ij}^{\mathcal{C}}\dfrac{1}{v_c}+ d_{ij}^{\mathcal{R}}\dfrac{1}{v_r},
\end{equation}
where the first term represents the travel time within the congestion area $\mathcal{C}$ and the second term denotes the travel time within the remote area $\mathcal{R}$. We assume that the average travel speed in the remote area does not depend on the ride-sourcing traffic.  Therefore, in equation (\ref{travel_time}),  $d_{ij}^{\mathcal{C}}$ , $d_{ij}^{\mathcal{R}}$ and $v_r$ are all exogenous variables. On the other hand, the traffic speed in the congestion area should depend on the number of ride-sourcing vehicles. Let $N_{\mathcal{C}}$ denote the total number of ride-sourcing vehicles in  the congestion area $\mathcal{C}$. With slight abuse of notation, the traffic speed in $N_{\mathcal{C}}$ can be then denoted as $v_c(\cdot)$: a decreasing function of $N_{\mathcal{C}}$. Based on equation (\ref{travel_time}),  $t_{ij}$  is also a function of $N_{\mathcal{C}}$.
 }

{ Finally, the number of ride-sourcing vehicles in $\mathcal{C}$ should consist of all vehicles traveling in the congestion area, including occupied vehicles that spend a proportion of trip time in $\mathcal{C}$, idle vehicles cruising in $\mathcal{C}$, and empty vehicles on the way to pick up passengers in $\mathcal{C}$. Similar to (\ref{conservation}), we have that:
\begin{equation}
\label{conservation2}
N_{\mathcal{C}}=\sum_{i=1}^M \sum_{j=1}^M \lambda_{ij}d_{ij}^{\mathcal{C}}\dfrac{1}{v_c}+\sum_{i\in \mathcal{C}}  \sum_{j=1}^M w_i^p \lambda_{ij}  +  \sum_{i\in \mathcal{C}} \sum_{j=1}^M   w_i^d \lambda_{ij},
\end{equation}
where the first term represents all trip times associated with the congestion area,\footnote{We remind that $d_{ij}^{\mathcal{C}}/v_c$ is the time spent in the congestion area for a trip from zone $i$ to zone $j$. If a trip does not pass through the congestion area, then we have $d_{ij}^{\mathcal{C}}=0$.} the second term represents all vehicles on the way to pick up the passenger in $\mathcal{C}$, and the third term represents all idle cruising vehicles in $\mathcal{C}$. Each term in (\ref{conservation2}) reflects a proportion of the corresponding terms in (\ref{conservation}), which is only associated with the congestion area. 
}

{
\begin{remark}
Ideally, traffic congestion can be modeled at the zonal level: each zone in the city has its own traffic speed, which depends on the total number of ride-sourcing vehicles in this zone. Our model does not capture congestion at this level of granularity because we primarily focus on congestion charges that differentiate the urban core and the suburb. Therefore, a congestion model at the urban/suburb level already suffices for our purpose. Furthermore, a zonal-level congestion model requires significantly more data (background traffic, road capacity, etc), which is unavailable to us. Therefore,  we leave this more detailed congestion model as future work. 
\end{remark}
}

\subsection{Network Flow Balance}
In the long run, the inflow and outflow of vehicles in each zone should be  balanced. This leads to the following constraints,
\begin{equation}
\label{rebalance}
\sum_{j=1}^M  (\lambda_{ji}+\tilde{f}_{ji})=\sum_{j=1}^M  (\lambda_{ij}+\tilde{f}_{ij}), \quad  \forall i\in \mathcal{V}.
\end{equation}
The left-hand side of (\ref{rebalance}) denotes the total inflow to zone $i$ consisting of vehicles with and without passengers, and the right-hand side of (\ref{rebalance}) denotes the total outflow from zone $i$  consisting of vehicles with and without passengers. These flows have to be balanced in  equilibrium.

\subsection{Platform Decision Model}
Consider a ride-sourcing platform that determines the ride fare $r_i$ and driver wage $q$ (or equivalently, a commission rate\footnote{In practice, ride-sourcing platforms often determine a commission rate which indirectly leads to a driver wage $q$. However, we note that optimizing over $q$ is equivalent to optimizing over the commission rate as there is a one-to-one mapping between driver wage $q$ and the commission rate.}) to maximize its profit. These decisions are subject to the passenger demand model, driver supply model, driver repositioning model, and the network flow balancing constraints. The optimal spatial pricing problem can be formulated as follows:
\begin{equation}
\label{optimalpricing_trip}
 \hspace{-5cm} \mathop {\max }\limits_{{\bf r}, q} \quad \sum_{i=1}^M \sum_{j=1}^M r_it_{ij}\lambda_{ij}-N_0F_d(q)q
\end{equation}
\begin{subnumcases}{\label{constraint_optimapricing}}
{
\lambda_{ij}  = {\lambda^0_{ij}}{F_p}\left(\alpha {w^p_i}(N^I_i)  + \beta t_{ij}+ r_it_{ij}\right)} \label{demand_constraint}\\
w_i^p(N_i^I)\leq w_{max} \label{upper_waiting}\\
N_i^I=w_i^d\sum_{j=1}^M \lambda_{ij} \label{demand_constraint2}\\
 {N_0}{F_d}\left(q\right)  =\sum_{i=1}^M \sum_{j=1}^M \lambda_{ij}t_{ij}+\sum_{i=1}^M \sum_{j=1}^M w_i^p \lambda_{ij}  +  \sum_{i=1}^M \sum_{j=1}^M   w_i^d \lambda_{ij} \label{suply_conservation_const}  \\
{N_{\mathcal{C}}=\sum_{i=1}^M \sum_{j=1}^M \lambda_{ij}d_{ij}^{\mathcal{C}}\dfrac{1}{v_c}+\sum_{i\in \mathcal{C}}  \sum_{j=1}^M w_i^p \lambda_{ij}  +  \sum_{i\in \mathcal{C}} \sum_{j=1}^M   w_i^d \lambda_{ij}} \label{_supply_conservation2} \\
\sum_{j=1}^M  (\lambda_{ji}+\tilde{f}_{ji})=\sum_{j=1}^M  (\lambda_{ij}+\tilde{f}_{ij}), \quad   i\in \mathcal{V}. \label{reballance_const} \\
\tilde{f}_{ij}= \sum_{k: j\in \mathcal{P}_{ik}}   \Pi(j|i\rightarrow k) \mathbb{P}_{ik} \sum_{m=1}^M \lambda_{mi}  \label{definitionofflow}
\end{subnumcases}
where ${\bf r}=(r_1,r_2, \ldots, r_M)$, $t_{ij}$ is defined by (\ref{travel_time}), $\Pi(j|i\rightarrow k)$ is given by (\ref{overall_intercept}), and $\mathbb{P}_{ik}$ is determined by the logit model (\ref{Discrete_choice}). The objective function (\ref{optimalpricing_trip}) defines the platform profit as the total revenue $\sum_{i=1}^M \sum_{j=1}^M r_it_{ij}\lambda_{ij}$ minus the total driver payment  $Nq=N_0F_d(q)q$. Constraints (\ref{demand_constraint})-(\ref{demand_constraint2}) specify the passenger demand. 
Constraint (\ref{upper_waiting}) requires  the waiting time in each zone to be smaller than an upper bound.\footnote{In certain zones of the city,  passenger demand is very low, and profit-maximizing decision is to offer no service to this zone. We impose an upper bound on the waiting time to avoid trivial solution arising from this scenario.}
Constraint (\ref{suply_conservation_const}) combines the driver supply model (\ref{driversupply}) and the vehicle hour conservation constraint (\ref{conservation}). Constraint (\ref{definitionofflow}) is obtained by plugging  (\ref{flow_intended}) into (\ref{flow_intercepted}). The optimal spatial pricing problem is a non-convex program. After the platform specifies the per-time rate $r_i$ and the driver payment $q$,  the constraint (\ref{constraint_optimapricing}) defines  $M$ endogenous variables $w_i^d$, $\forall i\in \mathcal{V}$. We can substitute (\ref{demand_constraint})-(\ref{demand_constraint2}) into (\ref{suply_conservation_const}) and substitute (\ref{definitionofflow}) into (\ref{reballance_const}). Since (\ref{reballance_const}) only has $M-1$ independent constraints,\footnote{If the balancing constraint holds for $M-1$ zones,  it automatically holds for the other zone. This is in same spirit as the Walras's law \cite{mas1995microeconomic}.}  overall (\ref{suply_conservation_const}) and (\ref{reballance_const}) constitute  $M$ independent constraints. {We further note that given the $M$ endogenous variables $w_i^d$,  $N_\mathcal{C}$ can be derived as a function of these variables through constraint (\ref{_supply_conservation2}). This further determines $t_{ij}$ as a function of $w_i^d$. Therefore, after plugging $t_{ij}$ into (\ref{suply_conservation_const}) and (\ref{reballance_const}), overall,  (\ref{suply_conservation_const}) and (\ref{reballance_const}) constitute  $M$ independent constraints for $M$ endogenous variables when $r_i$ and $q$ are given. }

\begin{remark}
Constraint (\ref{upper_waiting}) imposes an upper bound on the pickup time in each zone. This requires the platform to recruit enough drivers to maintain its service quality in all zones, which implicitly places a lower bound on the driver wage $q$. However, we comment that the imposition of (\ref{upper_waiting}) is to avoid the trivial case where the platform finds it more profitable to offer no service to some low-demand areas.  For this reason, the value of $w_{max}$ is large (e.g., 10 min) relative to the average pickup time of other zones (e.g., 5 min), thus for most zones the upper bound on the pickup time (\ref{upper_waiting}) is inactive. This implies that the imposition of (\ref{upper_waiting}) will not significantly affect the well-posedness of the optimization problem (\ref{optimalpricing_trip}).  
\end{remark}

\section{The Solution Algorithm}
 This section proposes an algorithm to approximately compute the optimal solutions to (\ref{optimalpricing_trip})  and to establish an upper bound to evaluate the performance of the proposed algorithm. Both the algorithm and the upper bound are validated through realistic synthetic ride-sourcing data for San Francisco. 

\subsection{Constraint Relaxation and Upper Bound}
\begin{figure*}%
\begin{minipage}[b]{0.5\linewidth}
\centering
\includegraphics[width = 1\linewidth]{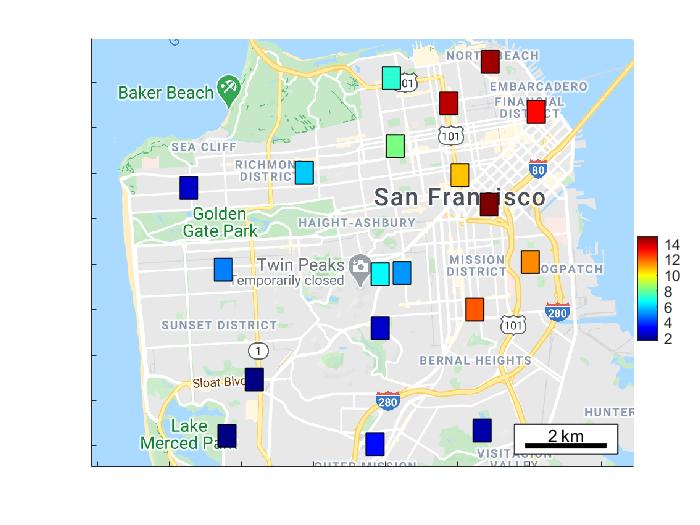}
\caption{Total inflow per minute to each zip code zone.}
\label{inflow_25}
\end{minipage}
\begin{minipage}[b]{0.5\linewidth}
\centering
\includegraphics[width = 1\linewidth]{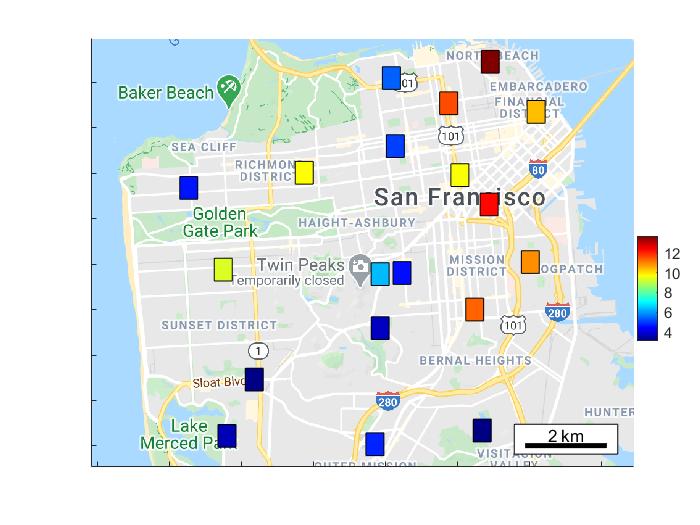}
\caption{Total outflow per minute from each zip code zone.}
\label{outflow_25}
\end{minipage}
\end{figure*}

The primary constraints in  (\ref{optimalpricing_trip}) are (\ref{suply_conservation_const}) and (\ref{reballance_const}). The vehicle hour balancing constraint (\ref{suply_conservation_const}) directly relates to driver payment, which has a significant impact on the platform profit. On the other hand, the network flow balance constraint (\ref{reballance_const}) does not directly affect the profit platform, but only indirectly affects spatial prices by dictating network balance. It is important to note that in a ride-sourcing market, a popular origin is typically also a popular destination (see for example, Figure \ref{inflow_25} and Figure \ref{outflow_25}). For this reason, the incoming passenger flow $\sum_{j=1}^M \lambda_{ji}$  and outgoing passenger flow $\sum_{j=1}^M \lambda_{ij}$  in each zone $i$ are naturally close even without vehicle repositioning. Therefore, we anticipate that relaxing (\ref{reballance_const}) will not significantly affect the platform's optimal profit. The subsequent discussion is built upon this intuition. 

By dropping the network flow balance constraint (\ref{reballance_const}), the optimal spatial pricing problem (\ref{optimalpricing_trip}) can be relaxed to:
\begin{equation}
\label{optimalpricing_relax}
 \hspace{-4cm} \mathop {\max }\limits_{{\bf r}, {\bf N^I}, q} \quad \sum_{i=1}^M \sum_{j=1}^M r_it_{ij}\lambda_{ij}-N_0F_d(q)q
\end{equation}
\begin{subnumcases}{\label{constraint_optimapricing_relax}}
{
\lambda_{ij}  = {\lambda^0_{ij}}{F_p}\left(\alpha {w^p_i}(N^I_i)  + \beta t_{ij}+ r_it_{ij}\right)} \\
w_i^p(N_i^I)\leq w_{max} \\
 {N_0}{F_d}\left(q\right)  =\sum_{i=1}^M \sum_{j=1}^M \lambda_{ij}t_{ij}+\sum_{i=1}^M \sum_{j=1}^M w_i^p \lambda_{ij}  +  \sum_{i=1}^M N_i^I, \label{suply_conservation_relax}  \\
 {N_{\mathcal{C}}=\sum_{i=1}^M \sum_{j=1}^M \lambda_{ij}d_{ij}^{\mathcal{C}}\dfrac{1}{v_c}+\sum_{i\in \mathcal{C}}  \sum_{j=1}^M w_i^p \lambda_{ij}  +  \sum_{i\in \mathcal{C}} N_i^I} \label{supply_conservation2_relax} 
\end{subnumcases}
where ${\bf N^I}=(N_1^I, N_2^I, \ldots, N_M^I)$,  (\ref{suply_conservation_relax}) is derived by substituting (\ref{demand_constraint2}) into (\ref{suply_conservation_const}), and   (\ref{definitionofflow}) is removed since  $\tilde{f}_{ij}$ only appears in (\ref{reballance_const}). { In the relaxed problem, since multiple constraints (\ref{reballance_const}) are removed, it releases more freedom for the platform to make decisions. In this case, }the platform has the freedom to place idle vehicles\footnote{It is easy to see that optimizing over ${\bf N^I}=(N_1^I, \ldots, N_M^I)$ is equivalent to optimizing over $(w_1^d, \ldots, w_M^d)$.}(note that this is not the case for the original problem) and optimize over ${\bf N^I}$. The optimal value of the relaxed problem (\ref{optimalpricing_relax}) provides an upper bound to the optimal spatial pricing problem (\ref{optimalpricing_trip}).

\subsection{The Solution Algorithm}
{
The relaxed spatial pricing problem (\ref{optimalpricing_relax}) has a special structure that is amenable to numerical computation: when $N_{\mathcal{C}}$ is given, the traffic speed $v_c$ and the travel times $t_{ij}$ are uniquely determined. In this case, the decision variables for each zone $i$, i.e., $r_i$ and $N_i^I$, are separable both in the objective function (\ref{optimalpricing_relax}) and in the constraints (\ref{suply_conservation_relax})-(\ref{supply_conservation2_relax}). This indicates that the overall problem can be solved by enumerating over $N_{\mathcal{C}}$ and then performing a dual decomposition under each $N_{\mathcal{C}}$. In particular, we first conduct a grid search over the scalar decision variable $N_{\mathcal{C}}$. For each $N_{\mathcal{C}}$, we dualize (\ref{suply_conservation_relax}) and (\ref{supply_conservation2_relax}) to decompose the overall problem into a few optimization problems of much smaller size  (e.g., two-dimensional optimization), which can be solved individually and in parallel. In each iteration of the dual-decomposition algorithm, we collect solutions from the decomposed problems to update the dual variable, repeat this calculation until the stopping criterion is satisfied. This entire process is then replicated for a different $N_{\mathcal{C}}$ until we find the optimal $N_{\mathcal{C}}$ that leads to the highest platform profit. Since the dual-decomposition step may introduce a non-zero duality gap, the derived solution can be either the exact solution to  (\ref{optimalpricing_relax}) or an upper bound to (\ref{optimalpricing_relax}). In either case, the algorithm provides an upper bound for the original optimal pricing problem  (\ref{optimalpricing_trip}). Details of the algorithm are summarized in Algorithm~\ref{algorithm1}.
}

\begin{algorithm}
\caption{The solution algorithm for the optimal spatial pricing problem (\ref{optimalpricing_trip})} \label{algorithm1}
\begin{algorithmic}[1]
\REQUIRE Initial guess of primal variables ${\bf \bar{r}}$, ${\bf \bar{N}^I}$, $\bar{q}$, and dual variables $\delta$, $\kappa$.
\STATE Setup stopping criterion: constraints (\ref{suply_conservation_relax})-(\ref{supply_conservation2_relax})  satisfied or maximum iteration reached. 
\FOR {each iteration} 
\STATE Choose a candidate $N_{\mathcal{C}}$,
\STATE Calculate $v_c$ and $t_{ij}$ based on $N_{\mathcal{C}}$,
\WHILE {stopping criterion not satisfied}
\STATE {If $i\in \mathcal{C}$, solve the decomposed Lagrangian for  zone $i$:
\begin{equation}
\label{primaldual4}
\hspace{-0.5cm}(r_i, N_i^I)=\argmax_{r_i, N_i^I}  \sum_{j=1}^M r_it_{ij}\lambda_{ij}-\delta \left(  \sum_{j=1}^M \lambda_{ij}t_{ij}+ \sum_{j=1}^M w_i^p \lambda_{ij}  +  N_i^I \right) -\kappa \left(  \sum_{j=1}^M \lambda_{ij}d_{ij}^{\mathcal{C}}\dfrac{1}{v_c}+ \sum_{j=1}^M w_i^p \lambda_{ij}  +  N_i^I \right) 
\end{equation} 
\begin{subnumcases}{\hspace{-7cm}}{\label{constraint_optimapricing_relax_2_case1}}
\lambda_{ij}  = {\lambda^0_{ij}}{F_p}\left(\alpha {w^p_i}(N^I_i)  +\beta t_{ij}+ r_it_{ij}\right) \\
w_i^p(N_i^I)\leq w_{max} 
\end{subnumcases}
}

{
If $i\in \mathcal{R}$, solve the decomposed Lagrangian for  zone $i$:
\begin{equation}
\label{primaldual1}
\hspace{1cm}(r_i, N_i^I)=\argmax_{r_i, N_i^I}  \sum_{j=1}^M r_it_{ij}\lambda_{ij}-\delta \left(  \sum_{j=1}^M \lambda_{ij}t_{ij}+ \sum_{j=1}^M w_i^p \lambda_{ij}  +  N_i^I \right) - \kappa  \sum_{j=1}^M \lambda_{ij}d_{ij}^{\mathcal{C}}\dfrac{1}{v_c}
\end{equation}
\begin{subnumcases}{\hspace{-3cm}}{\label{constraint_optimapricing_relax_2_case2}}
\lambda_{ij}  = {\lambda^0_{ij}}{F_p}\left(\alpha {w^p_i}(N^I_i)  +\beta t_{ij}+ r_it_{ij}\right) \\
w_i^p(N_i^I)\leq w_{max} 
\end{subnumcases}
}
\STATE Solve the decomposed Lagrangian over driver supply,
\begin{equation}
\label{primaldual2}
\bar{q}=\argmax_{q>0} \delta N_0F_d(q)-N_0F_d(q)q,
\end{equation}
\STATE Update the dual variable
{
\begin{align}
\label{primaldual3} 
\begin{cases}
\delta =\delta-\gamma_1  \left({N_0}{F_d}\left(q\right) -\sum_{i=1}^M \sum_{j=1}^M \lambda_{ij}t_{ij}-\sum_{i=1}^M \sum_{j=1}^M w_i^p \lambda_{ij}  -  \sum_{i=1}^M N_i^I\right) \\
\kappa =\kappa-\gamma_2  \left( N_{\mathcal{C}}-\sum_{i=1}^M \sum_{j=1}^M \lambda_{ij}d_{ij}^{\mathcal{C}}\dfrac{1}{v_c}-\sum_{i\in \mathcal{C}}  \sum_{j=1}^M w_i^p \lambda_{ij}  -  \sum_{i\in \mathcal{C}} N_i^I\right) \\
\end{cases}
\end{align}
}
\ENDWHILE 
\ENDFOR
\STATE Obtain the solutions $({\bf \bar{r}}, {\bf \bar{N}^I}, \bar{q})$ and the corresponding platform profit $\bar{R}$. 
\STATE Use $({\bf \bar{r}}, {\bf \bar{N}^I}, \bar{q})$ as the initial guess to solve  (\ref{optimalpricing_trip}) using  interior-point algorithm, and obtain the solutions $({\bf {r}}, {\bf {N}^I}, {q})$ and the corresponding platform profit ${R}$. 
\ENSURE the approximate solution $({\bf {r}}, {\bf {N}^I}, {q})$, the corresponding platform profit $R$ and its upper bound $\bar{R}$.
\end{algorithmic}
\end{algorithm}

In Algorithm \ref{algorithm1}, Step 5-Step 9 executes a standard dual-decomposition for the relaxed problem (\ref{optimalpricing_relax}). In each iteration, the Lagrangian of (\ref{optimalpricing_relax}) is decomposed into decoupled problems (\ref{primaldual4}), (\ref{primaldual1}) and (\ref{primaldual2}), which are small-scale non-convex programs (1D or 2D) that can be efficiently solved to global optimality by brute-force computation. Depending on the duality gap of (\ref{optimalpricing_relax}), the termination conditions of  the dual-decomposition iteration have two scenarios: (a) there is no duality gap of (\ref{optimalpricing_relax}), and the algorithm converges to a solution that satisfies  constraints (\ref{suply_conservation_relax})-(\ref{supply_conservation2_relax}); (b) there is a  duality gap of (\ref{optimalpricing_relax}) and the algorithm terminates to an infeasible point when the algorithm reaches the maximum number of iterations.\footnote{The platform can terminate the while loop at any time and implement the resulting solution $\bar{r}$ and $\bar{q}$, with $\bar{R}$ as an upper bound on its profit.} In either case, the algorithm output $\bar{R}$ provides an upper bound on the optimal value of (\ref{optimalpricing_trip}).

\begin{proposition}
\label{proposition2}
After Algorithm 1 terminates, (a) if constraints (\ref{suply_conservation_relax})-(\ref{supply_conservation2_relax}) are satisfied at $({\bf \bar{r}}, {\bf \bar{N}^I}, \bar{q})$,  then $({\bf \bar{r}}, {\bf \bar{N}^I}, \bar{q})$ is the globally optimal solution to  (\ref{optimalpricing_relax}), and $\bar{R}$ is an upper bound on the optimal value of (\ref{optimalpricing_trip});
 (b) if  constraints (\ref{suply_conservation_relax})-(\ref{supply_conservation2_relax}) are not satisfied at $({\bf \bar{r}}, {\bf \bar{N}^I}, \bar{q})$, then $({\bf \bar{r}}, {\bf \bar{N}^I}, \bar{q})$ is not the globally optimal solution to  (\ref{optimalpricing_relax}), but $\bar{R}$ is still an upper bound on the optimal value of (\ref{optimalpricing_trip}).
\end{proposition}

The proof of Proposition \ref{proposition2} can be found in \cite[p.385, p.563]{bertsekas1997nonlinear}, and is therefore omitted.  We note that in deriving the optimal spatial prices, the role of the relaxed problem (\ref{optimalpricing_relax}) is crucial: its optimizer provides a good initial guess for solving the optimal spatial pricing problem (\ref{optimalpricing_trip}), while its optimal value provides an upper bound on the optimality loss.

\begin{remark}
{
Although the rebalancing constraint (\ref{reballance_const}) is relaxed when we compute the upper bond, this relaxation is only performed to provide an initial guess and an upper bound. Ultimately, Algorithm \ref{algorithm1} provides a final solution that satisfies all constraints in (\ref{constraint_optimapricing}), including the rebalancing constraint (\ref{reballance_const}). We emphasize that although the rebalancing constraint does not significantly affect the optimal value of the spatial pricing problem, it does affect the optimal solution when we consider the impact of a congestion charge: if (\ref{reballance_const}) is completely ignored in the model, then very different conclusions will be drawn in Section~\ref{sec_congestion_charge}. 
}
\end{remark}

\subsection{Case Studies}
\label{numericalexample}
\begin{figure}
\centering
\includegraphics[width = 0.6\linewidth]{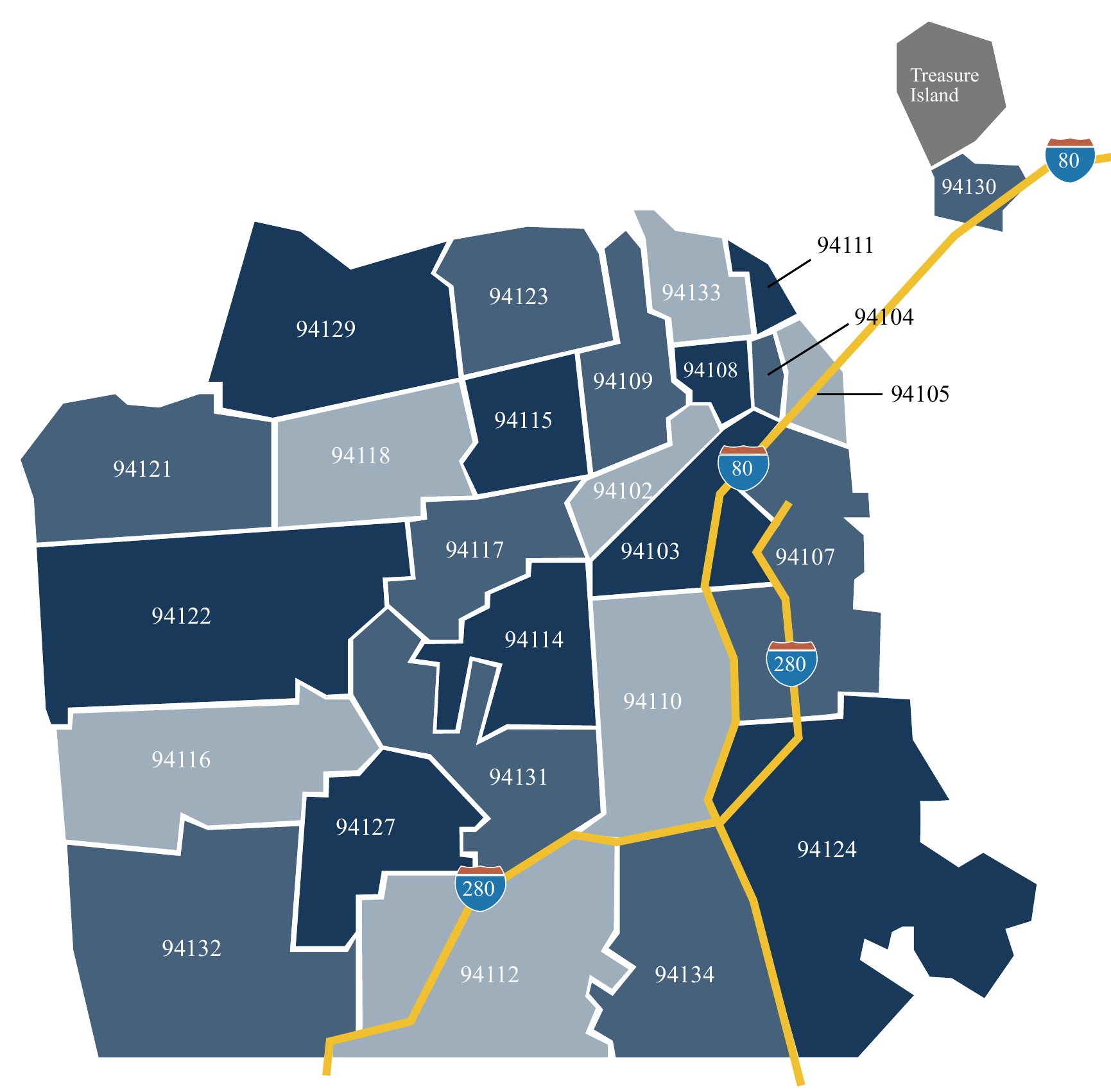}
\caption{Zip code zones of San Francisco County (figure courtesy: usmapguide.com).}
\label{SFzip}
\end{figure}

To test the proposed algorithm, we  conduct a numerical study using realistic synthetic data for San Francisco. The data consists of the origin and destination of each ride-sourcing trip at the zip-code granularity, which is synthesized based on the total pickup and dropoff in each zone \cite{SFCTA2016} combined with a  choice model calibrated by  survey data. The zip code zones of San Francisco are shown in Figure \ref{SFzip}, and the total passenger  inflow (e.g.,$\sum_{j=1}^M \lambda_{ji}$) and outflow  (e.g., $\sum_{j=1}^M \lambda_{ij}$) of each zone are shown  in Figure \ref{inflow_25} and Figure \ref{outflow_25}, respectively. Based on the data, we remove zip code zones 94127, 94129, 94130 and 94134  from our analysis since they all have negligible trip volumes. We also aggregate zip code zones 94111, 94104 and 94105 into a single zone, and aggregate 94133 and 94108 into a single zone, since each of these individual zones is very small.  

%

Passenger choice among different transport modes is captured by a logit model, described as
\begin{equation}
\label{logit_demand}
	\lambda_{ij} =\lambda^0_{ij} \frac{e^{-\epsilon c_{ij}}}{e^{-\epsilon c_{ij}}+e^{-\epsilon c^0_{ij}}},
\end{equation}
where $\epsilon$ and $c^0_{ij}$ are model parameters. For driver supply, the long-term driver subscription decisions is  captured by
\begin{equation}
\label{logit_supply}
	N =N_0 \frac{e^{\sigma q}}{e^{\sigma q}+e^{\sigma q_0}},
\end{equation}
where $\sigma$ and $q_0$ are  model  parameters.

Assume that the passenger waiting time  $w^p_i$ follows the ``square root law'' established in \cite{arnott1996taxi} and \cite{li2019regulating}:
\begin{equation}
	w^p_i(N_i^I) = \frac{L}{\sqrt {N_i^I } },
\label{pickuptime_func}
\end{equation}
where the constant $L$ depends on the size of the zone and demand/supply distribution. The square root law implies that the average passenger waiting time in zone $i$ is inversely proportional  to the square root of the number  of idle vehicles in zone $i$. The intuition behind (\ref{pickuptime_func}) is that if all idle vehicles are uniformly distributed in zone $i$,  then the distance between two neighboring idle vehicles is inversely proportional to the square root of the total number of idle vehicles in zone $i$. This distance is also proportional to the distance between a passenger and their closest driver, which determines the passenger waiting time.\footnote{{The total waiting time for the passenger is the sum of ride confirmation time (from ride request to confirmation) and the pickup time (from ride confirmation to pickup). { Typically, the ride confirmation time (e.g., 30 sec) is much shorter than the pickup time (e.g., 5 min) when there is no  maximum matching radius and an idle vehicle is always available somewhere.} Therefore, in this numerical study we ignore the ride confirmation time and focus on the pickup time.}} Detailed justification of the square root law can be found in \cite{li2019regulating}.

{We further assume that the average value of the inverse of the average traffic speed\footnote{We choose the inverse of traffic speed (instead of speed) because it is easier to calculate the average trip time $t_{ij}$ in (\ref{travel_time}), given that the actual trip time is random across different trips.} in the congestion area is a linear function of $N_{\mathcal{C}}$, i.e, $1/v_c=1/v_c^0+\rho N_{\mathcal{C}}$, where $v_c^0$ and $\rho$ are model parameters. We emphasize that the choice of this model is only for illustration purpose, which does not affect our methodology or conclusion.  
}

{
In summary, the model parameters are
\[{\Theta}=\{\lambda^0_{ij}, N_0,  L,  \alpha, \beta, \epsilon, \sigma, \eta, c^0_{ij},  q_0, d_{ij}, v_c^0, \rho, v_r\}.
\]
The values of these model parameters are set based on  data for San Francisco city. In particular,  $\lambda^0_{ij}$ is set to satisfy $0.15\lambda^0_{ij}=\lambda_{ij}$ ($\lambda_{ij}$ is observed from data) so that $15\%$ of the potential passengers take ride-sourcing trips \cite{castiglione2016tncs}. The travel distance $d_{ij}$ and traffic speed $v_r$ are obtained from Google map estimates. The traffic speed model in the congestion area is calibrated so that  1000 ride-sourcing vehicles in the congestion area contributes to $10\%$ reduction of traffic speed ({with respect to the traffic speed without ridesourcing vehicles}).  The travel cost of the alternative transport mode $c_{ij}^0$ (e.g., taxi, buses) is assumed to be proportional to $d_{ij}$\footnote{{By considering $c_{ij}^0$ as a function of $d_{ij}$, it is inherently assumed that alternative transport options are not subject to traffic congestion. This assumption applies to many transport modes, such as buses with dedicated lanes, subway, walking, biking, etc. These modes altogether  serve a lion's share of the overall mobility needs. Therefore, it is reasonable to assume that ride-hailing passengers are primarily choosing between TNC and these congestion-independent transport options}.}. To account for unavailability of adequate public transit in remote areas such as zip code zone 94121, 94118, 94125, 94110, 94132, 94112, 94124, we assume that the per-distance cost of alternative transport modes in these areas is 50\% higher than in the urban core.  The  rest of the model parameters are set as
}

\[ \, N_0=10000,  \, L = 43, \, \alpha=3, \beta=1, \epsilon=0.12, \sigma=0.17, \eta=0.1, q_0=\$29.34/\text{hour}, \]
\[w_{max}=10 \text{min}, \, v_c^0=15\text{mph}, \, \rho=10\%\dfrac{1}{v_c^0}\dfrac{1}{1000}, \, v_r=20\text{mph}.\]

These parameter values are adjusted so that the corresponding optimal solution is close to the real-world data of San Francisco (i.e., trip volumes, average trip fare, number of drivers, driver wage, etc). { In particular, at the solution to Algorithm 1, passenger arrival rate is 156/min, total number of drivers is 3687, average trip fare is \$17.2/trip,  driver wage is \$26.2/hour. The passenger arrival rate and driver supply are consistent with \cite{castiglione2016tncs} during working hours (e.g., 4PM-5PM) on a typical weekday. The driver wage is close to the estimates\footnote{It is estimated that the before expense earning of ride-sourcing drivers in NYC is \$25.76/hour after the minimum wage regulation \cite{castiglione2016tncs}. The minimum wage in SF is \$0.59 higher than that of NYC, so we estimate that driver wage in SF is $26.35$/hour. This is close to the result of our numerical study.} (before expense) of \cite{parrott2018earning}. The average trip fare is close to the fare estimates \cite{lyftprice} for a 2.6 mile trip \cite{castiglione2016tncs}. 
 }

 \begin{figure*}%
\begin{minipage}[b]{0.5\linewidth}
\centering
\includegraphics[width = 1\linewidth]{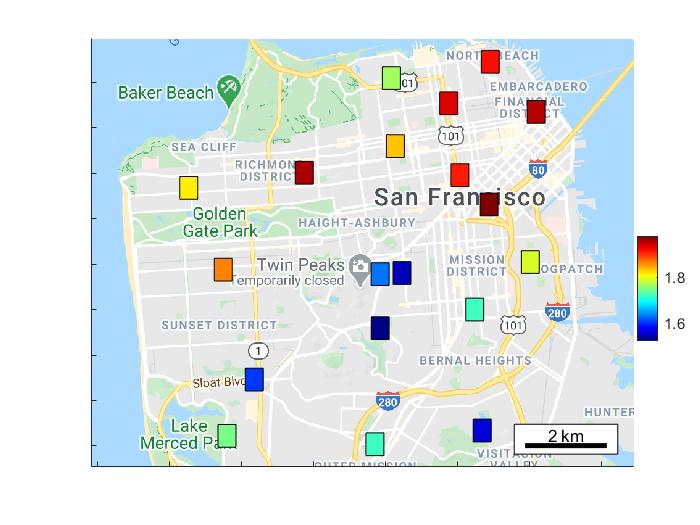}
\caption{Ride price (\$/min) of ride-sourcing trips.}
\label{price_unreg}
\end{minipage}
\begin{minipage}[b]{0.5\linewidth}
\centering
\includegraphics[width = 1\linewidth]{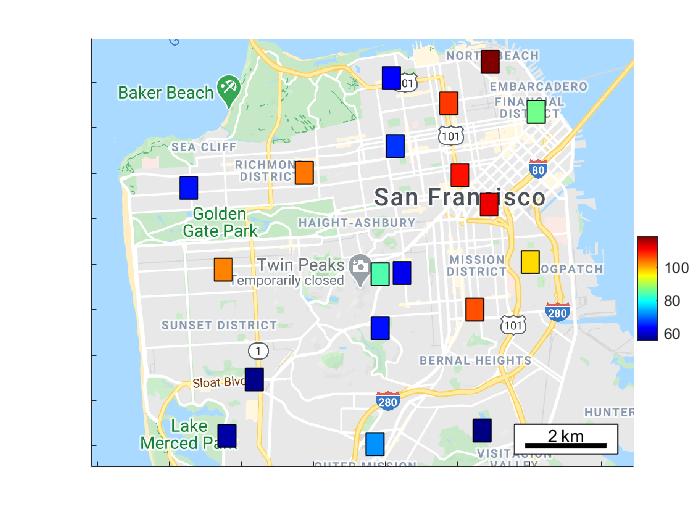}
\caption{Number of idle vehicles in each zone.}
\label{NI_unreg}
\end{minipage}
\end{figure*}

 \begin{figure*}%
\begin{minipage}[b]{0.5\linewidth}
\centering
\includegraphics[width = 1\linewidth]{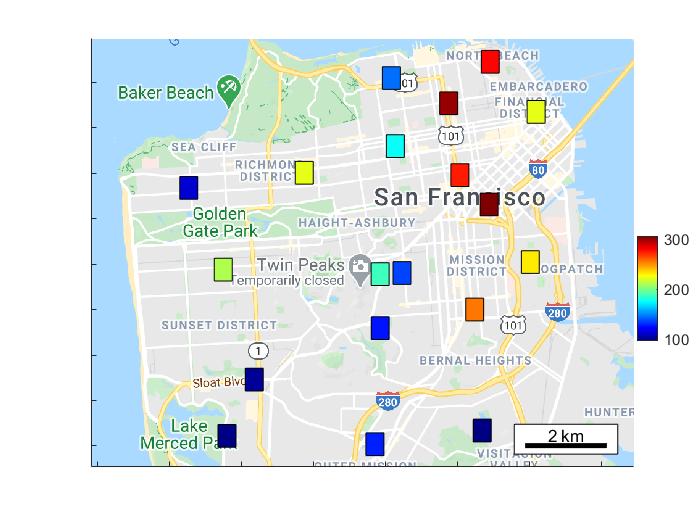}
\caption{Total number of vehicles in each zone.}
\label{N_total_mz}
\end{minipage}
\begin{minipage}[b]{0.5\linewidth}
\centering
\includegraphics[width = 1\linewidth]{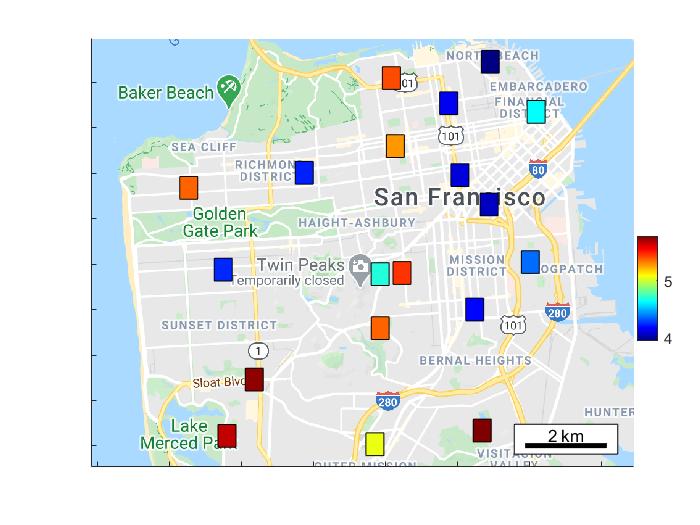}
\caption{Passenger waiting time (min) in each zone.}
\label{waiting_unreg}
\end{minipage}
\end{figure*}

{
{ We execute Algorithm \ref{algorithm1} using the aforementioned parameter values in Matlab on a Dell desktop machine with 8-core i7-9700 CPU (up to 3.00 GHz), which takes around 196 minutes. }  At the solution $({\bf \bar{r}}, {\bf \bar{N}^I}, \bar{q})$,  constraints (\ref{suply_conservation_relax})-(\ref{supply_conservation2_relax}) are satisfied. Based on Proposition \ref{proposition2}, $({\bf \bar{r}}, {\bf \bar{N}^I}, \bar{q})$ is the globally optimal solution to  (\ref{optimalpricing_relax}). The optimal value (platform profit) to (\ref{optimalpricing_relax}) is \$66,737/hour, while the optimal value of (\ref{optimalpricing_trip}), obtained using Algorithm \ref{algorithm1},  is \$64,749/hour. This indicates that the performance loss of the derived solution is at most $3.1\%$ compared with the globally optimal solution to (\ref{optimalpricing_trip}). 
 The ride fare, passenger waiting times, and the distribution of idle ride-sourcing vehicles are shown in Figure \ref{price_unreg}-Figure \ref{waiting_unreg}. The spatial distribution of the  ride-sourcing market is  unbalanced: the urban core with its high demand zones (zip code zones 94102, 94103, 94104, 94105, 94107, 94108, 94109, 94111, 94133) faces high trip fares (Figure \ref{price_unreg}) and low waiting times (Figure \ref{waiting_unreg}),  while  suffering greater  congestion (13.5mph in $\mathcal{C}$ v.s. 20mph in $\mathcal{R}$) from attracting idle ride-sourcing vehicles: based on Figure \ref{NI_unreg}, more than $44\%$ of the idle vehicles are concentrated in the northeast corner of the city, which only accounts for $18\%$ of the total area. 
 }

 To investigate how model parameters affect the performance of Algorithm 1, we conduct a sensitivity analysis by perturbing the parameter values of $\alpha$, $\epsilon$, $\sigma$, $\eta$, $N_0$ and $\lambda_0$ in both directions by 30\%. Each of the selected parameters is perturbed one by one  while other parameters are fixed at their nominal values. Since $\lambda_0$ is a matrix, we perturb all  entries of $\lambda_0$ proportionately. The upper bound on the performance loss is defined as $\dfrac{\bar{R}-R}{\bar{R}}$, where $R$ represents the (potentially suboptimal) platform profits of (\ref{optimalpricing_trip})  and $\bar{R}$ represents the platform profits of (\ref{optimalpricing_relax}), both evaluated at the solutions to Algorithm \ref{algorithm1}. The numerical results are shown in Figure \ref{sensitivity1_mz}-Figure \ref{sensitivity6_mz} (Appendix A). It is clear that the upper bound on the performance loss is consistently small for a large range of model parameters. 
 

\section{Impacts of Congestion Charges}
\label{sec_congestion_charge}
Our numerical study so far suggests that the spatial distribution of the ride-sourcing vehicles is highly unbalanced. Idle drivers prefer to cruise in the urban core to pick up the next passenger within a shorter time, which leads to  the flooding of ride-sourcing vehicles in the urban core and limited availability of ride-sourcing services in remote areas.  This section explores regulatory policies that limit congestion in the urban core and benefit travelers in remote areas. 

\subsection{Optimal Pricing under a Congestion Charge}
{
 To limit traffic in the congestion area, we consider three forms of congestion charge: (a) a one-directional charge on all vehicles that enter the congestion zone; (b) a bi-directional charge on all vehicles that either enter or exit the congestion zone; (c) a congestion charge on all ride-sourcing trips. The first two congestion charges are cordon-based, which impose a per-crossing fare on all vehicles that cross the cordon regardless of whether there is a passenger on-board. We assume that when an occupied vehicle passes the cordon, the congestion charge is paid by the passenger,\footnote{{The platform typically passes the congestion charge to passengers by adding an extra fee on top of the total trip fare. This is consistent with our formulation. We would like to emphasize that in the long term, the platform will need to adjust the ride fare and ultimately the tax burden is shared between the platform and the passengers depending on demand elasticity.} } whereas for idle vehicles  that pass the cordon, the extra charge is paid by the driver.\footnote{{We envision that when a congestion charge is imposed on idle vehicles, the ride-sourcing platform will also pass the extra charge to drivers, just like it did to passengers.} }  On the other hand, the third congestion charge is trip-based, which  imposes a charge only on passengers and does not penalize idle vehicles.  We define $\bar{p}_{ij}$ as the congestion charge on vehicles that travel from $i$ to $j$, and denote $\mathbbm{1}^p_{ij}$ and $\mathbbm{1}^d_{ij}$ as an indicator function to represent whether the charge is imposed on passengers or drivers, respectively. The optimal spatial pricing problem can be written as
\begin{equation}
\label{optimalpricing_trip_charge_multizone}
 \hspace{-2cm} \mathop {\max }\limits_{{\bf r}, q} \quad \sum_{i=1}^M \sum_{j=1}^M r_it_{ij}\lambda_{ij}-N_0F_d(q)q -\sum_{i=1}^M \sum_{j=1}^M\tilde{f}_{ij} \mathbbm{1}^d_{ij} \bar{p}_{ij}
\end{equation}
\begin{subnumcases}{\label{constraint_optimapricing_charge_multi}}
\lambda_{ij}  = {\lambda^0_{ij}}{F_p}\left(\alpha {w^p_i}(N^I_i)  + r_it_{ij}+ \beta t_{ij}+\bar{p}_{ij} \mathbbm{1}^p_{ij}\right) \label{demand_constraint_charge_multi}\\
w_i^p(N_i^I)\leq w_{max} \label{upper_waiting_charge_multi}\\
N_i^I=w_i^d\sum_{j=1}^M \lambda_{ij} \label{demand_constraint_charge2_multi}\\
 {N_0}{F_d}\left(q\right)  =\sum_{i=1}^M \sum_{j=1}^M \lambda_{ij}t_{ij}+\sum_{i=1}^M \sum_{j=1}^M w_i^p \lambda_{ij}  +  \sum_{i=1}^M \sum_{j=1}^M   w_i^d \lambda_{ij}, \label{suply_conservation_const_charge_multi}  \\
 {N_{\mathcal{C}}=\sum_{i=1}^M \sum_{j=1}^M \lambda_{ij}d_{ij}^{\mathcal{C}}\dfrac{1}{v_c}+\sum_{i\in \mathcal{C}}  \sum_{j=1}^M w_i^p \lambda_{ij}  +  \sum_{i\in \mathcal{C}} \sum_{j=1}^M   w_i^d \lambda_{ij}} \label{_supply_conservation2_charge_multi} \\
\sum_{j=1}^M  (\lambda_{ji}+\tilde{f}_{ji})=\sum_{j=1}^M  (\lambda_{ij}+\tilde{f}_{ij}), \quad  \forall i\in \mathcal{V}. \label{reballance_const_charge_multi} \\
\tilde{f}_{ij}= \sum_{k: j\in \mathcal{P}_{ik}}    \Pi(j|i\rightarrow k)\mathbb{\tilde{P}}_{ik} \sum_{m=1}^M \lambda_{mi}  \label{definitionofflow_charge_multi}
\end{subnumcases}
}
In (\ref{definitionofflow_charge_multi}), the modified driver repositioning probability $\mathbb{\tilde{P}}_{ij}$ under congestion charge is defined as
\begin{align}
\begin{cases}
&\mathbb{\tilde{P}}_{ij}=\dfrac{e^{\eta (\bar{e}_j-\mathbbm{1}^d_{ij}\bar{p}_{ij})/ (t_{ij}+w_j^d+\bar{t}_j)}}{\sum_{k\neq i} e^{\eta(\bar{e}_k-\mathbbm{1}^d_{ik}\bar{p}_{ik})/ (t_{ik}+w_k^d+\bar{t}_k)}+e^{\eta \bar{e}_i/ (w_i^d+\bar{t}_i)}}, \quad  j\neq i, \\
&\mathbb{\tilde{P}}_{ii}=\dfrac{e^{\eta \bar{e}_i/ (w_i^d+\bar{t}_i)}}{\sum_{k\neq i} e^{\eta(\bar{e}_k-\mathbbm{1}^d_{ik}\bar{p}_{ik})/ (t_{ik}+w_k^d+\bar{t}_k)}+e^{\eta \bar{e}_i/ (w_i^d+\bar{t}_i)}}.
\end{cases}
\label{modifiedlogit_charge}
\end{align}

{
Note that when the congestion charge is imposed on passengers, it enters the passenger travel cost and modifies the demand function (\ref{demand_constraint_charge_multi}), whereas when it is imposed on idle vehicles, it enters both driver repositioning model (\ref{modifiedlogit_charge}) and the driver supply model (\ref{suply_conservation_const_charge_multi}). However, to simplify the notation, we include the congestion charge (on idle drivers) as the third term of  the objective function (\ref{optimalpricing_trip_charge_multizone}) instead of in (\ref{suply_conservation_const_charge_multi}). We comment that these two formulations are mathematically  equivalent: by a change of variable, we can prove that subtracting a total congestion charge from the platform's objective (\ref{optimalpricing_trip_charge_multizone}) is equivalent to subtracting a per-driver congestion charge from the driver supply function (\ref{suply_conservation_const_charge_multi}) \cite[Chap. 16]{varian2014intermediate}.
}

Let $\mathbbm{1}^p$ and $\mathbbm{1}^d$ denote the matrix form of the indicator functions, where the $i$th row and $j$th column is $\mathbbm{1}^p_{ij}$ and $\mathbbm{1}^d_{ij}$, respectively. We have the following three cases:
\begin{itemize}
    \item for one-directional congestion charge, $\mathbbm{1}^p$ and $\mathbbm{1}^d$ are defined as
    \begin{equation}
        \mathbbm{1}^p_{ij}=\mathbbm{1}^d_{ij}=\begin{cases} 
        1, \quad \text{if } i \in \mathcal{R} \text{ and } j\in \mathcal{C}\\
        0, \quad \text{Otherwise} 
        \end{cases}
    \end{equation}
    \item for bi-directional congestion charge, $\mathbbm{1}^p$ and $\mathbbm{1}^d$ are defined as
    \begin{equation}
        \mathbbm{1}^p_{ij}=\mathbbm{1}^d_{ij}=\begin{cases} 
        1, \quad \text{if } i \in \mathcal{R} \text{ and } j\in \mathcal{C}\\
        1, \quad \text{if } i \in \mathcal{C} \text{ and } j\in \mathcal{R}\\
        0, \quad \text{Otherwise} 
        \end{cases}
    \end{equation}   
    \item for trip-based congestion charge, $\mathbbm{1}^p$ and $\mathbbm{1}^d$ are defined as
    \begin{equation}
        \mathbbm{1}^p_{ij}=1,\quad \mathbbm{1}^d=0.
    \end{equation}
\end{itemize}
In all cases, the optimal spatial pricing problem under the congestion charge is a non-convex program. The optimal solution to (\ref{optimalpricing_trip_charge_multizone}) and the upper bound can be obtained by Algorithm \ref{algorithm1}.

\subsection{Case Studies}
\label{simulation_section}
{
We present a case study for San Francisco to  evaluate the impacts of congestion charges on the ride-sourcing market. The northeast corner of San Francisco is the most congested area. Therefore, we define the congestion area $\mathcal{C}$ as consisting of zip code zones 94111, 94104, 94105, 94133, 94108, 94109, 94102, 94103 and 94107 (see Figure \ref{SFzip}).  In the case study, we will use the same passenger demand model, driver supply model, and  pickup time model as in (\ref{logit_demand}), (\ref{logit_supply}), and (\ref{pickuptime_func}), respectively.  We will also use the same set of model parameters as in Section \ref{numericalexample} and solve the optimal spatial pricing problem (\ref{optimalpricing_trip_charge_multizone}) based on Algorithm~\ref{algorithm1}. {The algorithm was executed in Matlab on a Dell desktop machine with 8-core i7-9700 CPU (up to 3.00 GHz)}. A few sets of simulation results will be presented, including  (a) the spatial distribution of ride-sourcing vehicles, passenger demand and ride fare before and after a congestion charge of $\bar{p}_{ij}=\$3$; (b) comparison of distinct congestion charge schemes in terms of congestion mitigation, tax revenue, and tax incidence; (c) sensitivity analysis that investigates how model parameters affect the insights and conclusions. 
}

\subsubsection{Spatial distribution before and after the congestion charge}
\label{casestudy_part1}
We first compare the spatial distribution of supply and demand in the ride-sourcing market before and after a congestion charge of $\bar{p}_{ij}=\$3$. The comparison is performed for each congestion charge scheme, and the {\em changes} (in percentage) in vehicle distributions, ride fare, and passenger demand are summarized in Figure \ref{figure1_mz}-\ref{figure12_mz}. The numerical values are shown in Table 2-4 in Appendix C. 

\begin{figure*}
\begin{minipage}[b]{0.445\linewidth}
\centering
\includegraphics[width = 1.0\linewidth]{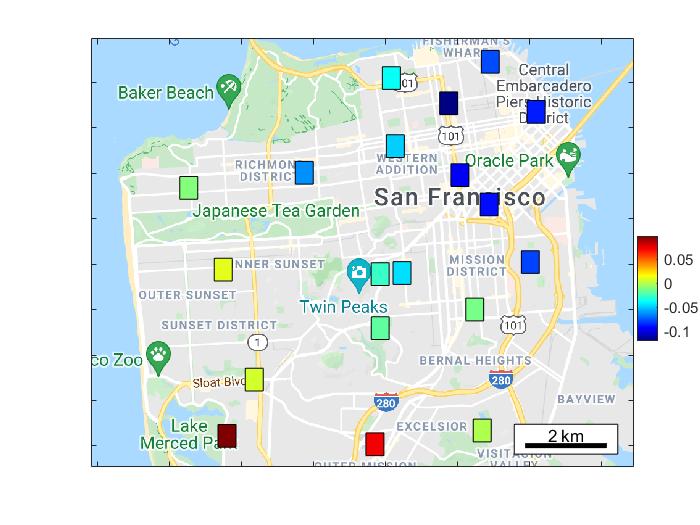}
\caption{\% Change in idle driver distribution  after the one-directional cordon charge of $\$3$.} 
\label{figure1_mz}
\end{minipage}
\begin{minipage}[b]{0.005\linewidth}
\hfill
\end{minipage}
\begin{minipage}[b]{0.445\linewidth}
\centering
\includegraphics[width = 1\linewidth]{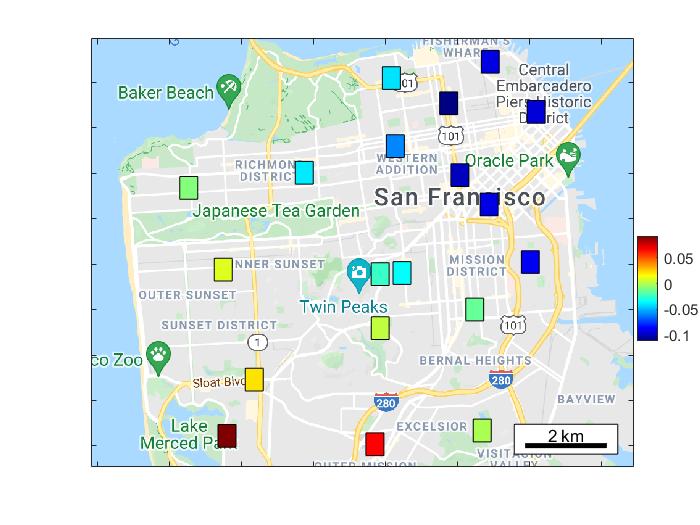}
\caption{\% Change in total number of vehicles after the one-directional cordon charge of $\$3$.} 
\label{N_total1}
\end{minipage}
\end{figure*}

\begin{figure*}
\begin{minipage}[b]{0.445\linewidth}
\centering
\includegraphics[width = 1.0\linewidth]{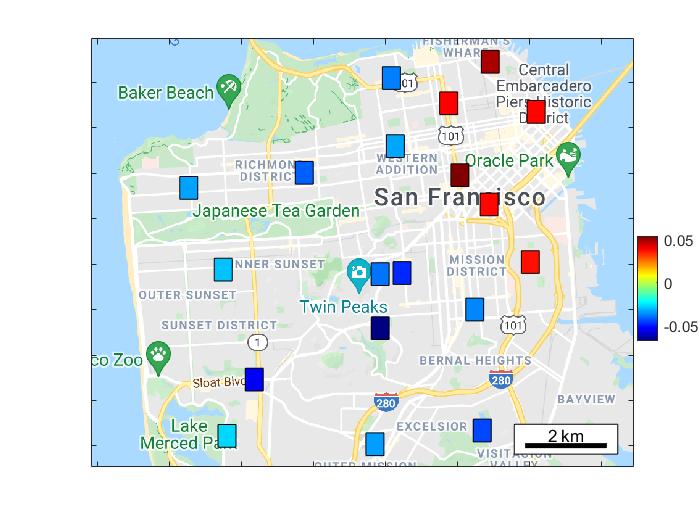}
\caption{\% Change in ride fare  after the one-directional cordon charge of $\$3$.} 
\label{figure2_mz}
\end{minipage}
\begin{minipage}[b]{0.005\linewidth}
\hfill
\end{minipage}
\begin{minipage}[b]{0.445\linewidth}
\centering
\includegraphics[width = 1\linewidth]{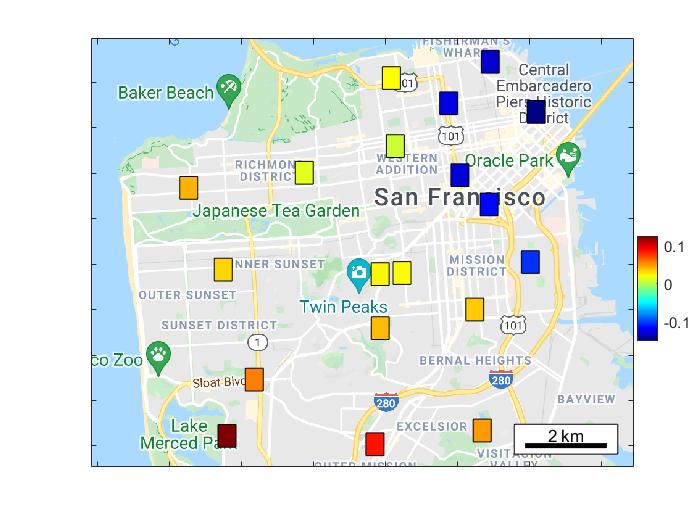}
\caption{\% Change in outflow to uncongested area  after the one-directional cordon charge of $\$3$.} 
\label{figure3_mz}
\end{minipage}
\end{figure*}

\begin{figure*}
\begin{minipage}[b]{0.445\linewidth}
\centering
\includegraphics[width = 1\linewidth]{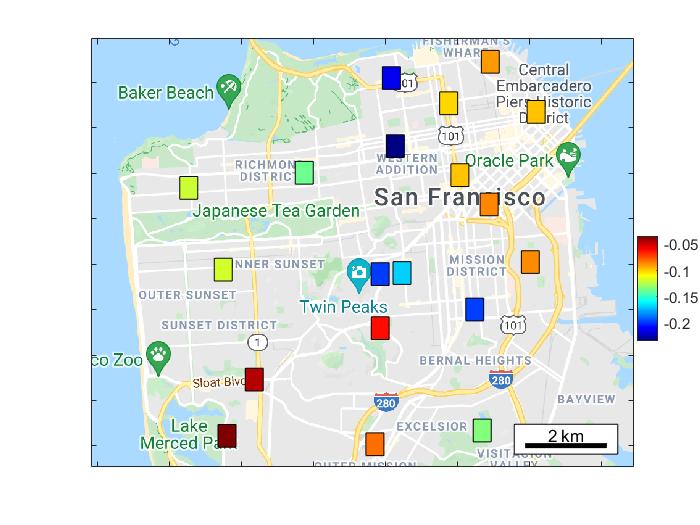}
\caption{\% Change in outflow to congested area after the one-directional cordon charge of $\$3$.} 
\label{figure4_mz}
\end{minipage}
\begin{minipage}[b]{0.005\linewidth}
\hfill
\end{minipage}
\begin{minipage}[b]{0.445\linewidth}
\centering
\includegraphics[width = 1\linewidth]{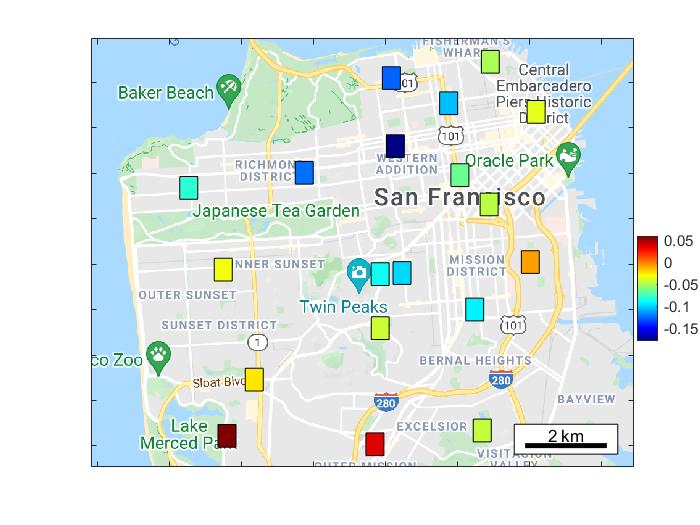}
\caption{\% Change in idle driver distribution  after the bi-directional cordon charge of $\$3$.} 
\label{figure5_mz}
\end{minipage}
\end{figure*}

\begin{figure*}
\begin{minipage}[b]{0.445\linewidth}
\centering
\includegraphics[width = 1\linewidth]{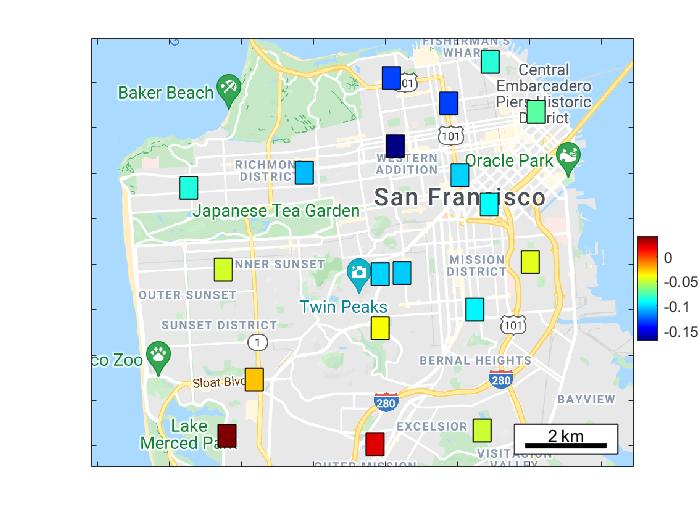}
\caption{\% Change in total number of vehicles after the bi-directional cordon charge of $\$3$.} 
\label{N_total2}
\end{minipage}
\begin{minipage}[b]{0.005\linewidth}
\hfill
\end{minipage}
\begin{minipage}[b]{0.445\linewidth}
\centering
\includegraphics[width = 1\linewidth]{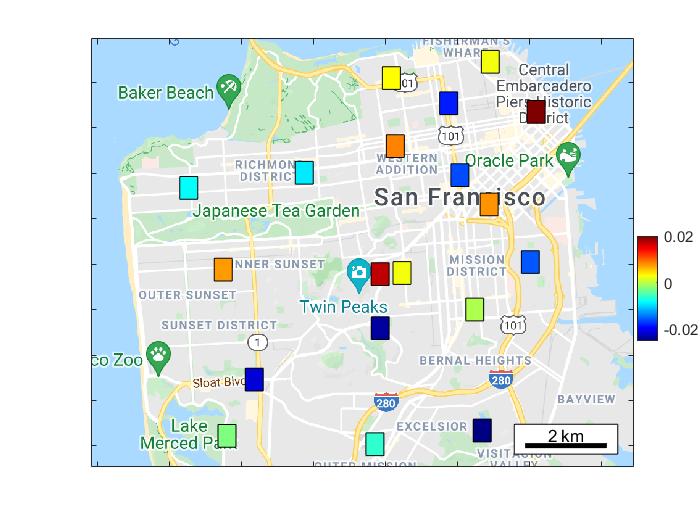}
\caption{\% Change in ride fare  after the bi-directional cordon charge of $\$3$.} 
\label{figure6_mz}
\end{minipage}
\end{figure*}

\begin{figure*}
\begin{minipage}[b]{0.445\linewidth}
\centering
\includegraphics[width = 1\linewidth]{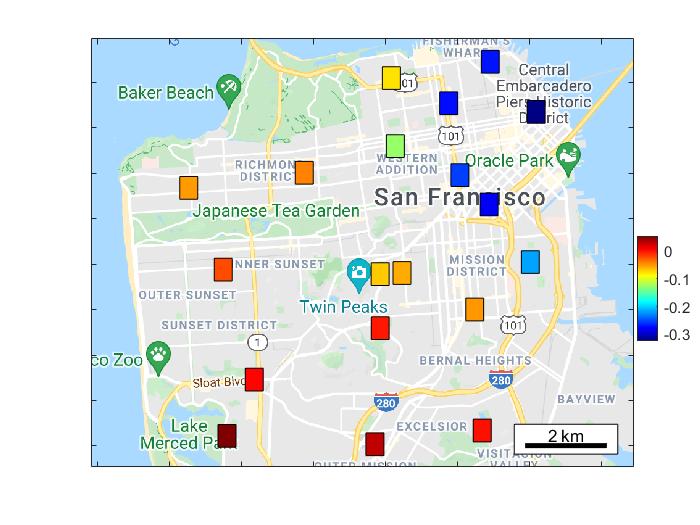}
\caption{\% Change in outflow to uncongested area  after the bi-directional cordon charge of $\$3$.} 
\label{figure7_mz}
\end{minipage}
\begin{minipage}[b]{0.005\linewidth}
\hfill
\end{minipage}
\begin{minipage}[b]{0.445\linewidth}
\centering
\includegraphics[width = 1\linewidth]{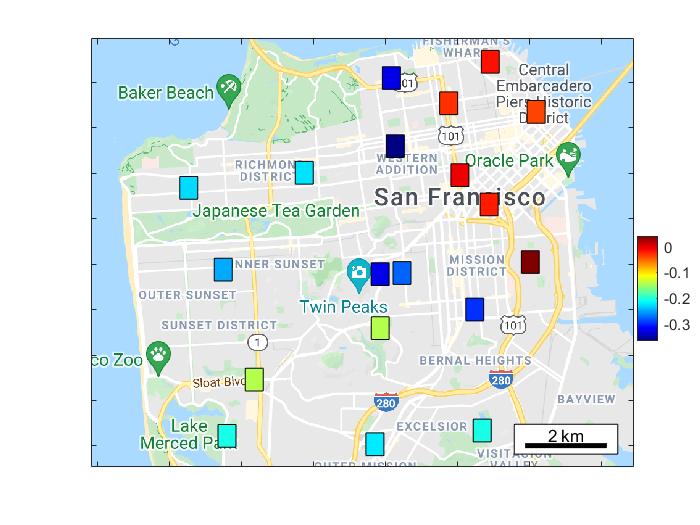}
\caption{\% Change in outflow to congested area after the bi-directional cordon charge of $\$3$.} 
\label{figure8_mz}
\end{minipage}
\end{figure*}

\begin{figure*}
\begin{minipage}[b]{0.445\linewidth}
\centering
\includegraphics[width = 1\linewidth]{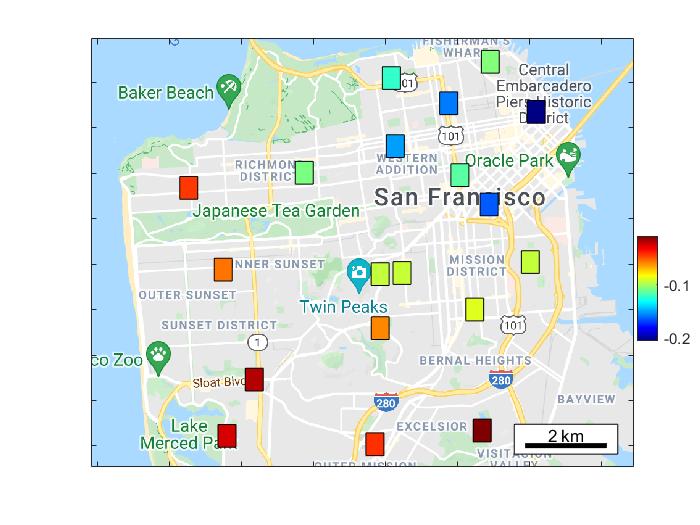}
\caption{\% Change in idle driver distribution  after the trip-based charge of $\$3$.} 
\label{figure9_mz}
\end{minipage}
\begin{minipage}[b]{0.005\linewidth}
\hfill
\end{minipage}
\begin{minipage}[b]{0.445\linewidth}
\centering
\includegraphics[width = 1\linewidth]{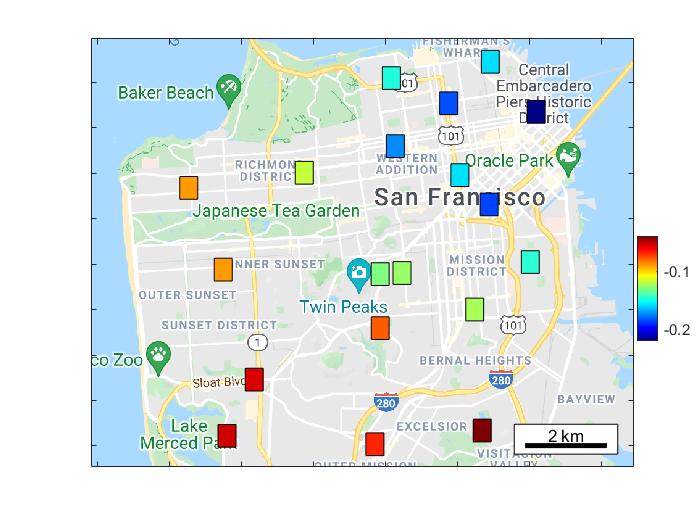}
\caption{\% Change in total number of vehicles after the trip-based charge of $\$3$.} 
\label{N_total3}
\end{minipage}
\end{figure*}

\begin{figure*}
\begin{minipage}[b]{0.445\linewidth}
\centering
\includegraphics[width = 1\linewidth]{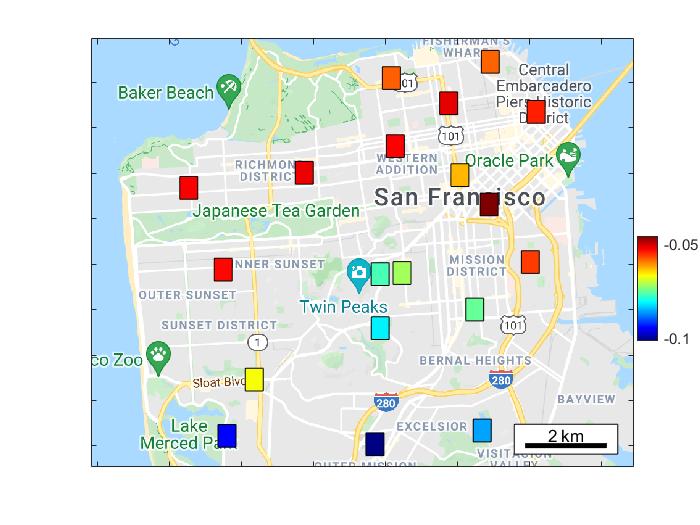}
\caption{\% Change in ride fare  after the trip-based charge of $\$3$.} 
\label{figure10_mz}
\end{minipage}
\begin{minipage}[b]{0.005\linewidth}
\hfill
\end{minipage}
\begin{minipage}[b]{0.445\linewidth}
\centering
\includegraphics[width = 1\linewidth]{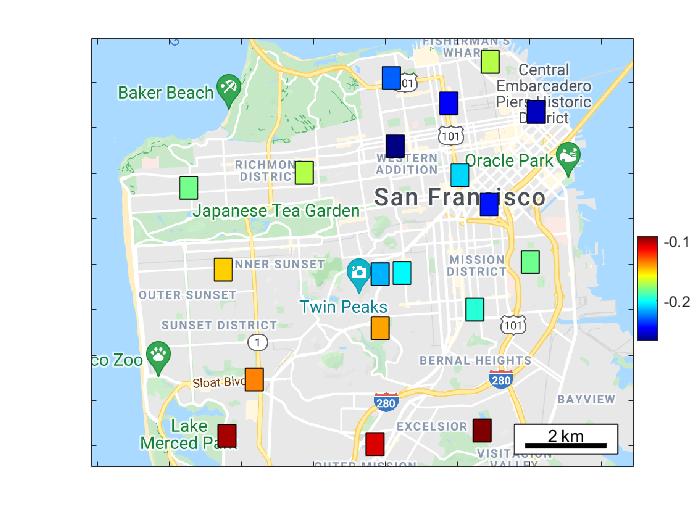}
\caption{\% Change in outflow to uncongested area  after the trip-based charge of $\$3$.} 
\label{figure11_mz}
\end{minipage}
\end{figure*}

\begin{figure}
\centering
\includegraphics[width = 0.5\linewidth]{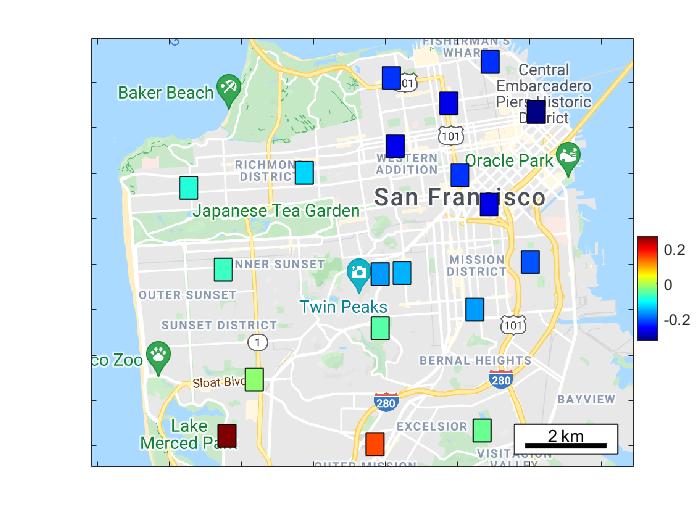}
\caption{\% Change in outflow to congested area after the trip-based charge of $\$3$.} 
\label{figure12_mz}
\end{figure}

{
Figure \ref{figure1_mz}-\ref{figure4_mz} show the impact of the one-directional cordon price. These results provide a few interesting insights, which are listed below:
\begin{itemize}
\item  It is clear that the one-directional cordon price reduces the number of idle vehicles in the congestion area (see Figure \ref{figure1_mz}) and thus reduces the total ride-sourcing traffic in the urban core (see Figure \ref{N_total1}).
\item The number of idle vehicles in some of the remote zones increases (Figure \ref{figure1_mz}). This indicates that the waiting time in those areas reduces and the service quality is improved. We point out that this is partially because the one-directional cordon price discourages idle drivers from entering the congestion area and motivates them to serve more passengers in the underserved zones. Furthermore, Figure \ref{figure2_mz} shows that under the one-directional cordon price, the ride fare in the remote area decreases.\footnote{Note that the ride fare $r_i$ shown in Figure \ref{figure2_mz}, Figure \ref{figure6_mz} and Figure \ref{figure10_mz} do not include the cordon charge for trips entering the congestion area.} Together with the increased number of idle vehicles, this indicates that the generalized travel cost for passengers starting from the remote area reduces, leading to higher passenger demand in $\mathcal{R}$ (see Figure \ref{figure3_mz}). 
\item In Figure \ref{figure1_mz} and  Figure \ref{figure3_mz}, we observe that zones further from the congestion area benefit more than zones closer to the congestion area. We validate that this is because those zones have fewer trips to the congestion area, so that they are less affected by the cordon charge. 
\item In Figure \ref{figure2_mz} and Figure \ref{figure3_mz}, the ride fare in the congestion area increases and the passenger trips from the congestion area to the remote area significantly reduce. This is somewhat surprising because these trips are not directly penalized under the one-directional cordon price. We conjecture that this is  because each trip from the congestion area to the remote area would necessitate some rebalancing flow back to the congestion area, which incurs an extra charge for passing the cordon in the future. 
\end{itemize}
}

{
Figures \ref{figure5_mz}-\ref{figure8_mz} show the impact of the bi-directional cordon price on the ride-sourcing market. The bi-directional cordon price discourages passengers and vehicles from moving between the congested and uncongested areas. Therefore, cross-cordon trips are severely penalized and the corresponding demand  reduces (Figure \ref{figure7_mz}-\ref{figure8_mz}). It is interesting to note that the impact of the bi-directional cordon price on each zone depends on how many trips starting from this zone cross the cordon. If the percentage of cross-cordon trips is small, then the congestion charge may have a smaller impact. For instance, we observe that for very remote zones that are far from the cordon, the number of idle vehicles may  increase (see Figure \ref{figure5_mz}) and the trip fare may reduce (see Figure \ref{figure6_mz}), which benefits passengers in these zones. It is also interesting to observe that although trips that cross the cordon are penalized and the corresponding demand reduces, trips that do not cross the cordon are not significantly affected, and the corresponding demands for these trips even increases in certain zones (Figure \ref{figure7_mz}-\ref{figure8_mz}).
}

{
Figures \ref{figure9_mz}-\ref{figure12_mz} show the impact of the trip-based congestion charge on the ride-sourcing market. The trip-based congestion charge penalizes all trips in the city, which curbs the overall ride-sourcing demand and reduces the market size. However, it is interesting to note that since the congestion area has a much larger demand, it assumes a larger tax incidence than other zones in the remote area. For instance, in the congestion area,  the number of ride-sourcing vehicles is reduceed (Figure \ref{figure9_mz} and Figure \ref{N_total3}), and the passenger demand is curbed (Figure \ref{figure11_mz} and Figure \ref{figure12_mz}). At the same time, in some zones of the remote area, the number of vehicles increases (Figure \ref{figure9_mz} and Figure \ref{N_total3}), the trip fare reduces (Figure \ref{figure10_mz}), and the trip volume also increases (Figure \ref{figure12_mz}). 
}

\begin{figure*}%
\begin{minipage}[b]{0.32\linewidth}
\centering
%
%
\definecolor{mycolor1}{rgb}{0.00000,0.44700,0.74100}%
\definecolor{mycolor2}{rgb}{0.85000,0.32500,0.09800}%
\begin{tikzpicture}

\begin{axis}[%
width=1.694in,
height=1.03in,
at={(1.358in,0.0in)},
scale only axis,
xmin=0,
xmax=3,
xlabel style={font=\color{white!15!black}},
xlabel={Congestion charge},
ymin=500,
ymax=650,
ytick={500,550, 600, 650},
yticklabels={{500},{550},{600}, {650}},
ylabel style={font=\color{white!15!black}},
ylabel={Idle driver in $\mathcal{C}$ },
axis background/.style={fill=white},
legend style={at={(0.0,0.0)}, anchor=south west, legend cell align=left, align=left, font=\small, draw=white!15!black, nodes={scale=0.8, transform shape} }
]
\addplot [color=mycolor1, line width=1.0pt]
  table[row sep=crcr]{%
0	618.662651954089\\
0.157894736842105	615.648248643792\\
0.315789473684211	612.632558711555\\
0.473684210526316	609.616316148785\\
0.631578947368421	606.599971386652\\
0.789473684210526	603.583625532213\\
0.947368421052632	600.567876133337\\
1.10526315789474	597.553147125733\\
1.26315789473684	594.539565261862\\
1.42105263157895	591.528051674041\\
1.57894736842105	588.518182435831\\
1.73684210526316	585.511113030239\\
1.89473684210526	582.506665388791\\
2.05263157894737	579.505061726784\\
2.21052631578947	576.506846422025\\
2.36842105263158	573.512702608561\\
2.52631578947368	570.522096981582\\
2.68421052631579	567.535666849029\\
2.84210526315789	564.553760775178\\
3	561.576331234704\\
};
\addlegendentry{Cordon 1}

\addplot [color=mycolor2, line width=1.0pt]
  table[row sep=crcr]{%
0	618.662651954089\\
0.157894736842105	616.460109716915\\
0.315789473684211	614.304998329263\\
0.473684210526316	612.197019641937\\
0.631578947368421	610.136590177002\\
0.789473684210526	608.123349833991\\
0.947368421052632	606.15743744911\\
1.10526315789474	604.238528383421\\
1.26315789473684	602.366239950006\\
1.42105263157895	600.540617465486\\
1.57894736842105	598.761099061734\\
1.73684210526316	597.027361614632\\
1.89473684210526	595.338971007116\\
2.05263157894737	593.695713692531\\
2.21052631578947	592.096875724595\\
2.36842105263158	590.542131801308\\
2.52631578947368	589.030771546312\\
2.68421052631579	587.562184630569\\
2.84210526315789	586.136418917652\\
3	584.751469261377\\
};
\addlegendentry{Cordon 2}

\addplot [color=black, line width=1.0pt]
  table[row sep=crcr]{%
0	618.662651954089\\
0.157894736842105	614.233103630081\\
0.315789473684211	609.797169784203\\
0.473684210526316	605.355709385811\\
0.631578947368421	600.908345366124\\
0.789473684210526	596.455441062282\\
0.947368421052632	591.997242152943\\
1.10526315789474	587.534243985861\\
1.26315789473684	583.066066762463\\
1.42105263157895	578.593297841928\\
1.57894736842105	574.11597530005\\
1.73684210526316	569.634492359063\\
1.89473684210526	565.148798985641\\
2.05263157894737	560.659152940265\\
2.21052631578947	556.166348100652\\
2.36842105263158	551.668814930368\\
2.52631578947368	547.168658909525\\
2.68421052631579	542.665083189333\\
2.84210526315789	538.158375560061\\
3	533.648660882528\\
};
\addlegendentry{Trip-based}

\end{axis}
\end{tikzpicture}%
\vspace*{-0.3in}
\caption{Number of idle drivers in $\mathcal{C}$ under different surcharge $\bar{p}_{ij}$.}
\label{figure1_2z}
\end{minipage}
\begin{minipage}[b]{0.005\linewidth}
\hfill
\end{minipage}
\begin{minipage}[b]{0.32\linewidth}
\centering
%
%
\definecolor{mycolor1}{rgb}{0.00000,0.44700,0.74100}%
\definecolor{mycolor2}{rgb}{0.85000,0.32500,0.09800}%
\begin{tikzpicture}

\begin{axis}[%
width=1.694in,
height=1.03in,
at={(1.358in,0.0in)},
scale only axis,
xmin=0,
xmax=3,
xlabel style={font=\color{white!15!black}},
xlabel={Congestion charge},
ymin=1.75,
ymax=2.15,
ylabel style={font=\color{white!15!black}},
ylabel={Ride fare in $\mathcal{C}$},
axis background/.style={fill=white},
legend style={at={(0.0,0.0)}, anchor=south west, legend cell align=left, align=left, font=\small, draw=white!15!black, nodes={scale=0.8, transform shape} }
]
\addplot [color=mycolor1, line width=1.0pt]
  table[row sep=crcr]{%
0	2.01851409725453\\
0.157894736842105	2.02300689995768\\
0.315789473684211	2.02749700206045\\
0.473684210526316	2.03198609533083\\
0.631578947368421	2.03647459622941\\
0.789473684210526	2.04096416246383\\
0.947368421052632	2.04545354677128\\
1.10526315789474	2.04994547757277\\
1.26315789473684	2.05443982290332\\
1.42105263157895	2.05893813527448\\
1.57894736842105	2.06344041011159\\
1.73684210526316	2.06794827858825\\
1.89473684210526	2.07246168711765\\
2.05263157894737	2.07698170394867\\
2.21052631578947	2.08150902014597\\
2.36842105263158	2.08604492994211\\
2.52631578947368	2.0905895798031\\
2.68421052631579	2.09514364193267\\
2.84210526315789	2.09970808063178\\
3	2.10428324354223\\
};
\addlegendentry{Cordon 1}

\addplot [color=mycolor2, line width=1.0pt]
  table[row sep=crcr]{%
0	2.01851409725453\\
0.157894736842105	2.01829849094418\\
0.315789473684211	2.01804259825883\\
0.473684210526316	2.0177513128229\\
0.631578947368421	2.0174301335436\\
0.789473684210526	2.01708342031318\\
0.947368421052632	2.0167162437503\\
1.10526315789474	2.01633299775534\\
1.26315789473684	2.01593777283254\\
1.42105263157895	2.01553490999155\\
1.57894736842105	2.01512816035416\\
1.73684210526316	2.01472123225027\\
1.89473684210526	2.0143177611428\\
2.05263157894737	2.01392114940457\\
2.21052631578947	2.01353429875984\\
2.36842105263158	2.01316058918113\\
2.52631578947368	2.01280258405455\\
2.68421052631579	2.01246283550257\\
2.84210526315789	2.01214345770497\\
3	2.01184725205871\\
};
\addlegendentry{Cordon 2}

\addplot [color=black, line width=1.0pt]
  table[row sep=crcr]{%
0	2.01851409725453\\
0.157894736842105	2.01253822961168\\
0.315789473684211	2.00657930229327\\
0.473684210526316	2.00063768564056\\
0.631578947368421	1.99471302004339\\
0.789473684210526	1.98880546419894\\
0.947368421052632	1.98291508068646\\
1.10526315789474	1.97704202267825\\
1.26315789473684	1.97118677577616\\
1.42105263157895	1.96534894285915\\
1.57894736842105	1.95952855687875\\
1.73684210526316	1.95372619298556\\
1.89473684210526	1.94794126808263\\
2.05263157894737	1.94217449629576\\
2.21052631578947	1.93642611755015\\
2.36842105263158	1.93069508714188\\
2.52631578947368	1.92498320012697\\
2.68421052631579	1.91928937401504\\
2.84210526315789	1.91361421378427\\
3	1.90795768935083\\
};
\addlegendentry{Trip-based}

\end{axis}
\end{tikzpicture}%
\vspace*{-0.3in}
\caption{Average ride fare in $\mathcal{C}$ under different surcharge $\bar{p}_{ij}$.} 
\label{figure2_2z}
\end{minipage}
\begin{minipage}[b]{0.005\linewidth}
\hfill
\end{minipage}
\begin{minipage}[b]{0.32\linewidth}
\centering
%
%
\definecolor{mycolor1}{rgb}{0.00000,0.44700,0.74100}%
\definecolor{mycolor2}{rgb}{0.85000,0.32500,0.09800}%
\begin{tikzpicture}

\begin{axis}[%
width=1.694in,
height=1.03in,
at={(1.358in,0.0in)},
scale only axis,
xmin=0,
xmax=3,
xlabel style={font=\color{white!15!black}},
xlabel={Congestion charge},
ymin=4.2,
ymax=4.8,
ylabel style={font=\color{white!15!black}},
ylabel={Waiting time in $\mathcal{C}$},
axis background/.style={fill=white},
legend style={at={(0.0,0.54)}, anchor=south west, legend cell align=left, align=left, font=\small, draw=white!15!black, nodes={scale=0.8, transform shape} }
]
\addplot [color=mycolor1, line width=1.0pt]
  table[row sep=crcr]{%
0	4.28240136407698\\
0.157894736842105	4.29266592365695\\
0.315789473684211	4.30301887014202\\
0.473684210526316	4.3134587310562\\
0.631578947368421	4.32398491102742\\
0.789473684210526	4.33459824645083\\
0.947368421052632	4.34529726677012\\
1.10526315789474	4.35608181580441\\
1.26315789473684	4.36695222416962\\
1.42105263157895	4.37790626763607\\
1.57894736842105	4.38894639493401\\
1.73684210526316	4.40006934110615\\
1.89473684210526	4.4112767224534\\
2.05263157894737	4.42256858580365\\
2.21052631578947	4.43394384577845\\
2.36842105263158	4.44540091782939\\
2.52631578947368	4.45694263969561\\
2.68421052631579	4.46856757099324\\
2.84210526315789	4.48027516037781\\
3	4.49206652837494\\
};
\addlegendentry{Cordon 1}

\addplot [color=mycolor2, line width=1.0pt]
  table[row sep=crcr]{%
0	4.28240136407698\\
0.157894736842105	4.28946901008936\\
0.315789473684211	4.29642268452766\\
0.473684210526316	4.30326145859842\\
0.631578947368421	4.30998197770084\\
0.789473684210526	4.31658364973483\\
0.947368421052632	4.3230639017887\\
1.10526315789474	4.32942213591265\\
1.26315789473684	4.33565771566448\\
1.42105263157895	4.34176858341042\\
1.57894736842105	4.34775485770505\\
1.73684210526316	4.35361585083744\\
1.89473684210526	4.35935132588171\\
2.05263157894737	4.36496012100971\\
2.21052631578947	4.37044306762161\\
2.36842105263158	4.37579960846831\\
2.52631578947368	4.38103056409383\\
2.68421052631579	4.38613637437846\\
2.84210526315789	4.39111564179067\\
3	4.39597327531296\\
};
\addlegendentry{Cordon 2}

\addplot [color=black, line width=1.0pt]
  table[row sep=crcr]{%
0	4.28240136407698\\
0.157894736842105	4.29832600343696\\
0.315789473684211	4.31446882224141\\
0.473684210526316	4.33083191898913\\
0.631578947368421	4.34742148281783\\
0.789473684210526	4.36424136767045\\
0.947368421052632	4.38129604377977\\
1.10526315789474	4.39858916969853\\
1.26315789473684	4.4161273668173\\
1.42105263157895	4.4339146339911\\
1.57894736842105	4.45195629419286\\
1.73684210526316	4.47025721161729\\
1.89473684210526	4.48882370384306\\
2.05263157894737	4.50766098171398\\
2.21052631578947	4.52677265960936\\
2.36842105263158	4.54617153105698\\
2.52631578947368	4.56585635368958\\
2.68421052631579	4.58583738053463\\
2.84210526315789	4.6061208258358\\
3	4.62671387158558\\
};
\addlegendentry{Trip-based}

\end{axis}
\end{tikzpicture}%
\vspace*{-0.3in}
\caption{Average passenger waiting time in $\mathcal{C}$ under different surcharge $\bar{p}_{ij}$.}
\label{figure3_2z}
\end{minipage}
\begin{minipage}[b]{0.32\linewidth}
\centering
%
%
\definecolor{mycolor1}{rgb}{0.00000,0.44700,0.74100}%
\definecolor{mycolor2}{rgb}{0.85000,0.32500,0.09800}%
\begin{tikzpicture}

\begin{axis}[%
width=1.694in,
height=1.03in,
at={(1.358in,0.0in)},
scale only axis,
xmin=0,
xmax=3,
xlabel style={font=\color{white!15!black}},
xlabel={Congestion charge},
ymin=800,
ymax=980,
ylabel style={font=\color{white!15!black}},
ylabel={Idle driver in $\mathcal{R}$},
axis background/.style={fill=white},
legend style={at={(0.0,0.0)}, anchor=south west, legend cell align=left, align=left, font=\small, draw=white!15!black, nodes={scale=0.8, transform shape} }
]
\addplot [color=mycolor1, line width=1.0pt]
  table[row sep=crcr]{%
0	959.858792163012\\
0.157894736842105	958.874099420824\\
0.315789473684211	957.90237968227\\
0.473684210526316	956.9441165231\\
0.631578947368421	955.999096900061\\
0.789473684210526	955.066308899492\\
0.947368421052632	954.146146685176\\
1.10526315789474	953.23863838503\\
1.26315789473684	952.342503430713\\
1.42105263157895	951.459137505661\\
1.57894736842105	950.586502901196\\
1.73684210526316	949.725890990332\\
1.89473684210526	948.876069027438\\
2.05263157894737	948.036721570629\\
2.21052631578947	947.207747861271\\
2.36842105263158	946.390352225723\\
2.52631578947368	945.581566231843\\
2.68421052631579	944.783238299587\\
2.84210526315789	943.994679278945\\
3	943.214723860882\\
};
\addlegendentry{Cordon 1}

\addplot [color=mycolor2, line width=1.0pt]
  table[row sep=crcr]{%
0	959.858792163012\\
0.157894736842105	955.440059687664\\
0.315789473684211	951.056805715625\\
0.473684210526316	946.709122794105\\
0.631578947368421	942.397854713331\\
0.789473684210526	938.123120985959\\
0.947368421052632	933.884927792018\\
1.10526315789474	929.683449479469\\
1.26315789473684	925.517950588148\\
1.42105263157895	921.389462726314\\
1.57894736842105	917.296771047958\\
1.73684210526316	913.24028498632\\
1.89473684210526	909.219642066969\\
2.05263157894737	905.23515220758\\
2.21052631578947	901.285681179707\\
2.36842105263158	897.371793922736\\
2.52631578947368	893.492614596203\\
2.68421052631579	889.647710826612\\
2.84210526315789	885.838013202516\\
3	882.061590484891\\
};
\addlegendentry{Cordon 2}

\addplot [color=black, line width=1.0pt]
  table[row sep=crcr]{%
0	959.858792163012\\
0.157894736842105	955.945605856114\\
0.315789473684211	952.00072768945\\
0.473684210526316	948.025390685713\\
0.631578947368421	944.019328775155\\
0.789473684210526	939.982056399803\\
0.947368421052632	935.914506048519\\
1.10526315789474	931.817290644043\\
1.26315789473684	927.689429838103\\
1.42105263157895	923.531450099757\\
1.57894736842105	919.343645156178\\
1.73684210526316	915.127791107487\\
1.89473684210526	910.881468203869\\
2.05263157894737	906.60620385758\\
2.21052631578947	902.303347204328\\
2.36842105263158	897.969912697886\\
2.52631578947368	893.609757142696\\
2.68421052631579	889.220703780115\\
2.84210526315789	884.803878166261\\
3	880.358710546468\\
};
\addlegendentry{Trip-based}

\end{axis}
\end{tikzpicture}%
\vspace*{-0.3in}
\caption{Number of idle drivers in $\mathcal{R}$ under different surcharge $\bar{p}_{ij}$.} 
\label{figure4_2z}
\end{minipage}
\begin{minipage}[b]{0.005\linewidth}
\hfill
\end{minipage}
\begin{minipage}[b]{0.32\linewidth}
\centering
%
%
\definecolor{mycolor1}{rgb}{0.00000,0.44700,0.74100}%
\definecolor{mycolor2}{rgb}{0.85000,0.32500,0.09800}%
\begin{tikzpicture}

\begin{axis}[%
width=1.694in,
height=1.03in,
at={(1.358in,0.0in)},
scale only axis,
xmin=0,
xmax=3,
xlabel style={font=\color{white!15!black}},
xlabel={Congestion surcharge},
ymin=1.60,
ymax=1.82,
ylabel style={font=\color{white!15!black}},
ylabel={Ride fare in $\mathcal{R}$},
axis background/.style={fill=white},
legend style={at={(0.0,0.0)}, anchor=south west, legend cell align=left, align=left, font=\small, draw=white!15!black, nodes={scale=0.8, transform shape} }
]
\addplot [color=mycolor1, line width=1.0pt]
  table[row sep=crcr]{%
0	1.80414519802732\\
0.157894736842105	1.80030729283088\\
0.315789473684211	1.79647324536852\\
0.473684210526316	1.79264402278723\\
0.631578947368421	1.78881943175777\\
0.789473684210526	1.78499939750561\\
0.947368421052632	1.78118424619414\\
1.10526315789474	1.77737452245614\\
1.26315789473684	1.77356962328602\\
1.42105263157895	1.76977079240564\\
1.57894736842105	1.76597743966299\\
1.73684210526316	1.76219028471919\\
1.89473684210526	1.75840913007969\\
2.05263157894737	1.75463460794048\\
2.21052631578947	1.75086654621819\\
2.36842105263158	1.747106311269\\
2.52631578947368	1.74335274315268\\
2.68421052631579	1.73960750017279\\
2.84210526315789	1.73587037995953\\
3	1.73214134492897\\
};
\addlegendentry{Cordon 1}

\addplot [color=mycolor2, line width=1.0pt]
  table[row sep=crcr]{%
0	1.80414519802732\\
0.157894736842105	1.80378853275547\\
0.315789473684211	1.80342663819405\\
0.473684210526316	1.80306049347435\\
0.631578947368421	1.80269143829064\\
0.789473684210526	1.80232068366709\\
0.947368421052632	1.80194877202629\\
1.10526315789474	1.8015771052433\\
1.26315789473684	1.80120604995389\\
1.42105263157895	1.80083720046468\\
1.57894736842105	1.80047060786969\\
1.73684210526316	1.8001077174856\\
1.89473684210526	1.79974923999438\\
2.05263157894737	1.79939607194005\\
2.21052631578947	1.79904840570832\\
2.36842105263158	1.7987077267872\\
2.52631578947368	1.79837427664893\\
2.68421052631579	1.79804855110097\\
2.84210526315789	1.79773114972949\\
3	1.79742371799321\\
};
\addlegendentry{Cordon 2}

\addplot [color=black, line width=1.0pt]
  table[row sep=crcr]{%
0	1.80414519802732\\
0.157894736842105	1.79720725957096\\
0.315789473684211	1.79030120986164\\
0.473684210526316	1.78342720621789\\
0.631578947368421	1.77658574962733\\
0.789473684210526	1.76977628060639\\
0.947368421052632	1.76299953896244\\
1.10526315789474	1.75625566447107\\
1.26315789473684	1.74954487299922\\
1.42105263157895	1.74286664649063\\
1.57894736842105	1.73622153784585\\
1.73684210526316	1.72961114031497\\
1.89473684210526	1.72303306594416\\
2.05263157894737	1.71648879517363\\
2.21052631578947	1.70997904061458\\
2.36842105263158	1.70350209687732\\
2.52631578947368	1.69706023941344\\
2.68421052631579	1.69065188592014\\
2.84210526315789	1.68427808263471\\
3	1.67793848925938\\
};
\addlegendentry{Trip-based}

\end{axis}
\end{tikzpicture}%
\vspace*{-0.3in}
\caption{Average ride fare in $\mathcal{R}$ under different surcharge  $\bar{p}_{ij}$.}
\label{figure5_2z}
\end{minipage}
\begin{minipage}[b]{0.005\linewidth}
\hfill
\end{minipage}
\begin{minipage}[b]{0.32\linewidth}
\centering
%
%
\definecolor{mycolor1}{rgb}{0.00000,0.44700,0.74100}%
\definecolor{mycolor2}{rgb}{0.85000,0.32500,0.09800}%
\begin{tikzpicture}

\begin{axis}[%
width=1.694in,
height=1.03in,
at={(1.358in,0.0in)},
scale only axis,
xmin=0,
xmax=3,
xlabel style={font=\color{white!15!black}},
xlabel={Congestion surcharge},
ymin=5.2,
ymax=5.5,
ylabel style={font=\color{white!15!black}},
ylabel={Waiting time in $\mathcal{R}$},
axis background/.style={fill=white},
legend style={at={(0.0,0.54)}, anchor=south west, legend cell align=left, align=left, font=\small, draw=white!15!black, nodes={scale=0.8, transform shape} }
]
\addplot [color=mycolor1, line width=1.0pt]
  table[row sep=crcr]{%
0	5.27785113689752\\
0.157894736842105	5.27785946368729\\
0.315789473684211	5.27787335244842\\
0.473684210526316	5.27788983047807\\
0.631578947368421	5.2779098452571\\
0.789473684210526	5.27793837287638\\
0.947368421052632	5.27796915159819\\
1.10526315789474	5.27800589371939\\
1.26315789473684	5.27805157316319\\
1.42105263157895	5.2781010352045\\
1.57894736842105	5.27816107656635\\
1.73684210526316	5.27822709534722\\
1.89473684210526	5.27830299241385\\
2.05263157894737	5.27838917571257\\
2.21052631578947	5.27848612045129\\
2.36842105263158	5.27858936827341\\
2.52631578947368	5.27870857166919\\
2.68421052631579	5.2788373111367\\
2.84210526315789	5.27897764294909\\
3	5.27913314998544\\
};
\addlegendentry{Cordon 1}

\addplot [color=mycolor2, line width=1.0pt]
  table[row sep=crcr]{%
0	5.27785113689752\\
0.157894736842105	5.28575072948116\\
0.315789473684211	5.29372878499283\\
0.473684210526316	5.30178332337729\\
0.631578947368421	5.30990930254827\\
0.789473684210526	5.31810404949641\\
0.947368421052632	5.32636575644665\\
1.10526315789474	5.33469143298244\\
1.26315789473684	5.34308126876105\\
1.42105263157895	5.35152906934762\\
1.57894736842105	5.36003708852698\\
1.73684210526316	5.36860113316824\\
1.89473684210526	5.37721971013719\\
2.05263157894737	5.38588944463629\\
2.21052631578947	5.39461191146464\\
2.36842105263158	5.40338211932037\\
2.52631578947368	5.41220092623464\\
2.68421052631579	5.4210670658672\\
2.84210526315789	5.42997526235887\\
3	5.43892903412175\\
};
\addlegendentry{Cordon 2}

\addplot [color=black, line width=1.0pt]
  table[row sep=crcr]{%
0	5.27785113689752\\
0.157894736842105	5.28588196030919\\
0.315789473684211	5.29408053849555\\
0.473684210526316	5.30244553128856\\
0.631578947368421	5.31097946392852\\
0.789473684210526	5.31968571955675\\
0.947368421052632	5.32856368465158\\
1.10526315789474	5.33761342121317\\
1.26315789473684	5.34684020058458\\
1.42105263157895	5.35624531297576\\
1.57894736842105	5.36582852266799\\
1.73684210526316	5.37558729168319\\
1.89473684210526	5.38553294311042\\
2.05263157894737	5.39566169024774\\
2.21052631578947	5.40597042039784\\
2.36842105263158	5.41647566365673\\
2.52631578947368	5.42716312288936\\
2.68421052631579	5.43804527096735\\
2.84210526315789	5.44911983660539\\
3	5.4603928831244\\
};
\addlegendentry{Trip-based}

\end{axis}
\end{tikzpicture}%
\vspace*{-0.3in}
\caption{Average passenger waiting time in $\mathcal{R}$ under different surcharge $\bar{p}_{ij}$.}
\label{figure6_2z}
\end{minipage}
\begin{minipage}[b]{0.32\linewidth}
\centering
%
%
\definecolor{mycolor1}{rgb}{0.00000,0.44700,0.74100}%
\definecolor{mycolor2}{rgb}{0.85000,0.32500,0.09800}%
\begin{tikzpicture}

\begin{axis}[%
width=1.694in,
height=1.03in,
at={(1.358in,0.0in)},
scale only axis,
xmin=0,
xmax=3,
xlabel style={font=\color{white!15!black}},
xlabel={Congestion charge},
ymin=35,
ymax=55,
ylabel style={font=\color{white!15!black}},
ylabel={$\lambda_{11}$},
axis background/.style={fill=white},
legend style={at={(0.0,0.0)}, anchor=south west, legend cell align=left, align=left, font=\small, draw=white!15!black, nodes={scale=0.8, transform shape} }
]
\addplot [color=mycolor1, line width=1.0pt]
  table[row sep=crcr]{%
0	51.8020213564446\\
0.157894736842105	51.5501208470144\\
0.315789473684211	51.2980562989856\\
0.473684210526316	51.0458048707291\\
0.631578947368421	50.7933712327562\\
0.789473684210526	50.5407081141476\\
0.947368421052632	50.2878774333462\\
1.10526315789474	50.0348142816872\\
1.26315789473684	49.7815286618761\\
1.42105263157895	49.5280150053179\\
1.57894736842105	49.2742566888721\\
1.73684210526316	49.0202557576683\\
1.89473684210526	48.7660096505329\\
2.05263157894737	48.5114982390457\\
2.21052631578947	48.2567291580848\\
2.36842105263158	48.0016993343848\\
2.52631578947368	47.7463831083776\\
2.68421052631579	47.4907954361354\\
2.84210526315789	47.2349293923264\\
3	46.978777122328\\
};
\addlegendentry{Cordon 1}

\addplot [color=mycolor2, line width=1.0pt]
  table[row sep=crcr]{%
0	51.8020213564446\\
0.157894736842105	51.7512207124355\\
0.315789473684211	51.7027632409081\\
0.473684210526316	51.6565147371297\\
0.631578947368421	51.6123535094477\\
0.789473684210526	51.5701566424752\\
0.947368421052632	51.5298070558776\\
1.10526315789474	51.4911806893828\\
1.26315789473684	51.4541663107711\\
1.42105263157895	51.4186605493109\\
1.57894736842105	51.3845508912896\\
1.73684210526316	51.3517384273146\\
1.89473684210526	51.3201191944401\\
2.05263157894737	51.2896066684026\\
2.21052631578947	51.2601068652225\\
2.36842105263158	51.2315267175314\\
2.52631578947368	51.2037813204128\\
2.68421052631579	51.1767921038396\\
2.84210526315789	51.1505106223369\\
3	51.1248019057998\\
};
\addlegendentry{Cordon 2}

\addplot [color=black, line width=1.0pt]
  table[row sep=crcr]{%
0	51.8020213564446\\
0.157894736842105	51.0734582200083\\
0.315789473684211	50.3489120281773\\
0.473684210526316	49.6284612129186\\
0.631578947368421	48.9121547457252\\
0.789473684210526	48.200055742311\\
0.947368421052632	47.4922241448457\\
1.10526315789474	46.7887240767804\\
1.26315789473684	46.089584167414\\
1.42105263157895	45.3948795113323\\
1.57894736842105	44.7046644138314\\
1.73684210526316	44.0189797262804\\
1.89473684210526	43.3378887973627\\
2.05263157894737	42.6614288377107\\
2.21052631578947	41.9896681082348\\
2.36842105263158	41.3226104227994\\
2.52631578947368	40.6603392797206\\
2.68421052631579	40.0028893148199\\
2.84210526315789	39.3502966969581\\
3	38.7026027832065\\
};
\addlegendentry{Trip-based}

\end{axis}
\end{tikzpicture}%
\vspace*{-0.3in}
\caption{Trips from $\mathcal{C}$ to $\mathcal{C}$, $\lambda_{11}$,  under different surcharge $\bar{p}_{ij}$.}
\label{figure7_2z}
\end{minipage}
\begin{minipage}[b]{0.005\linewidth}
\hfill
\end{minipage}
\begin{minipage}[b]{0.32\linewidth}
\centering
%
%
\definecolor{mycolor1}{rgb}{0.00000,0.44700,0.74100}%
\definecolor{mycolor2}{rgb}{0.85000,0.32500,0.09800}%
\begin{tikzpicture}

\begin{axis}[%
width=1.694in,
height=1.03in,
at={(1.358in,0.0in)},
scale only axis,
xmin=0,
xmax=3,
xlabel style={font=\color{white!15!black}},
xlabel={Congestion charge},
ymin=14,
ymax=20,
ylabel style={font=\color{white!15!black}},
ylabel={$\lambda_{12}$},
axis background/.style={fill=white},
legend style={at={(0.0,0.0)}, anchor=south west, legend cell align=left, align=left, font=\small, draw=white!15!black, nodes={scale=0.8, transform shape} }
]
\addplot [color=mycolor1, line width=1.0pt]
  table[row sep=crcr]{%
0	19.4906972927633\\
0.157894736842105	19.3620046409013\\
0.315789473684211	19.233651042641\\
0.473684210526316	19.105615507917\\
0.631578947368421	18.9778975052205\\
0.789473684210526	18.8504662391642\\
0.947368421052632	18.7233562281633\\
1.10526315789474	18.596522514531\\
1.26315789473684	18.4699708215923\\
1.42105263157895	18.3436889896087\\
1.57894736842105	18.2176710789104\\
1.73684210526316	18.0919086242336\\
1.89473684210526	17.9664009608705\\
2.05263157894737	17.841135032766\\
2.21052631578947	17.7161097965312\\
2.36842105263158	17.5913176465561\\
2.52631578947368	17.4667482693707\\
2.68421052631579	17.3424048877126\\
2.84210526315789	17.2182802483838\\
3	17.094370444111\\
};
\addlegendentry{Cordon 1}

\addplot [color=mycolor2, line width=1.0pt]
  table[row sep=crcr]{%
0	19.4906972927633\\
0.157894736842105	19.1803712437095\\
0.315789473684211	18.8750058953144\\
0.473684210526316	18.57446281187\\
0.631578947368421	18.2786078087548\\
0.789473684210526	17.9873171235158\\
0.947368421052632	17.7004687149492\\
1.10526315789474	17.4179459543185\\
1.26315789473684	17.1396412628142\\
1.42105263157895	16.8654517236976\\
1.57894736842105	16.5952777721686\\
1.73684210526316	16.3290272174758\\
1.89473684210526	16.0666091003168\\
2.05263157894737	15.8079426142768\\
2.21052631578947	15.5529490402803\\
2.36842105263158	15.3015489409935\\
2.52631578947368	15.0536728120383\\
2.68421052631579	14.8092547379941\\
2.84210526315789	14.5682400880958\\
3	14.33055040357\\
};
\addlegendentry{Cordon 2}

\addplot [color=black, line width=1.0pt]
  table[row sep=crcr]{%
0	19.4906972927633\\
0.157894736842105	19.2641824128883\\
0.315789473684211	19.0381996395251\\
0.473684210526316	18.8127701625777\\
0.631578947368421	18.5879122690094\\
0.789473684210526	18.3636446053938\\
0.947368421052632	18.1399862227741\\
1.10526315789474	17.9169579439235\\
1.26315789473684	17.6945644697796\\
1.42105263157895	17.4728351227775\\
1.57894736842105	17.2517899571804\\
1.73684210526316	17.0314377751918\\
1.89473684210526	16.8118057706028\\
2.05263157894737	16.5929038495663\\
2.21052631578947	16.3747566667549\\
2.36842105263158	16.1573717500418\\
2.52631578947368	15.9407709086347\\
2.68421052631579	15.7249748897493\\
2.84210526315789	15.5099955332208\\
3	15.2958511991009\\
};
\addlegendentry{Trip-based}

\end{axis}
\end{tikzpicture}%
\vspace*{-0.3in}
\caption{Trips from $\mathcal{C}$ to $\mathcal{R}$, $\lambda_{12}$,  under different surcharge $\bar{p}_{ij}$.}
\label{figure8_2z}
\end{minipage}
\begin{minipage}[b]{0.005\linewidth}
\hfill
\end{minipage}
\begin{minipage}[b]{0.32\linewidth}
\centering
%
%
\definecolor{mycolor1}{rgb}{0.00000,0.44700,0.74100}%
\definecolor{mycolor2}{rgb}{0.85000,0.32500,0.09800}%
\begin{tikzpicture}

\begin{axis}[%
width=1.694in,
height=1.03in,
at={(1.358in,0.0in)},
scale only axis,
xmin=0,
xmax=3,
xlabel style={font=\color{white!15!black}},
xlabel={Congestion charge},
ymin=20,
ymax=28,
ylabel style={font=\color{white!15!black}},
ylabel={$\lambda_{21}$},
axis background/.style={fill=white},
legend style={at={(0.0,0.0)}, anchor=south west, legend cell align=left, align=left, font=\small, draw=white!15!black, nodes={scale=0.8, transform shape} }
]
\addplot [color=mycolor1, line width=1.0pt]
  table[row sep=crcr]{%
0	26.799255369302\\
0.157894736842105	26.5523984563374\\
0.315789473684211	26.307104042281\\
0.473684210526316	26.0633702870319\\
0.631578947368421	25.8212050712527\\
0.789473684210526	25.5805931601194\\
0.947368421052632	25.3415659742403\\
1.10526315789474	25.1040962259507\\
1.26315789473684	24.8681945943961\\
1.42105263157895	24.6338576260532\\
1.57894736842105	24.4010745765137\\
1.73684210526316	24.1698512754779\\
1.89473684210526	23.9401826419425\\
2.05263157894737	23.7120561195807\\
2.21052631578947	23.4854759261775\\
2.36842105263158	23.2604304272775\\
2.52631578947368	23.0369146502698\\
2.68421052631579	22.814913420043\\
2.84210526315789	22.5944287131362\\
3	22.3754488519291\\
};
\addlegendentry{Cordon 1}

\addplot [color=mycolor2, line width=1.0pt]
  table[row sep=crcr]{%
0	26.799255369302\\
0.157894736842105	26.3729103408643\\
0.315789473684211	25.9514524112313\\
0.473684210526316	25.5348630441213\\
0.631578947368421	25.1231286930379\\
0.789473684210526	24.7162253319957\\
0.947368421052632	24.3141387567881\\
1.10526315789474	23.9168370762871\\
1.26315789473684	23.5243011036477\\
1.42105263157895	23.1365057279757\\
1.57894736842105	22.7534230090254\\
1.73684210526316	22.3750235768728\\
1.89473684210526	22.0012786471705\\
2.05263157894737	21.6321632911296\\
2.21052631578947	21.2676463191782\\
2.36842105263158	20.9076947020178\\
2.52631578947368	20.5522783510897\\
2.68421052631579	20.2013700631672\\
2.84210526315789	19.8549526991099\\
3	19.5129542557729\\
};
\addlegendentry{Cordon 2}

\addplot [color=black, line width=1.0pt]
  table[row sep=crcr]{%
0	26.799255369302\\
0.157894736842105	26.5972592143034\\
0.315789473684211	26.394679490625\\
0.473684210526316	26.1915403973759\\
0.631578947368421	25.9878408530283\\
0.789473684210526	25.7836049829266\\
0.947368421052632	25.5788388293915\\
1.10526315789474	25.3735602414435\\
1.26315789473684	25.1677681554389\\
1.42105263157895	24.9614910176621\\
1.57894736842105	24.7547423315486\\
1.73684210526316	24.5475176697022\\
1.89473684210526	24.339851923955\\
2.05263157894737	24.1317504008119\\
2.21052631578947	23.9232349363441\\
2.36842105263158	23.7142983666959\\
2.52631578947368	23.5049773592396\\
2.68421052631579	23.2952819249596\\
2.84210526315789	23.085222384615\\
3	22.8748112644201\\
};
\addlegendentry{Trip-based}

\end{axis}
\end{tikzpicture}%
\vspace*{-0.3in}
\caption{Trips from $\mathcal{R}$ to $\mathcal{C}$, $\lambda_{21}$,  under different surcharge $\bar{p}_{ij}$.}
\label{figure9_2z}
\end{minipage}
\begin{minipage}[b]{0.32\linewidth}
\centering
%
%
\definecolor{mycolor1}{rgb}{0.00000,0.44700,0.74100}%
\definecolor{mycolor2}{rgb}{0.85000,0.32500,0.09800}%
\begin{tikzpicture}

\begin{axis}[%
width=1.694in,
height=1.03in,
at={(1.358in,0.0in)},
scale only axis,
xmin=0,
xmax=3,
xlabel style={font=\color{white!15!black}},
xlabel={Congestion charge},
ymin=45,
ymax=70,
ylabel style={font=\color{white!15!black}},
ylabel={$\lambda_{22}$},
axis background/.style={fill=white},
legend style={at={(0.0,0.54)}, anchor=south west, legend cell align=left, align=left, font=\small, draw=white!15!black, nodes={scale=0.8, transform shape} }
]
\addplot [color=mycolor1, line width=1.0pt]
  table[row sep=crcr]{%
0	58.5373681759301\\
0.157894736842105	58.6389711202651\\
0.315789473684211	58.7406193745077\\
0.473684210526316	58.8423080536616\\
0.631578947368421	58.9440376961458\\
0.789473684210526	59.0457698592082\\
0.947368421052632	59.1475370016456\\
1.10526315789474	59.249305161649\\
1.26315789473684	59.3510597647947\\
1.42105263157895	59.4528139552078\\
1.57894736842105	59.5545173106246\\
1.73684210526316	59.6561963248957\\
1.89473684210526	59.7578171699758\\
2.05263157894737	59.8593561661616\\
2.21052631578947	59.9608137760196\\
2.36842105263158	60.0621915307725\\
2.52631578947368	60.1634306819543\\
2.68421052631579	60.2645378782795\\
2.84210526315789	60.3655037916252\\
3	60.4662903184067\\
};
\addlegendentry{Cordon 1}

\addplot [color=mycolor2, line width=1.0pt]
  table[row sep=crcr]{%
0	58.5373681759301\\
0.157894736842105	58.3990147997098\\
0.315789473684211	58.2609057807309\\
0.473684210526316	58.1230280577868\\
0.631578947368421	57.9853914942741\\
0.789473684210526	57.8479805620486\\
0.947368421052632	57.7107957053928\\
1.10526315789474	57.5738187981769\\
1.26315789473684	57.4370339593067\\
1.42105263157895	57.3004477045819\\
1.57894736842105	57.1640328818158\\
1.73684210526316	57.0277808625824\\
1.89473684210526	56.8916775178091\\
2.05263157894737	56.7557240399313\\
2.21052631578947	56.6198942940214\\
2.36842105263158	56.484182457611\\
2.52631578947368	56.3485683428888\\
2.68421052631579	56.2130417347294\\
2.84210526315789	56.0776362424327\\
3	55.9422445544458\\
};
\addlegendentry{Cordon 2}

\addplot [color=black, line width=1.0pt]
  table[row sep=crcr]{%
0	58.5373681759301\\
0.157894736842105	57.96123582357\\
0.315789473684211	57.3860850889104\\
0.473684210526316	56.8119992791425\\
0.631578947368421	56.2390038888986\\
0.789473684210526	55.6671514867307\\
0.947368421052632	55.0964959245675\\
1.10526315789474	54.5270992693238\\
1.26315789473684	53.9589736751957\\
1.42105263157895	53.3921911935091\\
1.57894736842105	52.8267992797663\\
1.73684210526316	52.2628464973479\\
1.89473684210526	51.70036657767\\
2.05263157894737	51.1394131576031\\
2.21052631578947	50.5800573428069\\
2.36842105263158	50.0222704720594\\
2.52631578947368	49.4661768374931\\
2.68421052631579	48.9117787771274\\
2.84210526315789	48.3591262414903\\
3	47.8082437719339\\
};
\addlegendentry{Trip-based}

\end{axis}
\end{tikzpicture}%
\vspace*{-0.3in}
\caption{Trips from $\mathcal{R}$ to $\mathcal{R}$, $\lambda_{22}$,  under different surcharge $\bar{p}_{ij}$.} 
\label{figure10_2z}
\end{minipage}
\begin{minipage}[b]{0.005\linewidth}
\hfill
\end{minipage}
\begin{minipage}[b]{0.32\linewidth}
\centering
%
%
\definecolor{mycolor1}{rgb}{0.00000,0.44700,0.74100}%
\definecolor{mycolor2}{rgb}{0.85000,0.32500,0.09800}%
\begin{tikzpicture}

\pgfplotsset{every axis y label/.style={
at={(-0.5,0.5)},
xshift=30pt,
rotate=90}}

\begin{axis}[%
width=1.694in,
height=1.03in,
at={(1.358in,0.0in)},
scale only axis,
xmin=0,
xmax=3,
xlabel style={font=\color{white!15!black}},
xlabel={Congestion charge},
ymin=1300,
ymax=1650,
ytick={1300,1400,  1500, 1600},
yticklabels={{1.3K},{1.4K},{1.5K}, {1.6K}},
ylabel style={font=\color{white!15!black}},
ylabel={All vehicles in $\mathcal{C}$},
axis background/.style={fill=white},
legend style={at={(0.0,0.0)}, anchor=south west, legend cell align=left, align=left, font=\small, draw=white!15!black, nodes={scale=0.8, transform shape} }
]
\addplot [color=mycolor1, line width=1.0pt]
  table[row sep=crcr]{%
0	1620.27020134007\\
0.157894736842105	1611.48440937909\\
0.315789473684211	1602.71342918838\\
0.473684210526316	1593.95712795177\\
0.631578947368421	1585.21589848813\\
0.789473684210526	1576.4886700668\\
0.947368421052632	1567.77728798781\\
1.10526315789474	1559.08039323546\\
1.26315789473684	1550.39830928906\\
1.42105263157895	1541.73138384489\\
1.57894736842105	1533.07890769252\\
1.73684210526316	1524.44159578162\\
1.89473684210526	1515.81922659771\\
2.05263157894737	1507.21134543315\\
2.21052631578947	1498.61843474002\\
2.36842105263158	1490.04066513767\\
2.52631578947368	1481.47716342953\\
2.68421052631579	1472.92839257281\\
2.84210526315789	1464.39439167064\\
3	1455.87485057663\\
};
\addlegendentry{Cordon 1}

\addplot [color=mycolor2, line width=1.0pt]
  table[row sep=crcr]{%
0	1620.27020134007\\
0.157894736842105	1611.0949372732\\
0.315789473684211	1602.0833541192\\
0.473684210526316	1593.23134016565\\
0.631578947368421	1584.53555007379\\
0.789473684210526	1575.9921552407\\
0.947368421052632	1567.59783647644\\
1.10526315789474	1559.34882419256\\
1.26315789473684	1551.24170570065\\
1.42105263157895	1543.27346293464\\
1.57894736842105	1535.44058874439\\
1.73684210526316	1527.74000995291\\
1.89473684210526	1520.16854779795\\
2.05263157894737	1512.72354440375\\
2.21052631578947	1505.40192621071\\
2.36842105263158	1498.20085632455\\
2.52631578947368	1491.11751219111\\
2.68421052631579	1484.14931115924\\
2.84210526315789	1477.29479362436\\
3	1470.54901643243\\
};
\addlegendentry{Cordon 2}

\addplot [color=black, line width=1.0pt]
  table[row sep=crcr]{%
0	1620.27020134007\\
0.157894736842105	1605.08104977503\\
0.315789473684211	1589.91038134612\\
0.473684210526316	1574.7600218383\\
0.631578947368421	1559.63025601002\\
0.789473684210526	1544.52232651306\\
0.947368421052632	1529.4371979962\\
1.10526315789474	1514.37619564189\\
1.26315789473684	1499.33914199581\\
1.42105263157895	1484.32789439682\\
1.57894736842105	1469.34322422993\\
1.73684210526316	1454.3857706482\\
1.89473684210526	1439.4567875227\\
2.05263157894737	1424.55681454181\\
2.21052631578947	1409.68759805723\\
2.36842105263158	1394.84789682476\\
2.52631578947368	1380.04065049525\\
2.68421052631579	1365.26592213333\\
2.84210526315789	1350.52441202049\\
3	1335.8169214204\\
};
\addlegendentry{Trip-based}

\end{axis}
\end{tikzpicture}%
\vspace*{-0.3in}
\caption{Total number of drivers in $\mathcal{C}$ under different surcharge $\bar{p}_{ij}$.}
\label{figure11_2z}
\end{minipage}
\begin{minipage}[b]{0.005\linewidth}
\hfill
\end{minipage}
\begin{minipage}[b]{0.32\linewidth}
\centering
%
%
\definecolor{mycolor1}{rgb}{0.00000,0.44700,0.74100}%
\definecolor{mycolor2}{rgb}{0.85000,0.32500,0.09800}%
\begin{tikzpicture}
\pgfplotsset{every axis y label/.style={
at={(-0.5,0.5)},
xshift=30pt,
rotate=90}}

\begin{axis}[%
width=1.694in,
height=1.03in,
at={(1.358in,0.0in)},
scale only axis,
xmin=0,
xmax=3,
xlabel style={font=\color{white!15!black}},
xlabel={Congestion charge},
ymin=1800,
ymax=2100,
ytick={1800,1900, 2000, 2100},
yticklabels={{1.6K},{1.9K},{2.0K},{2.1K}},
ylabel style={font=\color{white!15!black}},
ylabel={All vehicles in $\mathcal{R}$ },
axis background/.style={fill=white},
legend style={at={(0.0,0.0)}, anchor=south west, legend cell align=left, align=left, font=\small, draw=white!15!black, nodes={scale=0.8, transform shape} }
]
\addplot [color=mycolor1, line width=1.0pt]
  table[row sep=crcr]{%
0	2067.08581258234\\
0.157894736842105	2065.39128620044\\
0.315789473684211	2063.72089200264\\
0.473684210526316	2062.07471078977\\
0.631578947368421	2060.45262120228\\
0.789473684210526	2058.85301907349\\
0.947368421052632	2057.27696380034\\
1.10526315789474	2055.72356472791\\
1.26315789473684	2054.19171762811\\
1.42105263157895	2052.68251997936\\
1.57894736842105	2051.19353118999\\
1.73684210526316	2049.72607956607\\
1.89473684210526	2048.27867085431\\
2.05263157894737	2046.85049081913\\
2.21052631578947	2045.44148135688\\
2.36842105263158	2044.05236113047\\
2.52631578947368	2042.67996023866\\
2.68421052631579	2041.32555670073\\
2.84210526315789	2039.98846813434\\
3	2038.66711098681\\
};
\addlegendentry{Cordon 1}

\addplot [color=mycolor2, line width=1.0pt]
  table[row sep=crcr]{%
0	2067.08581258234\\
0.157894736842105	2057.19210066731\\
0.315789473684211	2047.37420119381\\
0.473684210526316	2037.63140036623\\
0.631578947368421	2027.96396836136\\
0.789473684210526	2018.37123238021\\
0.947368421052632	2008.85278235447\\
1.10526315789474	1999.40789971897\\
1.26315789473684	1990.03542060937\\
1.42105263157895	1980.73580556034\\
1.57894736842105	1971.50728885873\\
1.73684210526316	1962.34968460582\\
1.89473684210526	1953.26211152297\\
2.05263157894737	1944.24453919718\\
2.21052631578947	1935.29539459552\\
2.36842105263158	1926.41464549117\\
2.52631578947368	1917.6010726143\\
2.68421052631579	1908.85396299459\\
2.84210526315789	1900.17443956161\\
3	1891.55872961268\\
};
\addlegendentry{Cordon 2}

\addplot [color=black, line width=1.0pt]
  table[row sep=crcr]{%
0	2067.08581258234\\
0.157894736842105	2055.23479894616\\
0.315789473684211	2043.33419079601\\
0.473684210526316	2031.38602621598\\
0.631578947368421	2019.39014145259\\
0.789473684210526	2007.34681713927\\
0.947368421052632	1995.25729493343\\
1.10526315789474	1983.12280269533\\
1.26315789473684	1970.94239618545\\
1.42105263157895	1958.71756575844\\
1.57894736842105	1946.44904348837\\
1.73684210526316	1934.1385626965\\
1.89473684210526	1921.78489748028\\
2.05263157894737	1909.38971455736\\
2.21052631578947	1896.95507538571\\
2.36842105263158	1884.47795932369\\
2.52631578947368	1871.96323376343\\
2.68421052631579	1859.40920534848\\
2.84210526315789	1846.81728490754\\
3	1834.18738451779\\
};
\addlegendentry{Trip-based}

\end{axis}
\end{tikzpicture}%
\vspace*{-0.3in}
\caption{Total number of drivers in $\mathcal{R}$  under different surcharge $\bar{p}_{ij}$.}
\label{figure12_2z}
\end{minipage}
\end{figure*}

\begin{figure*}%
\begin{minipage}[b]{0.32\linewidth}
\centering
%
%
\definecolor{mycolor1}{rgb}{0.00000,0.44700,0.74100}%
\definecolor{mycolor2}{rgb}{0.85000,0.32500,0.09800}%
\begin{tikzpicture}

\pgfplotsset{every axis y label/.style={
at={(-0.5,0.5)},
xshift=30pt,
rotate=90}}

\begin{axis}[%
width=1.694in,
height=1.03in,
at={(1.358in,0.0in)},
scale only axis,
xmin=0,
xmax=3,
xlabel style={font=\color{white!15!black}},
xlabel={Congestion charge},
ymin=0.40,
ymax=0.435,
ylabel style={font=\color{white!15!black}},
ylabel={Occupancy in $\mathcal{C}$},
axis background/.style={fill=white},
legend style={at={(0.0,0.0)}, anchor=south west, legend cell align=left, align=left, font=\small, draw=white!15!black, nodes={scale=0.8, transform shape} }
]
\addplot [color=mycolor1, line width=1.0pt]
  table[row sep=crcr]{%
0	0.429745306195207\\
0.157894736842105	0.429066574926187\\
0.315789473684211	0.428387003152623\\
0.473684210526316	0.427706331067236\\
0.631578947368421	0.427024492058674\\
0.789473684210526	0.426341270080349\\
0.947368421052632	0.42565681488165\\
1.10526315789474	0.424970726828522\\
1.26315789473684	0.424283039147917\\
1.42105263157895	0.423593487620166\\
1.57894736842105	0.422902088692206\\
1.73684210526316	0.422208569761109\\
1.89473684210526	0.421512948532616\\
2.05263157894737	0.420815046366856\\
2.21052631578947	0.420114760670127\\
2.36842105263158	0.419411845099606\\
2.52631578947368	0.418706339763229\\
2.68421052631579	0.417998045483728\\
2.84210526315789	0.417286845819196\\
3	0.416572672816546\\
};
\addlegendentry{Cordon 1}

\addplot [color=mycolor2, line width=1.0pt]
  table[row sep=crcr]{%
0	0.429745306195207\\
0.157894736842105	0.428513519626948\\
0.315789473684211	0.427285150731791\\
0.473684210526316	0.42605994780823\\
0.631578947368421	0.424837534367594\\
0.789473684210526	0.42361768693551\\
0.947368421052632	0.422400106407356\\
1.10526315789474	0.421184527388587\\
1.26315789473684	0.419970783600852\\
1.42105263157895	0.418758598250779\\
1.57894736842105	0.417547812715716\\
1.73684210526316	0.416338244784146\\
1.89473684210526	0.415129726480916\\
2.05263157894737	0.413922098450344\\
2.21052631578947	0.412715290620888\\
2.36842105263158	0.411509121426455\\
2.52631578947368	0.410303550182617\\
2.68421052631579	0.409098546695426\\
2.84210526315789	0.407894039228849\\
3	0.40668999197576\\
};
\addlegendentry{Cordon 2}

\addplot [color=black, line width=1.0pt]
  table[row sep=crcr]{%
0	0.429745306195207\\
0.157894736842105	0.428958921724651\\
0.315789473684211	0.428165442274722\\
0.473684210526316	0.427364699194578\\
0.631578947368421	0.426556657674973\\
0.789473684210526	0.425741227252207\\
0.947368421052632	0.424918272594999\\
1.10526315789474	0.42408764833337\\
1.26315789473684	0.423249239009907\\
1.42105263157895	0.422402973763368\\
1.57894736842105	0.421548765266605\\
1.73684210526316	0.420686327884137\\
1.89473684210526	0.419815738144848\\
2.05263157894737	0.418936761284288\\
2.21052631578947	0.418049196434372\\
2.36842105263158	0.417153171079884\\
2.52631578947368	0.416248211819119\\
2.68421052631579	0.415334409908422\\
2.84210526315789	0.414411540381005\\
3	0.413479463837341\\
};
\addlegendentry{Trip-based}

\end{axis}

\end{tikzpicture}%
\vspace*{-0.3in}
\caption{Vehicle occupancy\protect\footnotemark in $\mathcal{R}$  under different surcharge $\bar{p}_{ij}$.} 
\label{occupancy1_2z}
\end{minipage}
\begin{minipage}[b]{0.005\linewidth}
\hfill
\end{minipage}
\begin{minipage}[b]{0.32\linewidth}
\centering
%
%
\definecolor{mycolor1}{rgb}{0.00000,0.44700,0.74100}%
\definecolor{mycolor2}{rgb}{0.85000,0.32500,0.09800}%
\begin{tikzpicture}

\pgfplotsset{every axis y label/.style={
at={(-0.5,0.5)},
xshift=30pt,
rotate=90}}

\begin{axis}[%
width=1.694in,
height=1.03in,
at={(1.358in,0.0in)},
scale only axis,
xmin=0,
xmax=3,
xlabel style={font=\color{white!15!black}},
xlabel={Congestion charge},
ymin=0.30,
ymax=0.33,
ylabel style={font=\color{white!15!black}},
ylabel={Occupancy in $\mathcal{R}$},
axis background/.style={fill=white},
legend style={at={(0.0,0.0)}, anchor=south west, legend cell align=left, align=left, font=\small, draw=white!15!black, nodes={scale=0.8, transform shape} }
]
\addplot [color=mycolor1, line width=1.0pt]
  table[row sep=crcr]{%
0	0.317757986061337\\
0.157894736842105	0.318045841978785\\
0.315789473684211	0.318331515544718\\
0.473684210526316	0.318614924359771\\
0.631578947368421	0.318896056236981\\
0.789473684210526	0.319174790696347\\
0.947368421052632	0.31945136679517\\
1.10526315789474	0.3197254669554\\
1.26315789473684	0.3199972235285\\
1.42105263157895	0.320266494090063\\
1.57894736842105	0.320533329026531\\
1.73684210526316	0.320797632127891\\
1.89473684210526	0.321059438014011\\
2.05263157894737	0.321318696304596\\
2.21052631578947	0.321575394600419\\
2.36842105263158	0.321829379231999\\
2.52631578947368	0.32208081086716\\
2.68421052631579	0.322329482888679\\
2.84210526315789	0.322575431273888\\
3	0.322818675065984\\
};
\addlegendentry{Cordon 1}

\addplot [color=mycolor2, line width=1.0pt]
  table[row sep=crcr]{%
0	0.317757986061337\\
0.157894736842105	0.317748048719561\\
0.315789473684211	0.317733812649369\\
0.473684210526316	0.317715100750194\\
0.631578947368421	0.317691743658041\\
0.789473684210526	0.317663600927155\\
0.947368421052632	0.317630569085537\\
1.10526315789474	0.317592496542166\\
1.26315789473684	0.317549361487726\\
1.42105263157895	0.317501024884044\\
1.57894736842105	0.317447464562729\\
1.73684210526316	0.317388597259213\\
1.89473684210526	0.317324391292988\\
2.05263157894737	0.317254789503278\\
2.21052631578947	0.317179827319049\\
2.36842105263158	0.317099419781305\\
2.52631578947368	0.317013602565392\\
2.68421052631579	0.316922411426244\\
2.84210526315789	0.316825831495\\
3	0.316723799353635\\
};
\addlegendentry{Cordon 2}

\addplot [color=black, line width=1.0pt]
  table[row sep=crcr]{%
0	0.317757986061337\\
0.157894736842105	0.3173958371632\\
0.315789473684211	0.317026627736019\\
0.473684210526316	0.316650247264717\\
0.631578947368421	0.316266637020429\\
0.789473684210526	0.315875826382795\\
0.947368421052632	0.315477673714685\\
1.10526315789474	0.315072106797223\\
1.26315789473684	0.314659080542546\\
1.42105263157895	0.314238552469597\\
1.57894736842105	0.313810501006983\\
1.73684210526316	0.313374628784635\\
1.89473684210526	0.312931134136342\\
2.05263157894737	0.312479831937079\\
2.21052631578947	0.312020638537099\\
2.36842105263158	0.31155350641503\\
2.52631578947368	0.311078288839483\\
2.68421052631579	0.310595018528313\\
2.84210526315789	0.3101035463062\\
3	0.309603821336183\\
};
\addlegendentry{Trip-based}

\end{axis}

\end{tikzpicture}%
\vspace*{-0.3in}
\caption{Vehicle occupancy in $\mathcal{C}$  under different surcharge $\bar{p}_{ij}$.}
\label{occupancy2_2z}
\end{minipage}
\begin{minipage}[b]{0.005\linewidth}
\hfill
\end{minipage}
\begin{minipage}[b]{0.32\linewidth}
\centering
%
%
\definecolor{mycolor1}{rgb}{0.00000,0.44700,0.74100}%
\definecolor{mycolor2}{rgb}{0.85000,0.32500,0.09800}%
\begin{tikzpicture}

\pgfplotsset{every axis y label/.style={
at={(-0.5,0.5)},
xshift=30pt,
rotate=90}}

\begin{axis}[%
width=1.694in,
height=1.03in,
at={(1.358in,0.0in)},
scale only axis,
xmin=0,
xmax=3,
xlabel style={font=\color{white!15!black}},
xlabel={Congestion charge},
ymin=45000,
ymax=65000,
ylabel style={font=\color{white!15!black}},
ylabel={TNC profit},
axis background/.style={fill=white},
legend style={at={(0.0,0.0)}, anchor=south west, legend cell align=left, align=left, font=\small, draw=white!15!black, nodes={scale=0.8, transform shape} }
]
\addplot [color=mycolor1, line width=1.0pt]
  table[row sep=crcr]{%
0	64748.6208855417\\
0.157894736842105	64355.2683610763\\
0.315789473684211	63964.640128053\\
0.473684210526316	63576.6853943822\\
0.631578947368421	63191.3536698617\\
0.789473684210526	62808.5947987218\\
0.947368421052632	62428.3589926274\\
1.10526315789474	62050.5968616058\\
1.26315789473684	61675.2594469892\\
1.42105263157895	61302.298252092\\
1.57894736842105	60931.6652731093\\
1.73684210526316	60563.3130293812\\
1.89473684210526	60197.1945929527\\
2.05263157894737	59833.2636175226\\
2.21052631578947	59471.4743665815\\
2.36842105263158	59111.7817408195\\
2.52631578947368	58754.1413046911\\
2.68421052631579	58398.509311949\\
2.84210526315789	58044.8427304451\\
3	57693.0992658287\\
};
\addlegendentry{Cordon 1}

\addplot [color=mycolor2, line width=1.0pt]
  table[row sep=crcr]{%
0	64748.6208855417\\
0.157894736842105	64012.3044270102\\
0.315789473684211	63292.217681577\\
0.473684210526316	62587.9141671576\\
0.631578947368421	61898.9606464216\\
0.789473684210526	61224.936938555\\
0.947368421052632	60565.4357157181\\
1.10526315789474	59920.0622854566\\
1.26315789473684	59288.4343601345\\
1.42105263157895	58670.1818145292\\
1.57894736842105	58064.9464328413\\
1.73684210526316	57472.3816460539\\
1.89473684210526	56892.152261073\\
2.05263157894737	56323.934182549\\
2.21052631578947	55767.4141286295\\
2.36842105263158	55222.2893416737\\
2.52631578947368	54688.2672950665\\
2.68421052631579	54165.0653971\\
2.84210526315789	53652.410693002\\
3	53150.0395660581\\
};
\addlegendentry{Cordon 2}

\addplot [color=black, line width=1.0pt]
  table[row sep=crcr]{%
0	64748.6208855417\\
0.157894736842105	63633.788448915\\
0.315789473684211	62530.3835630666\\
0.473684210526316	61438.400692786\\
0.631578947368421	60357.8329074349\\
0.789473684210526	59288.6718839794\\
0.947368421052632	58230.9079106586\\
1.10526315789474	57184.5298914602\\
1.26315789473684	56149.5253513415\\
1.42105263157895	55125.8804419384\\
1.57894736842105	54113.5799482936\\
1.73684210526316	53112.6072961353\\
1.89473684210526	52122.9445596274\\
2.05263157894737	51144.5724701343\\
2.21052631578947	50177.4704253703\\
2.36842105263158	49221.6164996752\\
2.52631578947368	48276.9874537166\\
2.68421052631579	47343.5587464813\\
2.84210526315789	46421.3045465792\\
3	45510.1977445475\\
};
\addlegendentry{Trip-based}

\end{axis}

\end{tikzpicture}%
\vspace*{-0.3in}
\caption{Platform profit per hour under different surcharge $\bar{p}_{ij}$.}
\label{profit_2z}
\end{minipage}
\end{figure*}
\footnotetext{Vehicle occupancy is calculated by dividing the number of occupied vehicles by the total number of vehicles (including all vehicles that are either occupied, idle, or on the way to pick up passengers.)}

\subsubsection{Comparing distinct schemes of congestion charges}
{
The simulation results in Section \ref{casestudy_part1} provide a snapshot of the ride-sourcing network at two specific values of $\bar{p}_{ij}$. To better compare distinct congestion charge schemes, here we will continuously adjust $\bar{p}_{ij}$  and calculate the aggregated market outcomes in the congestion area and the remote area for the three congestion charge schemes, respectively. 
}

Figures \ref{figure1_2z}-\ref{figure12_2z} show the simulation results as we vary the congestion charge $\bar{p}_{ij}$ from $\$0$/vehicle (or trip) to $\$3$/vehicle (or trip).  For notation convenience, we denote $\lambda_{11}$ as the total trip flow from $\mathcal{C}$ to $\mathcal{C}$, $\lambda_{12}$ as the total trip flow from $\mathcal{C}$ to $\mathcal{R}$, $\lambda_{21}$ as the total trip flow from $\mathcal{R}$ to $\mathcal{C}$, and $\lambda_{22}$ as the total trip flow from $\mathcal{R}$ to $\mathcal{R}$.  In the figure legend, ``Cordon 1'' represents the one-directional cordon price, ``Cordon 2'' represents the bi-directional cordon price, and ``Trip-based'' represents the trip-based congestion charge.   In all cases, the number of idle vehicles and the total number of  vehicles (including idle vehicles, occupied vehicles,\footnote{For occupied vehicles traveling from zone i to zone j, we evenly distribute the vehicle time to all zones along the shortest path from i to j.} and vehicles on the way to pick up passengers) in the congestion area is reduced (Figure \ref{figure1_2z} and Figure \ref{figure11_2z}). Passengers in $\mathcal{C}$ bear the congestion charge and shift to other transport modes due to the increased overall travel cost associated with ride-sourcing trips (Figure \ref{figure7_2z} and Figure \ref{figure8_2z}). This further confirms that all forms of congestion charge are effective in reducing the number of ride-sourcing vehicles and alleviating  traffic congestion in $\mathcal{C}$.  

{
Based on the simulation results, we note that the one-directional cordon price not only limits congestion in the congestion area but also benefits  passengers in the remote area. After the one-directional cordon price is imposed,  the average ride fare of the remote area reduces\footnote{Note that the ride fare $r_i$ shown in Figure \ref{figure2_2z} and Figure \ref{figure5_2z} do not include the cordon surcharge for trips entering the congestion area.} (Figure \ref{figure5_2z}), and the average waiting time also decreases (Figure \ref{figure6_2z}).\footnote{Note that in  Figure \ref{figure4_2z}, the total number of idle drivers in the remote area decreases. This does not conflict with the reduced average waiting time in Figure \ref{figure6_2z} because the waiting time in certain zones of the remote area reduces, while the waiting time in other zones increases. Its average value may either increase or decrease, depending on the model parameters.}  Therefore,  the number of trips from the remote area to the remote area increases (Figure \ref{figure10_2z}), and the surplus of passengers in these underserved zones is improved.  For this reason, the one-directional cordon price is distinguished from the other congestion charge schemes as it can reduce congestion and promote equity at the same time.
}

\begin{figure*}%
\begin{minipage}[b]{0.32\linewidth}
\centering
%
%
\definecolor{mycolor1}{rgb}{0.00000,0.44700,0.74100}%
\definecolor{mycolor2}{rgb}{0.85000,0.32500,0.09800}%
\begin{tikzpicture}

\begin{axis}[%
width=1.694in,
height=1.03in,
at={(1.358in,0.0in)},
scale only axis,
xmin=0,
xmax=150,
xlabel style={font=\color{white!15!black}},
xlabel={$\Delta N_\mathcal{C}$},
ymin=135,
ymax=160,
ylabel style={font=\color{white!15!black}},
ylabel={Passenger demand},
axis background/.style={fill=white},
legend style={at={(0.0,0.0)}, anchor=south west, legend cell align=left, align=left, font=\small, draw=white!15!black, nodes={scale=0.8, transform shape} }
]
\addplot [color=mycolor1, line width=1.0pt]
  table[row sep=crcr]{%
0	156.62934219444\\
7.8800623635603	156.15113096569\\
15.7601247271206	155.674379552742\\
23.6401870906809	155.199051138453\\
31.5202494542412	154.725149150534\\
39.4003118178015	154.252603793\\
47.2803741813618	153.781452441301\\
55.1604365449221	153.311690579426\\
63.0404989084824	152.843243958298\\
70.9205612720427	152.376117523755\\
78.800623635603	151.91029829901\\
86.6806859991633	151.445725695624\\
94.5607483627236	150.982555633283\\
102.440810726284	150.520624115394\\
110.320873089844	150.059875038048\\
118.200935453404	149.600309584173\\
126.080997816965	149.141918314102\\
133.961060180525	148.684649794948\\
141.841122544085	148.228471818759\\
149.721184907646	147.773382720705\\
};
\addlegendentry{Cordon 1}

\addplot [color=mycolor2, line width=1.0pt]
  table[row sep=crcr]{%
0	156.62934219444\\
7.8800623635603	155.703517096719\\
15.7601247271206	154.790127328185\\
23.6401870906809	153.888868650908\\
31.5202494542412	152.999481505514\\
39.4003118178015	152.121679660035\\
47.2803741813618	151.255210233008\\
55.1604365449221	150.399782518165\\
63.0404989084824	149.55514263654\\
70.9205612720427	148.721065705566\\
78.800623635603	147.897284554299\\
86.6806859991633	147.083570084246\\
94.5607483627236	146.279684459736\\
102.440810726284	145.48543661374\\
110.320873089844	144.700596518702\\
118.200935453404	143.924952818154\\
126.080997816965	143.15830082643\\
133.961060180525	142.40045863973\\
141.841122544085	141.651339651975\\
149.721184907646	140.910551119588\\
};
\addlegendentry{Cordon 2}

\addplot [color=black, line width=1.0pt]
  table[row sep=crcr]{%
0	156.62934219444\\
7.8800623635603	155.724468672992\\
15.7601247271206	154.819813621534\\
23.6401870906809	153.91752288515\\
31.5202494542412	153.015687382815\\
39.4003118178015	152.116087569097\\
47.2803741813618	151.217182751437\\
55.1604365449221	150.320321720331\\
63.0404989084824	149.424415341109\\
70.9205612720427	148.530375819515\\
78.800623635603	147.637560279977\\
86.6806859991633	146.746414778347\\
94.5607483627236	145.856781738272\\
102.440810726284	144.968616327942\\
110.320873089844	144.082229187828\\
118.200935453404	143.197067060032\\
126.080997816965	142.314009635298\\
133.961060180525	141.431957771991\\
141.841122544085	140.552327482102\\
149.721184907646	139.673456530448\\
};
\addlegendentry{Trip-based}

\end{axis}
\end{tikzpicture}%
\vspace*{-0.3in}
\caption{Total trip volumes per minute under different target of traffic reduction in $\mathcal{C}$.} 
\label{figure13_2z}
\end{minipage}
\begin{minipage}[b]{0.005\linewidth}
\hfill
\end{minipage}
\begin{minipage}[b]{0.32\linewidth}
\centering
%
%
\definecolor{mycolor1}{rgb}{0.00000,0.44700,0.74100}%
\definecolor{mycolor2}{rgb}{0.85000,0.32500,0.09800}%
\begin{tikzpicture}

\begin{axis}[%
width=1.694in,
height=1.03in,
at={(1.358in,0.0in)},
scale only axis,
xmin=0,
xmax=150,
xlabel style={font=\color{white!15!black}},
xlabel={$\Delta N_\mathcal{C}$},
ymin=3300,
ymax=3700,
ytick={3300,3400, 3500,3600,3700},
yticklabels={{3.3K},{3.4K},{3.5K},{3.6K},{3.7K}},
ylabel style={font=\color{white!15!black}},
ylabel={Driver supply},
axis background/.style={fill=white},
legend style={at={(0.0,0.0)}, anchor=south west, legend cell align=left, align=left, font=\small, draw=white!15!black, nodes={scale=0.8, transform shape} }
]
\addplot [color=mycolor1, line width=1.0pt]
  table[row sep=crcr]{%
0	3687.35601392241\\
7.8800623635603	3677.82509585276\\
15.7601247271206	3668.32606595928\\
23.6401870906809	3658.85888636263\\
31.5202494542412	3649.42372544537\\
39.4003118178015	3640.01916907228\\
47.2803741813618	3630.64581619305\\
55.1604365449221	3621.30389623596\\
63.0404989084824	3611.99217529381\\
70.9205612720427	3602.71049592758\\
78.800623635603	3593.45947812281\\
86.6806859991633	3584.23684638025\\
94.5607483627236	3575.04746203021\\
102.440810726284	3565.88684793139\\
110.320873089844	3556.75395229572\\
118.200935453404	3547.64893671886\\
126.080997816965	3538.57236620781\\
133.961060180525	3529.52261952658\\
141.841122544085	3520.49888899471\\
149.721184907646	3511.50177769184\\
};
\addlegendentry{Cordon 1}

\addplot [color=mycolor2, line width=1.0pt]
  table[row sep=crcr]{%
0	3687.35601392241\\
7.8800623635603	3668.2870379405\\
15.7601247271206	3649.45755531301\\
23.6401870906809	3630.86274053188\\
31.5202494542412	3612.49951843515\\
39.4003118178015	3594.36338762091\\
47.2803741813618	3576.45061883091\\
55.1604365449221	3558.75672391153\\
63.0404989084824	3541.27712631002\\
70.9205612720427	3524.00926849498\\
78.800623635603	3506.94787760312\\
86.6806859991633	3490.08969455873\\
94.5607483627236	3473.43065932092\\
102.440810726284	3456.96808360093\\
110.320873089844	3440.69732080623\\
118.200935453404	3424.61550181572\\
126.080997816965	3408.71858480541\\
133.961060180525	3393.00327415383\\
141.841122544085	3377.46923318597\\
149.721184907646	3362.10774604511\\
};
\addlegendentry{Cordon 2}

\addplot [color=black, line width=1.0pt]
  table[row sep=crcr]{%
0	3687.35601392241\\
7.8800623635603	3673.23886794355\\
15.7601247271206	3659.1203480482\\
23.6401870906809	3644.98695942218\\
31.5202494542412	3630.85116422955\\
39.4003118178015	3616.70355020748\\
47.2803741813618	3602.55234236668\\
55.1604365449221	3588.3905661229\\
63.0404989084824	3574.22426720137\\
70.9205612720427	3560.04912401358\\
78.800623635603	3545.86881464449\\
86.6806859991633	3531.68145638236\\
94.5607483627236	3517.48857517559\\
102.440810726284	3503.29033473322\\
110.320873089844	3489.08530410785\\
118.200935453404	3474.87559587063\\
126.080997816965	3460.65931013017\\
133.961060180525	3446.43988186397\\
141.841122544085	3432.21365145796\\
149.721184907646	3417.98528808963\\
};
\addlegendentry{Trip-based}

\end{axis}
\end{tikzpicture}%
\vspace*{-0.3in}
\caption{Total number of drivers  under different target of traffic reduction in the congestion area $\mathcal{C}$.}
\label{figure14_2z}
\end{minipage}
\begin{minipage}[b]{0.005\linewidth}
\hfill
\end{minipage}
\begin{minipage}[b]{0.32\linewidth}
\centering
%
%
\definecolor{mycolor1}{rgb}{0.00000,0.44700,0.74100}%
\definecolor{mycolor2}{rgb}{0.85000,0.32500,0.09800}%
\begin{tikzpicture}

\begin{axis}[%
width=1.694in,
height=1.03in,
at={(1.358in,0.0in)},
scale only axis,
xmin=0,
xmax=150,
xlabel style={font=\color{white!15!black}},
xlabel={$\Delta N_\mathcal{C}$},
ymin=53000,
ymax=65000,
ylabel style={font=\color{white!15!black}},
ylabel={Platform profit},
axis background/.style={fill=white},
legend style={at={(0.0,0.0)}, anchor=south west, legend cell align=left, align=left, font=\small, draw=white!15!black, nodes={scale=0.8, transform shape} }
]
\addplot [color=mycolor1, line width=1.0pt]
  table[row sep=crcr]{%
0	64748.6208855417\\
7.8800623635603	64390.901724681\\
15.7601247271206	64035.4132742556\\
23.6401870906809	63682.1185457711\\
31.5202494542412	63330.9807432297\\
39.4003118178015	62981.9632809349\\
47.2803741813618	62635.0298015029\\
55.1604365449221	62290.144193159\\
63.0404989084824	61947.2706072105\\
70.9205612720427	61606.3734756499\\
78.800623635603	61267.4175279583\\
86.6806859991633	60930.3678080056\\
94.5607483627236	60595.3841805161\\
102.440810726284	60262.2375123837\\
110.320873089844	59930.8860211889\\
118.200935453404	59601.2963469871\\
126.080997816965	59273.4355337539\\
133.961060180525	58947.2710439914\\
141.841122544085	58622.7707728796\\
149.721184907646	58299.9030620403\\
};
\addlegendentry{Cordon 1}

\addplot [color=mycolor2, line width=1.0pt]
  table[row sep=crcr]{%
0	64748.6208855417\\
7.8800623635603	64012.3044270102\\
15.7601247271206	63292.217681577\\
23.6401870906809	62587.9141671576\\
31.5202494542412	61898.9606464216\\
39.4003118178015	61224.936938555\\
47.2803741813618	60565.4357157181\\
55.1604365449221	59920.0622854566\\
63.0404989084824	59288.4343601345\\
70.9205612720427	58670.1818145292\\
78.800623635603	58064.9464328413\\
86.6806859991633	57472.3816460539\\
94.5607483627236	56892.152261073\\
102.440810726284	56323.934182549\\
110.320873089844	55767.4141286295\\
118.200935453404	55222.2893416737\\
126.080997816965	54688.2672950665\\
133.961060180525	54165.0653971\\
141.841122544085	53652.410693002\\
149.721184907646	53150.0395660581\\
};
\addlegendentry{Cordon 2}

\addplot [color=black, line width=1.0pt]
  table[row sep=crcr]{%
0	64748.6208855417\\
7.8800623635603	64166.5884854911\\
15.7601247271206	63585.060740055\\
23.6401870906809	63008.9944427007\\
31.5202494542412	62433.9369656922\\
39.4003118178015	61863.8338812083\\
47.2803741813618	61295.2431090466\\
55.1604365449221	60731.0996190811\\
63.0404989084824	60168.9710753237\\
70.9205612720427	59610.7828345821\\
78.800623635603	59055.1111321844\\
86.6806859991633	58502.8730703445\\
94.5607483627236	57953.6519144365\\
102.440810726284	57407.3582384854\\
110.320873089844	56864.5804297296\\
118.200935453404	56324.2246266896\\
126.080997816965	55787.8820651014\\
133.961060180525	55253.4569050098\\
141.841122544085	54723.5405943516\\
149.721184907646	54195.0381337743\\
};
\addlegendentry{Trip-based}

\end{axis}
\end{tikzpicture}%
\vspace*{-0.3in}
\caption{Platform profit per hour under different target of traffic reduction in $\mathcal{C}$.}
\label{figure15_2z}
\end{minipage}
\begin{minipage}[b]{0.32\linewidth}
\centering
%
%
\definecolor{mycolor1}{rgb}{0.00000,0.44700,0.74100}%
\definecolor{mycolor2}{rgb}{0.85000,0.32500,0.09800}%
\begin{tikzpicture}

\begin{axis}[%
width=1.694in,
height=1.03in,
at={(1.358in,0.0in)},
scale only axis,
xmin=0,
xmax=5000,
xlabel style={font=\color{white!15!black}},
xlabel={Tax Revenue/Hour},
ymin=145,
ymax=158,
ylabel style={font=\color{white!15!black}},
ylabel={Passegner Demand},
axis background/.style={fill=white},
legend style={at={(0.0,0.0)}, anchor=south west, legend cell align=left, align=left, font=\small, draw=white!15!black, nodes={scale=0.8, transform shape} }
]
\addplot [color=mycolor1, line width=1.0pt]
  table[row sep=crcr]{%
0	156.62934219444\\
274.83804424084	156.103495064518\\
549.67608848168	155.579430758415\\
824.51413272252	155.05709871934\\
1099.35217696336	154.536511505375\\
1374.1902212042	154.017537372639\\
1649.02826544504	153.500336637395\\
1923.86630968588	152.984738183818\\
2198.70435392672	152.470753842659\\
2473.54239816756	151.958375576188\\
2748.3804424084	151.447519654921\\
3023.21848664924	150.938211982275\\
3298.05653089008	150.430410423322\\
3572.89457513092	149.924045557554\\
3847.73261937176	149.419128656813\\
4122.5706636126	148.915638938991\\
4397.40870785344	148.413476709972\\
4672.24675209428	147.912651622171\\
4947.08479633512	147.413142145472\\
5221.92284057596	146.914886736775\\
};
\addlegendentry{Cordon 1}

\addplot [color=mycolor2, line width=1.0pt]
  table[row sep=crcr]{%
0	148.405560204055\\
0	156.62934219444\\
274.83804424084	156.17666610914\\
549.67608848168	155.723990023839\\
824.51413272252	155.277119126323\\
1099.35217696336	154.830523213346\\
1374.1902212042	154.389322203812\\
1649.02826544504	153.948657708026\\
1923.86630968588	153.513010168002\\
2198.70435392672	153.078150179562\\
2473.54239816756	152.647929996436\\
2748.3804424084	152.218734564103\\
3023.21848664924	151.793827064713\\
3298.05653089008	151.370172543763\\
3572.89457513092	150.950451791774\\
3847.73261937176	150.53219604329\\
4122.5706636126	150.117545062548\\
4397.40870785344	149.704563954082\\
4672.24675209428	149.294878888256\\
4947.08479633512	148.887062464957\\
5221.92284057596	148.482231044476\\
};
\addlegendentry{Cordon 2}

\addplot [color=black, line width=1.0pt]
  table[row sep=crcr]{%
0	156.62934219444\\
274.83804424084	156.294556819418\\
549.67608848168	155.959771444395\\
824.51413272252	155.624986069373\\
1099.35217696336	155.290200694351\\
1374.1902212042	154.955415319329\\
1649.02826544504	154.621416321693\\
1923.86630968588	154.287586526259\\
2198.70435392672	153.953756730824\\
2473.54239816756	153.619926935389\\
2748.3804424084	153.286097139954\\
3023.21848664924	152.952910359969\\
3298.05653089008	152.620076152904\\
3572.89457513092	152.287241945839\\
3847.73261937176	151.954407738774\\
4122.5706636126	151.621573531709\\
4397.40870785344	151.289214354735\\
4672.24675209428	150.957393443277\\
4947.08479633512	150.625572531818\\
5221.92284057596	150.293751620359\\
};
\addlegendentry{Trip-based}

\end{axis}

\end{tikzpicture}%
\vspace*{-0.3in}
\caption{Total trip volumes per minute under different target of tax revenue.} 
\label{supplement1_2z}
\end{minipage}
\begin{minipage}[b]{0.005\linewidth}
\hfill
\end{minipage}
\begin{minipage}[b]{0.32\linewidth}
\centering
%
%
\definecolor{mycolor1}{rgb}{0.00000,0.44700,0.74100}%
\definecolor{mycolor2}{rgb}{0.85000,0.32500,0.09800}%
\begin{tikzpicture}

\begin{axis}[%
width=1.694in,
height=1.03in,
at={(1.358in,0.0in)},
scale only axis,
xmin=0,
xmax=5000,
xlabel style={font=\color{white!15!black}},
xlabel={Tax Revenue/Hour},
ymin=3500,
ymax=3700,
ytick={3500,3600, 3700},
yticklabels={{3.5K},{3.6K},{3.7K}},
ylabel style={font=\color{white!15!black}},
ylabel={Driver Supply},
axis background/.style={fill=white},
legend style={at={(0.0,0.0)}, anchor=south west, legend cell align=left, align=left, font=\small, draw=white!15!black, nodes={scale=0.8, transform shape} }
]
\addplot [color=mycolor1, line width=1.0pt]
  table[row sep=crcr]{%
0	3687.35601392241\\
274.83804424084	3676.87569557953\\
549.67608848168	3666.43432119102\\
824.51413272252	3656.03183874154\\
1099.35217696336	3645.66851969041\\
1374.1902212042	3635.34168914029\\
1649.02826544504	3625.05425178815\\
1923.86630968588	3614.80395796337\\
2198.70435392672	3604.59002691716\\
2473.54239816756	3594.41390382425\\
2748.3804424084	3584.27243888251\\
3023.21848664924	3574.16767534769\\
3298.05653089008	3564.09789745202\\
3572.89457513092	3554.06183625228\\
3847.73261937176	3544.0599160969\\
4122.5706636126	3534.09302626814\\
4397.40870785344	3524.15712366819\\
4672.24675209428	3514.25394927353\\
4947.08479633512	3504.38285980497\\
5221.92284057596	3494.54196156344\\
};
\addlegendentry{Cordon 1}

\addplot [color=mycolor2, line width=1.0pt]
  table[row sep=crcr]{%
0	3687.35601392241\\
274.83804424084	3678.03236367851\\
549.67608848168	3668.70871343461\\
824.51413272252	3659.49686593228\\
1099.35217696336	3650.29031438704\\
1374.1902212042	3641.18812364045\\
1649.02826544504	3632.09631139348\\
1923.86630968588	3623.10237111989\\
2198.70435392672	3614.123794591\\
2473.54239816756	3605.23616601748\\
2748.3804424084	3596.36862427585\\
3023.21848664924	3587.58559771394\\
3298.05653089008	3578.82726736331\\
3572.89457513092	3570.14691398079\\
3847.73261937176	3561.49560056799\\
4122.5706636126	3552.9158948676\\
4397.40870785344	3544.3693607173\\
4672.24675209428	3535.88889743941\\
4947.08479633512	3527.44589206347\\
5221.92284057596	3519.06274658789\\
};
\addlegendentry{Cordon 2}

\addplot [color=black, line width=1.0pt]
  table[row sep=crcr]{%
0	3687.35601392241\\
274.83804424084	3682.13294804773\\
549.67608848168	3676.90988217305\\
824.51413272252	3671.68681629837\\
1099.35217696336	3666.46375042369\\
1374.1902212042	3661.24068454901\\
1649.02826544504	3656.01267329556\\
1923.86630968588	3650.78359796144\\
2198.70435392672	3645.55452262733\\
2473.54239816756	3640.32544729321\\
2748.3804424084	3635.09637195909\\
3023.21848664924	3629.86389736362\\
3298.05653089008	3624.62955891324\\
3572.89457513092	3619.39522046285\\
3847.73261937176	3614.16088201247\\
4122.5706636126	3608.92654356209\\
4397.40870785344	3603.68974873512\\
4672.24675209428	3598.45017054174\\
4947.08479633512	3593.21059234835\\
5221.92284057596	3587.97101415496\\
};
\addlegendentry{Trip-based}

\end{axis}

\end{tikzpicture}%
\vspace*{-0.3in}
\caption{Total number of drivers  under different target of tax revenue.} 
\label{supplement2_2z}
\end{minipage}
\begin{minipage}[b]{0.005\linewidth}
\hfill
\end{minipage}
\begin{minipage}[b]{0.32\linewidth}
\centering
%
%
\definecolor{mycolor1}{rgb}{0.00000,0.44700,0.74100}%
\definecolor{mycolor2}{rgb}{0.85000,0.32500,0.09800}%
\begin{tikzpicture}

\begin{axis}[%
width=1.694in,
height=1.03in,
at={(1.358in,0.0in)},
scale only axis,
xmin=0,
xmax=5000,
xlabel style={font=\color{white!15!black}},
xlabel={Tax Revenue/Hour},
ymin=57000,
ymax=65000,
ylabel style={font=\color{white!15!black}},
ylabel={Platform profit},
axis background/.style={fill=white},
legend style={at={(0.0,0.0)}, anchor=south west, legend cell align=left, align=left, font=\small, draw=white!15!black, nodes={scale=0.8, transform shape} }
]
\addplot [color=mycolor1, line width=1.0pt]
  table[row sep=crcr]{%
0	64748.6208855417\\
274.83804424084	64355.2683610763\\
549.67608848168	63964.640128053\\
824.51413272252	63576.6853943822\\
1099.35217696336	63191.3536698617\\
1374.1902212042	62808.5947987218\\
1649.02826544504	62428.3589926274\\
1923.86630968588	62050.5968616058\\
2198.70435392672	61675.2594469892\\
2473.54239816756	61302.298252092\\
2748.3804424084	60931.6652731093\\
3023.21848664924	60563.3130293812\\
3298.05653089008	60197.1945929527\\
3572.89457513092	59833.2636175226\\
3847.73261937176	59471.4743665815\\
4122.5706636126	59111.7817408195\\
4397.40870785344	58754.1413046911\\
4672.24675209428	58398.509311949\\
4947.08479633512	58044.8427304451\\
5221.92284057596	57693.0992658287\\
};
\addlegendentry{Cordon 1}

\addplot [color=mycolor2, line width=1.0pt]
  table[row sep=crcr]{%
0	64748.6208855417\\
274.83804424084	64355.2683610763\\
549.67608848168	63964.640128053\\
824.51413272252	63576.6853943822\\
1099.35217696336	63191.3536698617\\
1374.1902212042	62808.5947987218\\
1649.02826544504	62428.3589926274\\
1923.86630968588	62050.5968616058\\
2198.70435392672	61675.2594469892\\
2473.54239816756	61302.298252092\\
2748.3804424084	60931.6652731093\\
3023.21848664924	60563.3130293812\\
3298.05653089008	60197.1945929527\\
3572.89457513092	59833.2636175226\\
3847.73261937176	59471.4743665815\\
4122.5706636126	59111.7817408195\\
4397.40870785344	58754.1413046911\\
4672.24675209428	58398.509311949\\
4947.08479633512	58044.8427304451\\
5221.92284057596	57693.0992658287\\
};
\addlegendentry{Cordon 2}

\addplot [color=black, line width=1.0pt]
  table[row sep=crcr]{%
0	64748.6208855417\\
274.83804424084	64533.2803647867\\
549.67608848168	64317.9398440317\\
824.51413272252	64102.5993232767\\
1099.35217696336	63887.2588025217\\
1374.1902212042	63671.9182817667\\
1649.02826544504	63458.3942529507\\
1923.86630968588	63245.2610726486\\
2198.70435392672	63032.1278923465\\
2473.54239816756	62818.9947120443\\
2748.3804424084	62605.8615317422\\
3023.21848664924	62394.1533043204\\
3298.05653089008	62183.2263952924\\
3572.89457513092	61972.2994862645\\
3847.73261937176	61761.3725772366\\
4122.5706636126	61550.4456682086\\
4397.40870785344	61340.5524252764\\
4672.24675209428	61131.8304488037\\
4947.08479633512	60923.108472331\\
5221.92284057596	60714.3864958583\\
};
\addlegendentry{Trip-based}

\end{axis}

\end{tikzpicture}%
\vspace*{-0.3in}
\caption{Platform profit per hour under different target of tax revenue.} 
\label{supplement3_2z}
\end{minipage}
\end{figure*}

A congestion charge can achieve two objectives: (a) reduce the traffic congestion in the urban core, and (b) raise revenue to subsidize public transit. To evaluate the effectiveness of different congestion charge schemes in achieving these two goals, we conduct the following two sets of experiments:
\begin{itemize}
\item First, we measure traffic congestion in the congestion area by the total number of ride-sourcing vehicles within $\mathcal{C}$, i.e., $N_{\mathcal{C}}$. To compare the performance of different congestion charges in congestion mitigation, we set a  reduction target of $\Delta N_{\mathcal{C}}$ for $N_{\mathcal{C}}$,  and compare the passenger surplus, driver surplus and platform profit of these congestion charge schemes when they achieve the same level of traffic reduction $\Delta N_{\mathcal{C}}$. Figure  \ref{figure13_2z}-\ref{figure15_2z} show the comparison of passenger demand, driver supply, and platform profit, respectively, for different values of $\Delta N_{\mathcal{C}}$. It is evident that to achieve the same traffic reduction in $\mathcal{C}$, the one-directional cordon price leads to consistently  higher passenger surplus, higher driver surplus, and higher platform profit than other congestion charges. {We comment that this is because the one-directional cordon price can more precisely target at reducing traffic in the congestion area: (a) compared to the bi-directional cordon price, the one-directional cordon price only penalizes vehicle flow that deteriorates the traffic congestion in $\mathcal{C}$ (i.e., from $\mathcal{R}$ to $\mathcal{C}$), and does not penalize vehicle flow that mitigates the traffic congestion in $\mathcal{C}$ (i.e., from $\mathcal{C}$ to $\mathcal{R}$); (b) compared to the trip-based congestion charge, the one-directional cordon price not only penalizes occupied flows that contribute to the traffic congestion in $\mathcal{C}$, but also penalizes idle vehicles that intend to enter the congestion area. For these reasons, the one-directional cordon price can very precisely target at reducing $N_\mathcal{C}$ and achieve  the same level of traffic reduction at a smaller cost to passengers, drivers and the platform.  }

\item Second, to compare the performances of distinct congestion charges in revenue-raising, we set a tax revenue target  and compare the passenger surplus, driver surplus and platform profit of these congestion charge schemes when they achieve the same tax revenue. Figure  \ref{supplement1_2z}-\ref{supplement3_2z} show the comparison of passenger demand, driver supply, and platform profit, respectively, for different values of the tax target. It is clear that the trip-based congestion charge is most effective in revenue-raising: to achieve the same tax revenue, the trip-based charge leads to higher passenger surplus, higher driver surplus, and higher platform profit than other congestion charges. {This is intuitive as the trip-based congestion charge is imposed on the largest number of trips and enjoys the largest base. Therefore, it is more effective in revenue-raising than the other two congestion charge schemes. }
\end{itemize}

To further investigate how the tax burden is distributed among passengers, drivers, and the ride-sourcing platform, we fix the congestion charge level at $\$3$/trip for all the three schemes, and compute the passenger surplus, driver surplus, platform profit, and tax revenue under the three schemes. The results are summarized in Table \ref{tax_inci}.\footnote{The unit of passenger surplus, driver surplus, platform profit and tax revenue in Table \ref{tax_inci} is \$/hour. Passenger surplus is calculated as $\sum_{i=1}^M \sum_{j=1}^M \int_{c_{ij}}^{\infty} \lambda_{ij}^0 F_p(x) dx$, while driver surplus is calculated as $\int_{0}^{q} N_0 F_d(x) dx$.} Under the same congestion charge level, the trip-based congestion charge has the largest tax revenue and the largest impact on passenger surplus, driver surplus and platform profit. It is also interesting to note that the platform assumes significantly higher tax burden than  passengers and drivers under all congestion charge schemes. 

\begin{table}[ht]
\begin{center}
\begin{tabular}{ |p{1.5cm}||p{1.3cm}|p{1.3cm}|p{1.35cm}||p{1.3cm}|p{1.3cm}| p{1.35cm}||p{1.3cm}|p{1.3cm}|p{1.35cm}|    }
\hline
&\multicolumn{3}{|c|}{One-directional Cordon}& \multicolumn{3}{|c|}{Bi-directional Cordon} & \multicolumn{3}{|c|}{Trip-based Charge}\\ 
\hline
 &  Before charge & After charge &\% change & Before charge & After charge &\% change& Before charge & After charge &\% change  \\
\hline
passenger surplus & 94,234 & 88,668 &-5.91\% & 94,234 & 83,981 &-10.88\% & 94,234 & 73,040 &-22.49\%   \\
\hline
driver surplus  & 26,661 & 24,891 &-6.64\% &26,661 & 23,706 &-11.08\% & 26,661 & 22,027 &-17.38\%   \\
\hline
platform profit  & 64,749 & 57,693 &-10.90\% & 64,749 & 53,150 &-17.91\% & 64,749 & 45,510 &-29.71\%   \\
\hline
tax revenue   & 0 & 5222 &N.A. & 0 & 8861 &N.A.& 0& 22443 &N.A.  \\
\hline
\end{tabular}
\end{center}
\vspace{-0.5cm}
\caption{Distribution of tax burdens under distinct congestion charge schemes.}
\label{tax_inci}
\end{table}

\subsubsection{Sensitivity analysis}
To understand how parameter values affect the aforementioned results, we conduct a sensitivity analysis by perturbing the key parameters in both directions by 30\%. The following model parameters are selected for sensitivity study, including: $\alpha, \epsilon, \sigma, \eta, N_0$ and $\lambda_0$.\footnote{Since $\lambda_0$ is a matrix, we perturb each entry of the matrices proportionately.} We perturb these parameter values one by one while fixing other parameters at the nominal value, and solve the  optimal spatial pricing problem  under different values of the one-directional cordon price.\footnote{We focus on the one-directional cordon price in the sensitivity study since the major insights are 
derived for this form of congestion charge.} The total number of idle vehicles in the southwest corner of the city (including zip code zones 94122, 94116, 94132, 94112, 94124), $N_\mathcal{S}$, and the passenger arrival rates for trips that start from the remote area and end in the remote area, $\lambda_{22}$, under distinct cordon prices are shown in Figure \ref{figure16_2z}-\ref{figure27_2z} of Appendix B.  As the congestion charge is increased from $\$0$ to $\$3$, both $N_\mathcal{S}$ and $\lambda_{22}$ increase for all values of the model parameters. It is clear that the insights derived in Section \ref{simulation_section} always hold for all parameter values under the sensitivity test. Lastly, we point out that when the model parameters are perturbed by $30\%$, the total number of ride-sourcing trips at the optimal solution is perturbed by almost $50\%$.  This indicates that the sensitivity study covers a large range of model parameters of practical interest, and the  insights derived from the two-zone model are robust with respect to the variation of model parameters.

\section{Conclusion}
This paper evaluates the impact of congestion charges on the ride-sourcing market based on a network equilibrium model. The formulated model captures the intimate interactions among passenger waiting time, driver waiting time, passenger demand, driver supply, idle driver repositioning, traffic congestion, and network flow balance. The platform determines the location-differentiated price and driver payment, which induces a market equilibrium that in turn affects the platform profit. The overall problem is cast as a large-scale non-convex program. An algorithm is developed to approximately compute its optimal solution, and a tight upper bound is derived to characterize its performance. Based on the proposed model, we compared three forms of congestion charge: (a) a one-directional cordon-based charge that penalizes vehicles for entering the congestion area, (b) a bi-directional cordon-based charge that penalizes vehicles for entering or exiting the congestion area, (c) a trip-based congestion charge on all ride-sourcing trips in the city. Through numerical study, we showed that the one-directional congestion charge not only reduces the ride-sourcing traffic in the congestion area, but also reduces the generalized travel cost outside the congestion area and benefits passengers in these underserved zones. {We further show that compared to other congestion charges,  the one-directional cordon charge is more effective in congestion mitigation, while the trip-based congestion charge is more effective in revenue-raising. The toll revenues collected from the congestion charges can be used to subsidize public transit,  upgrade road infrastructure, or incentivize ridesourcing-transit integration to promote a more efficient multi-modal transportation system. } Future study includes extending the framework to capture inter-zone matching, studying the network dynamics, and applying the analysis to other cities based on the ride-sourcing data.

\section*{Acknowledgments} This research was supported by Hong Kong Research Grants Council under project HKUST26200420, the National Science Foundation EAGER award 1839843, and California Department of Transportation. We are very grateful to Joe Castiglione and Drew Cooper of San Francisco
County Transportation Authority for guidance and data.

\bibliographystyle{unsrt}
\bibliography{resourceprocurement}

\providecommand{\url}[1]{#1}
\begin{thebibliography}{10}

\bibitem{NYC2019surcharge}
NYCTLC.
\newblock {New York State's Congestion Surcharge}.
\newblock New York City Taxi and Limousine Commission, 2019.
\newblock \url{https://www1.nyc.gov/site/tlc/about/congestion-surcharge.page}.

\bibitem{NYC_cordon}
New York~State Assembly.
\newblock {Traffic mobility act A09633B}.
\newblock
  \url{https://assembly.state.ny.us/leg/?default_fld=&bn=9633&term=2015&Summary=Y&Memo=Y/}.

\bibitem{chicago_surcharge}
K.~Pierog.
\newblock Chicago approves traffic congestion tax on ride-hailing services.
\newblock Reuters, 2019.
\newblock
  \url{https://www.reuters.com/article/us-chicago-ridehailing-tax/chicago-approves-traffic-congestion-tax-on-ride-hailing-services-idUSKBN1Y02BV}.

\bibitem{SFmayor}
{Treasurer \& Tax Collector}.
\newblock {Traffic Congestion Mitigation Tax}.
\newblock
  \url{https://sftreasurer.org/business/taxes-fees/traffic-congestion-mitigation-tax-tcm#:~:text=The%20City%20imposes%20a%20Traffic,or%20private%20transit%20services%20vehicle.}

\bibitem{Massa_cordon}
Tod Newcombe.
\newblock {Massachusetts Bets Big by Raising Ride-Sharing Surcharges}.
\newblock
  \url{https://www.governing.com/finance/Massachusetts-Bets-Big-by-Raising-Ride-Sharing-Surcharges.html/}.

\bibitem{castiglione2018tncs}
J.~Castiglione, D.~Cooper, B.~Sana, D.~Tischler, T.~Chang, G.~D. Erhardt,
  S.~Roy, M.~Chen, and A.~Mucci.
\newblock {TNCs \& Congestion. Draft Report.}
\newblock San Francisco County Transportation Authority, 2018.

\bibitem{schaller2017empty}
B.~Schaller.
\newblock {Empty Seats, Full Streets. Fixing Manhattan's Traffic Problem}.
\newblock Schaller Consulting, December 2017.

\bibitem{qian2020impact}
X.~Qian, T.~Lei, J.~Xue, Z.~Lei, and S.~V. Ukkusuri.
\newblock Impact of transportation network companies on urban congestion:
  Evidence from large-scale trajectory data.
\newblock {\em Sustainable Cities and Society}, 55:102053, 2020.

\bibitem{balding2019}
M.~Balding, T.~Whinery, E.~Leshner, and E.~Womeldorff.
\newblock Estimated {TNC} share of {VMT} in six {US} metropolitan regions.
\newblock {\em Fehr \& Peers}, 2019.

\bibitem{yang1998network}
H.~Yang and S.C. Wong.
\newblock A network model of urban taxi services.
\newblock {\em Transportation Research Part B: Methodological}, 32(4):235--246,
  1998.

\bibitem{wong2001modeling}
K.~Wong, S.~C. Wong, and H.~Yang.
\newblock Modeling urban taxi services in congested road networks with elastic
  demand.
\newblock {\em Transportation Research Part B: Methodological}, 35(9):819--842,
  2001.

\bibitem{yang2002demand}
H.~Yang, S.~C. Wong, and K.~I. Wong.
\newblock Demand--supply equilibrium of taxi services in a network under
  competition and regulation.
\newblock {\em Transportation Research Part B: Methodological}, 36(9):799--819,
  2002.

\bibitem{yang2010equilibria}
H.~Yang, C.~Leung, S.~C. Wong, and M.~Bell.
\newblock Equilibria of bilateral taxi--customer searching and meeting on
  networks.
\newblock {\em Transportation Research Part B: Methodological},
  44(8-9):1067--1083, 2010.

\bibitem{wong2008modeling}
K.~Wong, S.~C. Wong, H.~Yang, and J.~Wu.
\newblock Modeling urban taxi services with multiple user classes and vehicle
  modes.
\newblock {\em Transportation Research Part B: Methodological},
  42(10):985--1007, 2008.

\bibitem{xu2019equilibrium}
Z.~Xu, Z.~Chen, and Y.~Yin.
\newblock Equilibrium analysis of urban traffic networks with ride-sourcing
  services.
\newblock {\em Available at SSRN 3422294}, 2019.

\bibitem{ban2019general}
X.~Ban, M.~Dessouky, J.~Pang, and R.~Fan.
\newblock A general equilibrium model for transportation systems with e-hailing
  services and flow congestion.
\newblock {\em Transportation Research Part B: Methodological}, 129:273--304,
  2019.

\bibitem{ghili2020spatial}
S.~Ghili, V.~Kumar, et~al.
\newblock Spatial distribution of supply and the role of market thickness:
  Theory and evidence from ride sharing.
\newblock Technical report, Cowles Foundation for Research in Economics, Yale
  University, 2020.

\bibitem{bimpikis2019spatial}
K.~Bimpikis, O.~Candogan, and D.~Saban.
\newblock Spatial pricing in ride-sharing networks.
\newblock {\em Operations Research}, 67(3):744--769, 2019.

\bibitem{zha2018geometric}
L.~Zha, Y.~Yin, and Z.~Xu.
\newblock Geometric matching and spatial pricing in ride-sourcing markets.
\newblock {\em Transportation Research Part C: Emerging Technologies},
  92:58--75, 2018.

\bibitem{chen2019optimal}
X.~Chen, C.~Chen, and W.~Xie.
\newblock Optimal spatial pricing for an on-demand ride-sourcing service
  platform.
\newblock {\em Available at SSRN 3464228}, 2019.

\bibitem{guda2019your}
H.~Guda and U.~Subramanian.
\newblock Your uber is arriving: Managing on-demand workers through surge
  pricing, forecast communication, and worker incentives.
\newblock {\em Management Science}, 65(5):1995--2014, 2019.

\bibitem{ma2018spatio}
H.~Ma, F.~Fang, and D.~C. Parkes.
\newblock Spatio-temporal pricing for ridesharing platforms.
\newblock {\em arXiv preprint arXiv:1801.04015}, 2018.

\bibitem{afifah2020spatial}
F.~Afifah and Z.~Guo.
\newblock Spatial pricing of ride-sourcing services in a congested
  transportation network.
\newblock {\em arXiv preprint arXiv:2006.00164}, 2020.

\bibitem{he2015modeling}
F.~He and Z.~M. Shen.
\newblock Modeling taxi services with smartphone-based e-hailing applications.
\newblock {\em Transportation Research Part C: Emerging Technologies},
  58:93--106, 2015.

\bibitem{qinpiggyback}
J.~Qin, J.~Porter, K.~Poolla, and P.~Varaiya.
\newblock Piggyback on {TNCs} for electricity services: Spatial pricing and
  synergetic value.
\newblock In {\em American Control Conference}. IEEE, 2020.

\bibitem{pigou2017economics}
A.~Pigou.
\newblock {\em The economics of welfare}.
\newblock MacMillion, London, 1920.

\bibitem{walters1961theory}
A.~Walters.
\newblock The theory and measurement of private and social cost of highway
  congestion.
\newblock {\em Econometrica: Journal of the Econometric Society}, pages
  676--699, 1961.

\bibitem{vickrey1963pricing}
W.~S. Vickrey.
\newblock Pricing in urban and suburban transport.
\newblock {\em The American Economic Review}, 53(2):452--465, 1963.

\bibitem{beckmann1967optimal}
M.~Beckmann.
\newblock On optimal tolls for highways, tunnels, and bridges in vehicular
  traffic science.
\newblock In {\em Proceedings of 3rd Symposium on the Theory of Traffic Flow},
  1965.

\bibitem{yang2005mathematical}
H.~Yang and H.~Huang.
\newblock {\em Mathematical and economic theory of road pricing}.
\newblock 2005.

\bibitem{lindsey2000traffic}
C.~R. Lindsey and E.~T. Verhoef.
\newblock Traffic congestion and congestion pricing.
\newblock Technical report, Tinbergen Institute Discussion Paper, 2000.

\bibitem{may2000effects}
A.~D. May and D.~S. Milne.
\newblock Effects of alternative road pricing systems on network performance.
\newblock {\em Transportation Research Part A: Policy and Practice},
  34(6):407--436, 2000.

\bibitem{zhang2004optimal}
X.~Zhang and H.~Yang.
\newblock The optimal cordon-based network congestion pricing problem.
\newblock {\em Transportation Research Part B: Methodological}, 38(6):517--537,
  2004.

\bibitem{yang2010road}
H.~Yang, W.~Xu, B.~He, and Q.~Meng.
\newblock Road pricing for congestion control with unknown demand and cost
  functions.
\newblock {\em Transportation Research Part C: Emerging Technologies},
  18(2):157--175, 2010.

\bibitem{de2005congestion}
A.~Palma, M.~Kilani, and R.~Lindsey.
\newblock Congestion pricing on a road network: A study using the dynamic
  equilibrium simulator metropolis.
\newblock {\em Transportation Research Part A: Policy and Practice},
  39(7-9):588--611, 2005.

\bibitem{wu2011pareto}
D.~Wu, Y.~Yin, and S.~Lawphongpanich.
\newblock Pareto-improving congestion pricing on multimodal transportation
  networks.
\newblock {\em European Journal of Operational Research}, 210(3):660--669,
  2011.

\bibitem{wu2012design}
D.~Wu, Y.~Yin, S.~Lawphongpanich, and H.~Yang.
\newblock Design of more equitable congestion pricing and tradable credit
  schemes for multimodal transportation networks.
\newblock {\em Transportation Research Part B: Methodological},
  46(9):1273--1287, 2012.

\bibitem{simoni2019congestion}
M.~D. Simoni, K.~M. Kockelman, K.~M. Gurumurthy, and J.~Bischoff.
\newblock Congestion pricing in a world of self-driving vehicles: An analysis
  of different strategies in alternative future scenarios.
\newblock {\em Transportation Research Part C: Emerging Technologies},
  98:167--185, 2019.

\bibitem{mehr2019pricing}
N.~Mehr and R.~Horowitz.
\newblock Pricing traffic networks with mixed vehicle autonomy.
\newblock {\em arXiv preprint arXiv:1904.01226}, 2019.

\bibitem{salazar2019intermodal}
M.~Salazar, N.~Lanzetti, F.~Rossi, M.~Schiffer, and M.~Pavone.
\newblock Intermodal autonomous mobility-on-demand.
\newblock {\em IEEE Transactions on Intelligent Transportation Systems}, 2019.

\bibitem{li2019regulating}
S.~Li, H.~Tavafoghi, K.~Poolla, and P.~Varaiya.
\newblock Regulating {TNCs}: Should {Uber} and {Lyft} set their own rules?
\newblock {\em Transportation Research Part B: Methodological}, 2019.

\bibitem{li2020impact}
S.~Li, K.~Poolla, and P.~Varaiya.
\newblock Impact of congestion charge and minimum wage on tncs: A case study
  for {San Francisco}.
\newblock {\em arXiv preprint arXiv:2003.02550}, 2020.

\bibitem{mohring1987values}
H.~Mohring, J.~Schroeter, and P.~Wiboonchutikula.
\newblock The values of waiting time, travel time, and a seat on a bus.
\newblock {\em The RAND Journal of Economics}, pages 40--56, 1987.

\bibitem{nourinejad2020ride}
M.~Nourinejad and M.~Ramezani.
\newblock Ride-sourcing modeling and pricing in non-equilibrium two-sided
  markets.
\newblock {\em Transportation Research Part B: Methodological}, 132:340--357,
  2020.

\bibitem{banerjee2015pricing}
S.~Banerjee, C.~Riquelme, and R.~Johari.
\newblock Pricing in ride-share platforms: A queueing-theoretic approach.
\newblock {https://ssrn.com/abstract=2568258}, 2015.

\bibitem{harrison1993response}
P.~G. Harrison.
\newblock Response time distributions in queueing network models.
\newblock In {\em Performance Evaluation of Computer and Communication
  Systems}, pages 147--164. Springer, 1993.

\bibitem{mas1995microeconomic}
A.~Mas-Colell, M.~D. Whinston, J.~R. Green, et~al.
\newblock {\em Microeconomic theory}, volume~1.
\newblock Oxford university press New York, 1995.

\bibitem{bertsekas1997nonlinear}
D.~P. Bertsekas.
\newblock {\em Nonlinear programming}.
\newblock Athena Scientific Belmont, MA, 1998.

\bibitem{SFCTA2016}
SFCTA.
\newblock {TNC pickup and dropoff data in San Francisco}.
\newblock {http://tncsandcongestion.sfcta.org/}, 2016.

\bibitem{arnott1996taxi}
R.~Arnott.
\newblock Taxi travel should be subsidized.
\newblock {\em Journal of Urban Economics}, 40:316--333, 1996.

\bibitem{castiglione2016tncs}
J.~Castiglione, T.~Chang, D.~Cooper, J.~Hobson, W.~Logan, E.~Young,
  B.~Charlton, C.~Wilson, A.~Mislove, L.~Chen, and Jiang S.
\newblock {TNCs today: a profile of San Francisco transportation network
  company activity. Final Report}.
\newblock San Francisco County Transportation Authority, 2017.

\bibitem{parrott2018earning}
J.~A. Parrott and M.~Reich.
\newblock {An earnings standard for {New York City's} app-based drivers:
  Economic analysis and policy assessment}.
\newblock {The New School, Center for New York City Affairs }, 2018.

\bibitem{lyftprice}
Lyft.
\newblock {San Francisco Lyft rates}.
\newblock {https://estimatefares.com/rates/san-francisco}.

\bibitem{varian2014intermediate}
H.~R. Varian.
\newblock {\em Intermediate microeconomics with calculus: a modern approach}.
\newblock WW Norton \& Company, 2014.

\end{thebibliography}

\newpage
\section*{Appendix}

\subsection*{\bf{A: Sensitivity Analysis for Algorithm 1}}
\label{AppendixB}
\begin{figure*}[h]%
\begin{minipage}[b]{0.32\linewidth}
\centering
%
%
\definecolor{mycolor1}{rgb}{0.00000,0.44700,0.74100}%
\definecolor{mycolor2}{rgb}{0.85000,0.32500,0.09800}%
\begin{tikzpicture}

\begin{axis}[%
width=1.694in,
height=1.03in,
at={(1.358in,0.0in)},
scale only axis,
unbounded coords=jump,
xmin=0.7,
xmax=1.3,
xtick={0.7,0.8,0.9, 1.0, 1.1,1.2, 1.3},
xticklabels={{0.7},{0.8},{0.9},{1.0}, {1.1},{1.2}, {1.3}},
xlabel style={font=\color{white!15!black}},
xlabel={$\text{scale of }\alpha$},
ymin=0.0,
ymax=0.10,
ytick={0.02,0.04,0.06,0.08, 0.10},
yticklabels={{2\%},{4\%},{6\%},{8\%}, {10\%}},
ylabel style={font=\color{white!15!black}},
ylabel={Upper bound},
axis background/.style={fill=white},
]
\addplot [color=black, line width=1.0pt]
  table[row sep=crcr]{%
0.7	0.0152850336045008\\
0.703030303030303	0.0153493350284144\\
0.706060606060606	0.0154149543996855\\
0.709090909090909	0.015481900476177\\
0.712121212121212	0.0155501820969892\\
0.715151515151515	0.0156198081880064\\
0.718181818181818	0.0156907877590159\\
0.721212121212121	0.0157631298988263\\
0.724242424242424	0.0158368437923134\\
0.727272727272727	0.0159119386946641\\
0.73030303030303	0.0159884239635122\\
0.733333333333333	0.0160663090335856\\
0.736363636363636	0.0161456034287222\\
0.739393939393939	0.0162263167615341\\
0.742424242424242	0.0163084587333259\\
0.745454545454545	0.0163920391343347\\
0.748484848484848	0.0164770678452375\\
0.751515151515151	0.0165635548372961\\
0.754545454545454	0.0166515101737707\\
0.757575757575757	0.0167409440097877\\
0.76060606060606	0.016831866594308\\
0.763636363636364	0.01692428827015\\
0.766666666666667	0.0170182194754468\\
0.76969696969697	0.0171136707428746\\
0.772727272727273	0.0172106527042674\\
0.775757575757576	0.0173091760874974\\
0.778787878787879	0.0174092517202512\\
0.781818181818182	0.0175108906018291\\
0.784848484848485	0.0176141035440499\\
0.787878787878788	0.0177189018936307\\
0.790909090909091	0.0178252968121118\\
0.793939393939394	0.0179332996373807\\
0.796969696969697	0.0180429218130174\\
0.8	0.0181541748889395\\
0.803030303030303	0.0182670705239713\\
0.806060606060606	0.0183816204854439\\
0.809090909090909	0.0184978366519616\\
0.812121212121212	0.0186157310136205\\
0.815151515151515	0.0187353156744476\\
0.818181818181818	0.0188566028529003\\
0.821212121212121	0.0189796048839896\\
0.824242424242424	0.0191043342201538\\
0.827272727272727	0.0192308034333255\\
0.83030303030303	0.0193590252166629\\
0.833333333333333	0.0194890123844718\\
0.836363636363636	0.0196207778776782\\
0.839393939393939	0.0197543347583537\\
0.842424242424242	0.0198896962213393\\
0.845454545454545	0.0200268755875145\\
0.848484848484849	0.020165886309097\\
0.851515151515151	0.0203067419713718\\
0.854545454545454	0.0204494562939067\\
0.857575757575758	0.0205940431327324\\
0.86060606060606	0.0207405164826716\\
0.863636363636364	0.0208888904782041\\
0.866666666666667	0.0210391793965193\\
0.86969696969697	0.0211913976591757\\
0.872727272727273	0.0213455598339778\\
0.875757575757576	0.0215016806374278\\
0.878787878787879	0.0216597749367383\\
0.881818181818182	0.0218198577517459\\
0.884848484848485	0.0219819442578958\\
0.887878787878788	0.0221460497878695\\
0.890909090909091	0.0223121898341527\\
0.893939393939394	0.0224803800518125\\
0.896969696969697	0.0226506362596916\\
0.9	0.0228229744445976\\
0.903030303030303	0.0229974107621717\\
0.906060606060606	0.0231739615411484\\
0.909090909090909	0.0233526432844043\\
0.912121212121212	0.0235334726727682\\
0.915151515151515	0.0237164665668647\\
0.918181818181818	0.0239016420109554\\
0.921212121212121	0.0240890162348803\\
0.924242424242424	0.0242786066570443\\
0.927272727272727	0.0244704308882724\\
0.93030303030303	0.02466450673347\\
0.933333333333333	0.0248608521959494\\
0.936363636363636	0.0250594854799578\\
0.939393939393939	0.0252604249938086\\
0.942424242424242	0.0254636893532548\\
0.945454545454545	0.025669297386036\\
0.948484848484848	0.0258772681316805\\
0.951515151515152	0.0260876208503742\\
0.954545454545455	0.0263003750223883\\
0.957575757575758	0.0265155503529202\\
0.960606060606061	0.0267331667753813\\
0.963636363636364	0.0269532444579574\\
0.966666666666667	0.027175803803102\\
0.96969696969697	0.0274008654543602\\
0.972727272727273	0.0276284503002536\\
0.975757575757576	0.0278585794779554\\
0.978787878787879	0.0280912743781409\\
0.981818181818182	0.0283265566456588\\
0.984848484848485	0.0285644481928888\\
0.987878787878788	0.0288049711951002\\
0.990909090909091	0.0290481480976764\\
0.993939393939394	0.0292940016256504\\
0.996969696969697	0.0295425547828922\\
1	0.0297938308590001\\
1.0030303030303	0.0300478534360902\\
1.00606060606061	0.0303046463925846\\
1.00909090909091	0.0305642339087756\\
1.01212121212121	0.0308266404731151\\
1.01515151515152	0.0310918908882115\\
1.01818181818182	0.0313600102755757\\
1.02121212121212	0.0316310240832848\\
1.02424242424242	0.0319049580987658\\
1.02727272727273	0.0321818384187851\\
1.03030303030303	0.0324616915301722\\
1.03333333333333	0.032744544243071\\
1.03636363636364	0.0330304237343122\\
1.03939393939394	0.0333193575483135\\
1.04242424242424	0.0336113736044988\\
1.04545454545455	0.0339065002052029\\
1.04848484848485	0.0342047660443516\\
1.05151515151515	0.0345062002153665\\
1.05454545454545	0.0348108322212198\\
1.05757575757576	0.0351186919809189\\
1.06060606060606	0.0354298098439199\\
1.06363636363636	0.0357442165950515\\
1.06666666666667	0.0360619434686593\\
1.06969696969697	0.036383022157312\\
1.07272727272727	0.0367074848240903\\
1.07575757575758	0.0370353641137176\\
1.07878787878788	0.0373666931662401\\
1.08181818181818	0.0377015056292801\\
1.08484848484848	0.0380398356722376\\
1.08787878787879	0.0383817179971366\\
1.09090909090909	0.0387271878707558\\
1.09393939393939	0.0390762810860296\\
1.0969696969697	0.039429034076237\\
1.1	0.0397854838395914\\
1.1030303030303	0.0401456680028695\\
1.10606060606061	0.0405096248340121\\
1.10909090909091	0.0408773932624179\\
1.11212121212121	0.0412490129027888\\
1.11515151515152	0.0416245240784285\\
1.11818181818182	0.0420039678486306\\
1.12121212121212	0.0423873860361606\\
1.12424242424242	0.0427756292741803\\
1.12727272727273	0.0431663169599233\\
1.13030303030303	0.043561917458087\\
1.13333333333333	0.043961667934881\\
1.13636363636364	0.0443656145748688\\
1.13939393939394	0.0447738045173431\\
1.14242424242424	0.0451862859555367\\
1.14545454545455	0.0456031081888368\\
1.14848484848485	0.0460243216801317\\
1.15151515151515	0.0464499781471105\\
1.15454545454545	0.0468801305923919\\
1.15757575757576	0.0473148334948564\\
1.16060606060606	0.0477541428216926\\
1.16363636363636	0.0481981162060623\\
1.16666666666667	0.0486468130824039\\
1.16969696969697	0.0491002948638769\\
1.17272727272727	0.049558625155133\\
1.17575757575758	0.0500218700097167\\
1.17878787878788	0.0504900982541353\\
1.18181818181818	0.0509633818866443\\
1.18484848484848	0.0514417966165738\\
1.18787878787879	0.0519254225503855\\
1.19090909090909	0.0524143451558823\\
1.19393939393939	0.0529086566244014\\
1.1969696969697	0.0534084579049023\\
1.2	0.0539138619195157\\
1.2030303030303	0.0544249990931329\\
1.20606060606061	0.0549420279272364\\
1.20909090909091	0.0554651591315054\\
1.21212121212121	0.0559947318262328\\
1.21515151515152	0.0565423914685467\\
1.21818181818182	0.0569196393032256\\
1.22121212121212	0.0576167083953441\\
1.22424242424242	0.0578856777198837\\
1.22727272727273	0.0587227354910789\\
1.23030303030303	0.0592048206242681\\
1.23333333333333	0.059808571436913\\
1.23636363636364	0.0601666260070697\\
1.23939393939394	0.0608905176821388\\
1.24242424242424	0.0612635829650435\\
1.24545454545455	0.0623129132410086\\
1.24848484848485	0.0629203988039249\\
1.25151515151515	0.0635288193074188\\
1.25454545454545	0.0641664776957231\\
1.25757575757576	0.0647932871089464\\
1.26060606060606	0.0654208242285159\\
1.26363636363636	0.0660818369850919\\
1.26666666666667	0.0665504206074537\\
1.26969696969697	0.0673783479464742\\
1.27272727272727	0.0678995801150144\\
1.27575757575758	0.0687398535607444\\
1.27878787878788	0.0694696128645711\\
1.28181818181818	0.0699891952865896\\
1.28484848484848	0.0708742690674339\\
1.28787878787879	0.0714639124464553\\
1.29090909090909	0.0722769410829606\\
1.29393939393939	0.0729079701069327\\
1.2969696969697	0.0738812818757699\\
1.3	0.0744596092176833\\
};

\end{axis}
\end{tikzpicture}%
\vspace*{-0.3in}
\caption{Upper bound on performance loss of Algorithm \ref{algorithm1} by comparing (\ref{optimalpricing_trip}) and (\ref{optimalpricing_relax}) under distinct values of $\alpha$.} 
\label{sensitivity1_mz}
\end{minipage}
\begin{minipage}[b]{0.005\linewidth}
\hfill
\end{minipage}
\begin{minipage}[b]{0.32\linewidth}
\centering
%
%
\definecolor{mycolor1}{rgb}{0.00000,0.44700,0.74100}%
\definecolor{mycolor2}{rgb}{0.85000,0.32500,0.09800}%

\pgfplotsset{scaled y ticks=false}

\begin{tikzpicture}

\begin{axis}[%
width=1.694in,
height=1.03in,
at={(1.358in,0.0in)},
scale only axis,
unbounded coords=jump,
xmin=0.7,
xmax=1.3,
xtick={0.7,0.8,0.9, 1.0, 1.1,1.2, 1.3},
xticklabels={{0.7},{0.8},{0.9},{1.0}, {1.1},{1.2}, {1.3}},
xlabel style={font=\color{white!15!black}},
xlabel={$\text{scale of } \epsilon$},
ymin=0.01,
ymax=0.06,
ytick={0.01,0.02,0.03,0.04, 0.05, 0.06},
yticklabels={{1\%}, {2\%},{3\%},{4\%}, {5\%}, {6\%}},
ylabel style={font=\color{white!15!black}},
ylabel={Upper bound},
axis background/.style={fill=white},
]
\addplot [color=black, line width=1.0pt]
  table[row sep=crcr]{%
0.7	0.0152831028034302\\
0.703030303030303	0.0153873416152291\\
0.706060606060606	0.0154927481479073\\
0.709090909090909	0.015599314115913\\
0.712121212121212	0.015707031153242\\
0.715151515151515	0.0158158908162372\\
0.718181818181818	0.0159258845848152\\
0.721212121212121	0.0160370038664922\\
0.724242424242424	0.0161492399957243\\
0.727272727272727	0.0162625842382633\\
0.73030303030303	0.0163770277920888\\
0.733333333333333	0.016492561789743\\
0.736363636363636	0.0166091773003906\\
0.739393939393939	0.0167268653314355\\
0.742424242424242	0.0168456168304966\\
0.745454545454545	0.0169654226878518\\
0.748484848484848	0.0170862737376345\\
0.751515151515151	0.0172081607597653\\
0.754545454545454	0.017331074482389\\
0.757575757575757	0.0174550055825687\\
0.76060606060606	0.0175799446893196\\
0.763636363636364	0.0177058823841569\\
0.766666666666667	0.017832809203824\\
0.76969696969697	0.0179607156412638\\
0.772727272727273	0.0180895921475925\\
0.775757575757576	0.018219429133559\\
0.778787878787879	0.0183502169716267\\
0.781818181818182	0.018481945996963\\
0.784848484848485	0.0186146065211416\\
0.787878787878788	0.0187481887748449\\
0.790909090909091	0.0188826830269897\\
0.793939393939394	0.0190180794733608\\
0.796969696969697	0.019154368273122\\
0.8	0.0192915395861757\\
0.803030303030303	0.0194295835269867\\
0.806060606060606	0.0195684901898701\\
0.809090909090909	0.0197082496455744\\
0.812121212121212	0.0198488519445718\\
0.815151515151515	0.0199902871091406\\
0.818181818181818	0.0201325451535556\\
0.821212121212121	0.0202756160663784\\
0.824242424242424	0.0204194898210239\\
0.827272727272727	0.0205641563765531\\
0.83030303030303	0.0207096056773862\\
0.833333333333333	0.0208558276546183\\
0.836363636363636	0.0210028122289782\\
0.839393939393939	0.0211505493100787\\
0.842424242424242	0.0212990287991933\\
0.845454545454545	0.0214482405894142\\
0.848484848484849	0.0215981745674594\\
0.851515151515151	0.021748820614859\\
0.854545454545454	0.0219001686083991\\
0.857575757575758	0.0220522084237084\\
0.86060606060606	0.0222049299326228\\
0.863636363636364	0.0223583230067443\\
0.866666666666667	0.0225123775192912\\
0.86969696969697	0.0226670833447768\\
0.872727272727273	0.022822430360053\\
0.875757575757576	0.0229784084461814\\
0.878787878787879	0.0231350074897593\\
0.881818181818182	0.0232922173830843\\
0.884848484848485	0.0234500280253322\\
0.887878787878788	0.0236084293252845\\
0.890909090909091	0.0237674111999916\\
0.893939393939394	0.0239269635777565\\
0.896969696969697	0.0240870763980373\\
0.9	0.0242477396131304\\
0.903030303030303	0.0244089431892738\\
0.906060606060606	0.0245706771068086\\
0.909090909090909	0.0247329313623156\\
0.912121212121212	0.0248956959695232\\
0.915151515151515	0.0250589609595787\\
0.918181818181818	0.0252227163824942\\
0.921212121212121	0.0253869523081174\\
0.924242424242424	0.0255516588281812\\
0.927272727272727	0.0257168260546479\\
0.93030303030303	0.0258824441243378\\
0.933333333333333	0.0260485031973207\\
0.936363636363636	0.0262149934581698\\
0.939393939393939	0.0263819051182429\\
0.942424242424242	0.0265492284156201\\
0.945454545454545	0.0267169536160058\\
0.948484848484848	0.0268850710151584\\
0.951515151515152	0.0270535709382794\\
0.954545454545455	0.0272224437420072\\
0.957575757575758	0.0273916798154737\\
0.960606060606061	0.0275612695777241\\
0.963636363636364	0.0277312034871889\\
0.966666666666667	0.0279014720342307\\
0.96969696969697	0.0280720657458508\\
0.972727272727273	0.028242975186539\\
0.975757575757576	0.0284141909589111\\
0.978787878787879	0.0285857037050334\\
0.981818181818182	0.0287575041073776\\
0.984848484848485	0.0289295828893236\\
0.987878787878788	0.029101930817145\\
0.990909090909091	0.029274538700913\\
0.993939393939394	0.0294473973947191\\
0.996969696969697	0.029620497798415\\
1	0.0297938308590001\\
1.0030303030303	0.0299673875711126\\
1.00606060606061	0.0301411589778607\\
1.00909090909091	0.0303151362098033\\
1.01212121212121	0.0304893103027295\\
1.01515151515152	0.0306636725627443\\
1.01818181818182	0.0308382142049281\\
1.02121212121212	0.0310129265349081\\
1.02424242424242	0.0311878009122791\\
1.02727272727273	0.0313628287580541\\
1.03030303030303	0.0315380015477713\\
1.03333333333333	0.0317133108183881\\
1.03636363636364	0.0318887481671132\\
1.03939393939394	0.0320643052534567\\
1.04242424242424	0.032239973800242\\
1.04545454545455	0.0324157455949437\\
1.04848484848485	0.0325916124912951\\
1.05151515151515	0.0327675664106052\\
1.05454545454545	0.0329435993428593\\
1.05757575757576	0.0331197033474418\\
1.06060606060606	0.0332958705575338\\
1.06363636363636	0.0334720931780137\\
1.06666666666667	0.0336483634888296\\
1.06969696969697	0.0338246738484845\\
1.07272727272727	0.0340010166838427\\
1.07575757575758	0.0341773845160222\\
1.07878787878788	0.0343537699374579\\
1.08181818181818	0.0345301656256821\\
1.08484848484848	0.0347065643431773\\
1.08787878787879	0.0348829589386183\\
1.09090909090909	0.0350593423490611\\
1.09393939393939	0.0352357076021523\\
1.0969696969697	0.0354120478167948\\
1.1	0.0355883562064704\\
1.1030303030303	0.0357646260815981\\
1.10606060606061	0.0359408508495975\\
1.10909090909091	0.0361170240198513\\
1.11212121212121	0.0362931392037994\\
1.11515151515152	0.0364691901181366\\
1.11818181818182	0.0366451705877094\\
1.12121212121212	0.0368210745471648\\
1.12424242424242	0.0369968960443605\\
1.12727272727273	0.0371726292434867\\
1.13030303030303	0.0373482684234069\\
1.13333333333333	0.0375238079896932\\
1.13636363636364	0.0376992424691298\\
1.13939393939394	0.0378745665167978\\
1.14242424242424	0.0380497749182834\\
1.14545454545455	0.0382248625951155\\
1.14848484848485	0.0383998246048882\\
1.15151515151515	0.0385746561481677\\
1.15454545454545	0.0387493525715652\\
1.15757575757576	0.0389239093700458\\
1.16060606060606	0.0390983221956474\\
1.16363636363636	0.0392725868573862\\
1.16666666666667	0.0394466993283635\\
1.16969696969697	0.0396206557511704\\
1.17272727272727	0.0397944524417533\\
1.17575757575758	0.039968085897339\\
1.17878787878788	0.0401415527992197\\
1.18181818181818	0.0403148500227919\\
1.18484848484848	0.0404879746411863\\
1.18787878787879	0.040660923934064\\
1.19090909090909	0.0408336953946784\\
1.19393939393939	0.0410062867382884\\
1.1969696969697	0.0411786959102533\\
1.2	0.0413509210956199\\
1.2030303030303	0.041522960728081\\
1.20606060606061	0.04169481350985\\
1.20909090909091	0.0418664783851148\\
1.21212121212121	0.0420379546271813\\
1.21515151515152	0.042209241764395\\
1.21818181818182	0.0423803396629887\\
1.22121212121212	0.0425512485046653\\
1.22424242424242	0.0427219688045386\\
1.22727272727273	0.0428925014577744\\
1.23030303030303	0.0430628477268416\\
1.23333333333333	0.043233009278662\\
1.23636363636364	0.0434029882059281\\
1.23939393939394	0.0435727870728532\\
1.24242424242424	0.0437424088401534\\
1.24545454545455	0.0439118571030977\\
1.24848484848485	0.044081135920533\\
1.25151515151515	0.0442502499552955\\
1.25454545454545	0.0444192044975603\\
1.25757575757576	0.0445880055124025\\
1.26060606060606	0.0447566596950812\\
1.26363636363636	0.0449251745307219\\
1.26666666666667	0.0450935583642499\\
1.26969696969697	0.0452618204812549\\
1.27272727272727	0.0454299711962655\\
1.27575757575758	0.0455980219584782\\
1.27878787878788	0.0457659854736758\\
1.28181818181818	0.0459338758459979\\
1.28484848484848	0.0461017087450149\\
1.28787878787879	0.0462695016045374\\
1.29090909090909	0.0464372738616932\\
1.29393939393939	0.0466050472449597\\
1.2969696969697	0.0467728461282882\\
1.3	0.0469406979718433\\
};

\end{axis}
\end{tikzpicture}%
\vspace*{-0.3in}
\caption{Upper bound on performance loss of Algorithm \ref{algorithm1} by comparing (\ref{optimalpricing_trip}) and (\ref{optimalpricing_relax}) under distinct values of $\epsilon$.} 
\label{sensitivity2_mz}
\end{minipage}
\begin{minipage}[b]{0.005\linewidth}
\hfill
\end{minipage}
\begin{minipage}[b]{0.32\linewidth}
\centering
%
%
\definecolor{mycolor1}{rgb}{0.00000,0.44700,0.74100}%
\definecolor{mycolor2}{rgb}{0.85000,0.32500,0.09800}%

\pgfplotsset{scaled y ticks=false}

\begin{tikzpicture}

\begin{axis}[%
width=1.694in,
height=1.03in,
at={(1.358in,0.0in)},
scale only axis,
unbounded coords=jump,
xmin=0.7,
xmax=1.3,
xtick={0.7,0.8,0.9, 1.0, 1.1,1.2, 1.3},
xticklabels={{0.7},{0.8},{0.9},{1.0}, {1.1},{1.2}, {1.3}},
xlabel style={font=\color{white!15!black}},
xlabel={$\text{scale of } \sigma$},
ymin=0.025,
ymax=0.035,
ytick={0.025,0.027,0.029, 0.031, 0.033, 0.035},
yticklabels={{2.5\%},{2.7\%},{23.9\%}, {3.1\%}, {3.3\%}, {3.5\%}},
ylabel style={font=\color{white!15!black}},
ylabel={Upper bound},
axis background/.style={fill=white},
legend style={at={(0.0,0.64)}, anchor=south west, legend cell align=left, align=left, font=\small, draw=white!15!black, nodes={scale=0.9, transform shape} }
]
\addplot [color=black, line width=1.0pt]
  table[row sep=crcr]{%
0.7	0.0268541208400551\\
0.703030303030303	0.0268912191677692\\
0.706060606060606	0.0269281276935995\\
0.709090909090909	0.0269648479219555\\
0.712121212121212	0.0270013813392235\\
0.715151515151515	0.0270377294178943\\
0.718181818181818	0.0270738936141542\\
0.721212121212121	0.0271098753687874\\
0.724242424242424	0.027145676107727\\
0.727272727272727	0.0271812972420114\\
0.73030303030303	0.0272167401676724\\
0.733333333333333	0.0272520062660229\\
0.736363636363636	0.0272870969045484\\
0.739393939393939	0.027322013436264\\
0.742424242424242	0.0273567571997916\\
0.745454545454545	0.0273913295206852\\
0.748484848484848	0.0274257317096594\\
0.751515151515151	0.0274599650667539\\
0.754545454545454	0.0274940308719995\\
0.757575757575757	0.0275279304006467\\
0.76060606060606	0.0275616649104483\\
0.763636363636364	0.0275952356465888\\
0.766666666666667	0.0276290810333254\\
0.76969696969697	0.0276618907164204\\
0.772727272727273	0.0276949774787876\\
0.775757575757576	0.0277279053239162\\
0.778787878787879	0.0277606754358914\\
0.781818181818182	0.0277932889873342\\
0.784848484848485	0.0278257471370256\\
0.787878787878788	0.0278580510336739\\
0.790909090909091	0.0278902018142978\\
0.793939393939394	0.0279222006045215\\
0.796969696969697	0.0279540485185171\\
0.8	0.0279857466598745\\
0.803030303030303	0.0280172961207187\\
0.806060606060606	0.0280486979826939\\
0.809090909090909	0.028079953316987\\
0.812121212121212	0.0281110631831172\\
0.815151515151515	0.0281420286322863\\
0.818181818181818	0.0281728507035621\\
0.821212121212121	0.0282035304263451\\
0.824242424242424	0.0282340688208212\\
0.827272727272727	0.0282644668962402\\
0.83030303030303	0.0282947256570739\\
0.833333333333333	0.028324846079089\\
0.836363636363636	0.0283548291577212\\
0.839393939393939	0.02838467586848\\
0.842424242424242	0.0284143871452444\\
0.845454545454545	0.0284439639689719\\
0.848484848484849	0.0284734072732206\\
0.851515151515151	0.0285027179929895\\
0.854545454545454	0.0285318970541567\\
0.857575757575758	0.0285609453731644\\
0.86060606060606	0.0285898638584186\\
0.863636363636364	0.0286186534098008\\
0.866666666666667	0.0286473149179221\\
0.86969696969697	0.0286758492660794\\
0.872727272727273	0.0287042573286653\\
0.875757575757576	0.0287325399739318\\
0.878787878787879	0.0287606980541538\\
0.881818181818182	0.0287887324257016\\
0.884848484848485	0.0288166439289281\\
0.887878787878788	0.028844433398357\\
0.890909090909091	0.02887210166098\\
0.893939393939394	0.0288996495355831\\
0.896969696969697	0.0289270778345966\\
0.9	0.0289543873617808\\
0.903030303030303	0.0289815789139561\\
0.906060606060606	0.0290086532811851\\
0.909090909090909	0.0290356112455398\\
0.912121212121212	0.0290624535824362\\
0.915151515151515	0.0290891810602328\\
0.918181818181818	0.0291157944403445\\
0.921212121212121	0.0291422944767847\\
0.924242424242424	0.0291686819185214\\
0.927272727272727	0.0291949575036181\\
0.93030303030303	0.0292211219693198\\
0.933333333333333	0.0292471760424481\\
0.936363636363636	0.0292731204440363\\
0.939393939393939	0.0292989558890526\\
0.942424242424242	0.0293246830860442\\
0.945454545454545	0.0293503027368597\\
0.948484848484848	0.0293758155382101\\
0.951515151515152	0.0294012221818027\\
0.954545454545455	0.0294265233443253\\
0.957575757575758	0.0294517197110844\\
0.960606060606061	0.0294768119512211\\
0.963636363636364	0.0295018007310164\\
0.966666666666667	0.0295266867106617\\
0.96969696969697	0.0295514705438382\\
0.972727272727273	0.0295761528804255\\
0.975757575757576	0.029600734363198\\
0.978787878787879	0.0296252156293352\\
0.981818181818182	0.0296495973127978\\
0.984848484848485	0.0296738800363014\\
0.987878787878788	0.029698064424898\\
0.990909090909091	0.0297221510932685\\
0.993939393939394	0.0297461406523321\\
0.996969696969697	0.0297700337074388\\
1	0.0297938308590001\\
1.0030303030303	0.0298175327025489\\
1.00606060606061	0.0298411398278261\\
1.00909090909091	0.0298646528205957\\
1.01212121212121	0.0298880722610908\\
1.01515151515152	0.0299113987247451\\
1.01818181818182	0.0299346327822575\\
1.02121212121212	0.0299577749994448\\
1.02424242424242	0.0299808259372238\\
1.02727272727273	0.0300037861525017\\
1.03030303030303	0.0300266561968976\\
1.03333333333333	0.0300494366176653\\
1.03636363636364	0.0300721279576916\\
1.03939393939394	0.0300947307551989\\
1.04242424242424	0.0301172455437972\\
1.04545454545455	0.0301396728536728\\
1.04848484848485	0.0301620132086905\\
1.05151515151515	0.0301842671302959\\
1.05454545454545	0.0302064351358234\\
1.05757575757576	0.0302285177364988\\
1.06060606060606	0.0302505154411883\\
1.06363636363636	0.0302724287531867\\
1.06666666666667	0.0302942581733203\\
1.06969696969697	0.0303160041968552\\
1.07272727272727	0.0303376673164129\\
1.07575757575758	0.0303592480197398\\
1.07878787878788	0.0303807467909245\\
1.08181818181818	0.0304021641095886\\
1.08484848484848	0.0304235004529316\\
1.08787878787879	0.0304447562930672\\
1.09090909090909	0.0304659320989433\\
1.09393939393939	0.0304870283355813\\
1.0969696969697	0.0305080454644371\\
1.1	0.0305289839433695\\
1.1030303030303	0.0305498442263787\\
1.10606060606061	0.0305706267640131\\
1.10909090909091	0.0305913320035214\\
1.11212121212121	0.0306119603884133\\
1.11515151515152	0.0306325123583087\\
1.11818181818182	0.0306529883504793\\
1.12121212121212	0.0306733887977737\\
1.12424242424242	0.0306937141299976\\
1.12727272727273	0.0307139647739537\\
1.13030303030303	0.0307341411527714\\
1.13333333333333	0.0307542436864673\\
1.13636363636364	0.0307742727915443\\
1.13939393939394	0.0307942288818013\\
1.14242424242424	0.0308141123673667\\
1.14545454545455	0.0308339236555575\\
1.14848484848485	0.0308536631504627\\
1.15151515151515	0.030873331253725\\
1.15454545454545	0.0308929283623959\\
1.15757575757576	0.0309124548705803\\
1.16060606060606	0.0309319111774846\\
1.16363636363636	0.0309512976552631\\
1.16666666666667	0.0309706147060204\\
1.16969696969697	0.0309898627074744\\
1.17272727272727	0.031009042040191\\
1.17575757575758	0.0310281530809469\\
1.17878787878788	0.03104719620486\\
1.18181818181818	0.0310661717837133\\
1.18484848484848	0.0310850801860425\\
1.18787878787879	0.0311039217786277\\
1.19090909090909	0.0311226969246446\\
1.19393939393939	0.0311414059852024\\
1.1969696969697	0.0311600493183644\\
1.2	0.0311786272796833\\
1.2030303030303	0.0311971402224152\\
1.20606060606061	0.0312155884958719\\
1.20909090909091	0.0312339724488886\\
1.21212121212121	0.0312522924261109\\
1.21515151515152	0.0312705487704793\\
1.21818181818182	0.0312887418218101\\
1.22121212121212	0.0313068719184491\\
1.22424242424242	0.0313249393950174\\
1.22727272727273	0.0313429445848542\\
1.23030303030303	0.0313608878177686\\
1.23333333333333	0.0313787694219931\\
1.23636363636364	0.0313965897237242\\
1.23939393939394	0.0314143490457043\\
1.24242424242424	0.0314320477087353\\
1.24545454545455	0.0314496860325699\\
1.24848484848485	0.0314672643325342\\
1.25151515151515	0.031484782923249\\
1.25454545454545	0.0315022421165201\\
1.25757575757576	0.0315196422214959\\
1.26060606060606	0.0315369835476797\\
1.26363636363636	0.0315542663984508\\
1.26666666666667	0.0315714910771836\\
1.26969696969697	0.0315886578853048\\
1.27272727272727	0.0316057671276644\\
1.27575757575758	0.0316228190834318\\
1.27878787878788	0.0316398140647091\\
1.28181818181818	0.0316567523583341\\
1.28484848484848	0.0316736342550167\\
1.28787878787879	0.0316904600435476\\
1.29090909090909	0.0317072300103563\\
1.29393939393939	0.0317239444397855\\
1.2969696969697	0.0317406036146545\\
1.3	0.031757207815652\\
};

\end{axis}
\end{tikzpicture}%
\vspace*{-0.3in}
\caption{Upper bound on performance loss of Algorithm \ref{algorithm1} by comparing (\ref{optimalpricing_trip}) and (\ref{optimalpricing_relax}) under distinct values of $\sigma$.} 
\label{sensitivity3_mz}
\end{minipage}
\begin{minipage}[b]{0.32\linewidth}
\centering
%
%
\definecolor{mycolor1}{rgb}{0.00000,0.44700,0.74100}%
\definecolor{mycolor2}{rgb}{0.85000,0.32500,0.09800}%

\pgfplotsset{scaled y ticks=false}
\begin{tikzpicture}

\begin{axis}[%
width=1.694in,
height=1.03in,
at={(1.358in,0.0in)},
scale only axis,
unbounded coords=jump,
xmin=0.7,
xmax=1.3,
xtick={0.7,0.8,0.9, 1.0, 1.1,1.2, 1.3},
xticklabels={{0.7},{0.8},{0.9},{1.0}, {1.1},{1.2}, {1.3}},
xlabel style={font=\color{white!15!black}},
xlabel={$\text{scale of } \eta$},
ymin=0.02,
ymax=0.05,
ytick={0.02, 0.03,0.04,0.05},
yticklabels={{2\%},{3\%},{4\%},{5\%}},
ylabel style={font=\color{white!15!black}},
ylabel={Upper bound},
axis background/.style={fill=white},
legend style={at={(0.235,0.7)}, anchor=south west, legend cell align=left, align=left, font=\small, draw=white!15!black, nodes={scale=0.9, transform shape} }
]
\addplot [color=black, line width=1.0pt]
  table[row sep=crcr]{%
0.7	0.0478714703699095\\
0.703030303030303	0.0476410616006416\\
0.706060606060606	0.0474119494601114\\
0.709090909090909	0.0471832784740684\\
0.712121212121212	0.0469559020082479\\
0.715151515151515	0.0467295342060645\\
0.718181818181818	0.0465041740202515\\
0.721212121212121	0.0462798204066063\\
0.724242424242424	0.0460564723225154\\
0.727272727272727	0.0458341287360119\\
0.73030303030303	0.0456127885832066\\
0.733333333333333	0.045392450850848\\
0.736363636363636	0.0451731145134455\\
0.739393939393939	0.044954778474517\\
0.742424242424242	0.0447374417709252\\
0.745454545454545	0.044521103331302\\
0.748484848484848	0.0443057620742353\\
0.751515151515151	0.0440914170328195\\
0.754545454545454	0.0438780671377079\\
0.757575757575757	0.043665711345941\\
0.76060606060606	0.0434543486136337\\
0.763636363636364	0.0432439778936496\\
0.766666666666667	0.0430345981358495\\
0.76969696969697	0.0428262082875919\\
0.772727272727273	0.0426188072922273\\
0.775757575757576	0.0424123940894961\\
0.778787878787879	0.0422069676151146\\
0.781818181818182	0.0420025267996928\\
0.784848484848485	0.0417990705701568\\
0.787878787878788	0.0415965978475278\\
0.790909090909091	0.041395107547802\\
0.793939393939394	0.0411945985811623\\
0.796969696969697	0.0409950698517964\\
0.8	0.0407965202577068\\
0.803030303030303	0.0405989486904666\\
0.806060606060606	0.0404023540346721\\
0.809090909090909	0.0402067351679356\\
0.812121212121212	0.0400120909604971\\
0.815151515151515	0.0398184202750051\\
0.818181818181818	0.0396257219657464\\
0.821212121212121	0.0394339948799273\\
0.824242424242424	0.0392432378555874\\
0.827272727272727	0.0390534497225773\\
0.83030303030303	0.0388646293013838\\
0.833333333333333	0.0386767754041235\\
0.836363636363636	0.0384898868332232\\
0.839393939393939	0.0383039623818335\\
0.842424242424242	0.03811900083322\\
0.845454545454545	0.0379350009608139\\
0.848484848484849	0.0377519615281944\\
0.851515151515151	0.0375698812883818\\
0.854545454545454	0.0373887589840858\\
0.857575757575758	0.0372085933474884\\
0.86060606060606	0.0370293830994236\\
0.863636363636364	0.0368511269505514\\
0.866666666666667	0.0366738235999859\\
0.86969696969697	0.036497471735773\\
0.872727272727273	0.0363220700344001\\
0.875757575757576	0.0361476171598023\\
0.878787878787879	0.0359741117664604\\
0.881818181818182	0.0358015524954777\\
0.884848484848485	0.0356299379762795\\
0.887878787878788	0.035459266826438\\
0.890909090909091	0.0352895376508401\\
0.893939393939394	0.0351207490429972\\
0.896969696969697	0.0349528995837537\\
0.9	0.0347859878402843\\
0.903030303030303	0.0346200123693918\\
0.906060606060606	0.0344549717138131\\
0.909090909090909	0.0342908644040176\\
0.912121212121212	0.03412768895728\\
0.915151515151515	0.0339654438789844\\
0.918181818181818	0.0338041276607168\\
0.921212121212121	0.0336437387817333\\
0.924242424242424	0.0334842757079518\\
0.927272727272727	0.0333257368962176\\
0.93030303030303	0.0331681207755903\\
0.933333333333333	0.0330114257839056\\
0.936363636363636	0.0328556503313561\\
0.939393939393939	0.0327007928187422\\
0.942424242424242	0.0325468516335972\\
0.945454545454545	0.0323938251503316\\
0.948484848484848	0.0322417117302881\\
0.951515151515152	0.032090509721638\\
0.954545454545455	0.0319402174592944\\
0.957575757575758	0.0317908332652336\\
0.960606060606061	0.03164235544844\\
0.963636363636364	0.0314947823031143\\
0.966666666666667	0.0313481121136528\\
0.96969696969697	0.0312023431483975\\
0.972727272727273	0.0310574736641027\\
0.975757575757576	0.0309135019290459\\
0.978787878787879	0.0307704260993157\\
0.981818181818182	0.0306282444666644\\
0.984848484848485	0.03048695521114\\
0.987878787878788	0.0303465565246489\\
0.990909090909091	0.0302070465860677\\
0.993939393939394	0.0300684235618802\\
0.996969696969697	0.0299306856058155\\
1	0.0297938308590001\\
1.0030303030303	0.0296578574499801\\
1.00606060606061	0.029522763494966\\
1.00909090909091	0.029388547097659\\
1.01212121212121	0.0292552063493157\\
1.01515151515152	0.0291227393294442\\
1.01818181818182	0.0289911441049973\\
1.02121212121212	0.0288604187310419\\
1.02424242424242	0.0287305612503425\\
1.02727272727273	0.0286015696962375\\
1.03030303030303	0.0284734420820859\\
1.03333333333333	0.028346176421393\\
1.03636363636364	0.0282197707082092\\
1.03939393939394	0.0280942229270596\\
1.04242424242424	0.0279695310509731\\
1.04545454545455	0.0278456930416591\\
1.04848484848485	0.0277227068498764\\
1.05151515151515	0.0276005704150018\\
1.05454545454545	0.027479281665508\\
1.05757575757576	0.0273588385190398\\
1.06060606060606	0.0272392388824703\\
1.06363636363636	0.0271204806518201\\
1.06666666666667	0.0270025617308613\\
1.06969696969697	0.0268854799408309\\
1.07272727272727	0.0267692332005437\\
1.07575757575758	0.0266538193467917\\
1.07878787878788	0.0265392362242199\\
1.08181818181818	0.026425481667466\\
1.08484848484848	0.0263134337953667\\
1.08787878787879	0.026200449541453\\
1.09090909090909	0.0260891675927845\\
1.09393939393939	0.0259787054518044\\
1.0969696969697	0.0258690609053554\\
1.1	0.0257602317309678\\
1.1030303030303	0.0256522156973181\\
1.10606060606061	0.0255450105637856\\
1.10909090909091	0.0254386140812207\\
1.11212121212121	0.0253330239917579\\
1.11515151515152	0.0252282380290534\\
1.11818181818182	0.0251242539179852\\
1.12121212121212	0.0250210693755417\\
1.12424242424242	0.0249186821102957\\
1.12727272727273	0.0248170898229015\\
1.13030303030303	0.0247162902060475\\
1.13333333333333	0.0246162809447458\\
1.13636363636364	0.0245170597161418\\
1.13939393939394	0.0244186241902168\\
1.14242424242424	0.0243209720293467\\
1.14545454545455	0.0242241008885912\\
1.14848484848485	0.0241280084162153\\
1.15151515151515	0.0240326922532048\\
1.15454545454545	0.0239381500337696\\
1.15757575757576	0.0238443793856819\\
1.16060606060606	0.0237513779295045\\
1.16363636363636	0.0236591432801214\\
1.16666666666667	0.0235676730454519\\
1.16969696969697	0.0234769648290766\\
1.17272727272727	0.0233870162227005\\
1.17575757575758	0.0232978248204554\\
1.17878787878788	0.0232093882052463\\
1.18181818181818	0.0231217039555694\\
1.18484848484848	0.0230347696445575\\
1.18787878787879	0.0229485828395783\\
1.19090909090909	0.0228631411031862\\
1.19393939393939	0.0227784419925172\\
1.1969696969697	0.0226944830601823\\
1.2	0.022611261851835\\
1.2030303030303	0.0225287759117992\\
1.20606060606061	0.022447022776897\\
1.20909090909091	0.0223659999807819\\
1.21212121212121	0.0222857050525634\\
1.21515151515152	0.0222061355162909\\
1.21818181818182	0.0221272888929329\\
1.22121212121212	0.0220491626989317\\
1.22424242424242	0.021971754446831\\
1.22727272727273	0.0218950616451265\\
1.23030303030303	0.0218190817992905\\
1.23333333333333	0.0217438124110631\\
1.23636363636364	0.0216692509778953\\
1.23939393939394	0.0215953949953406\\
1.24242424242424	0.02152224195478\\
1.24545454545455	0.021449789344862\\
1.24848484848485	0.0213780346512437\\
1.25151515151515	0.0213069753567319\\
1.25454545454545	0.0212366089413539\\
1.25757575757576	0.0211669586697929\\
1.26060606060606	0.0210979446556078\\
1.26363636363636	0.0210296417328506\\
1.26666666666667	0.0209620215841241\\
1.26969696969697	0.0208950816785967\\
1.27272727272727	0.0208288194814523\\
1.27575757575758	0.0207632324571566\\
1.27878787878788	0.020698318068116\\
1.28181818181818	0.0206341225231025\\
1.28484848484848	0.0205704970362527\\
1.28787878787879	0.0205075853097679\\
1.29090909090909	0.0204453360513938\\
1.29393939393939	0.020383746715782\\
1.2969696969697	0.0203228147564454\\
1.3	0.0202625376258232\\
};

\end{axis}
\end{tikzpicture}%
\vspace*{-0.3in}
\caption{Upper bound on performance loss of Algorithm \ref{algorithm1} by comparing (\ref{optimalpricing_trip}) and (\ref{optimalpricing_relax}) under distinct values of $\eta$.} 
\label{sensitivity4_mz}
\end{minipage}
\begin{minipage}[b]{0.005\linewidth}
\hfill
\end{minipage}
\begin{minipage}[b]{0.32\linewidth}
\centering
%
%
\definecolor{mycolor1}{rgb}{0.00000,0.44700,0.74100}%
\definecolor{mycolor2}{rgb}{0.85000,0.32500,0.09800}%
\pgfplotsset{scaled y ticks=false}

\begin{tikzpicture}

\begin{axis}[%
width=1.694in,
height=1.03in,
at={(1.358in,0.0in)},
scale only axis,
unbounded coords=jump,
xmin=0.7,
xmax=1.3,
xtick={0.7,0.8,0.9, 1.0, 1.1,1.2, 1.3},
xticklabels={{0.7},{0.8},{0.9},{1.0}, {1.1},{1.2}, {1.3}},
xlabel style={font=\color{white!15!black}},
xlabel={$\text{scale of N}_\text{0}$},
ymin=0.025,
ymax=0.04,
ytick={0.025, 0.030, 0.035, 0.04},
yticklabels={{2.5\%},{3\%}, {3.5\%}, {4\%}},
ylabel style={font=\color{white!15!black}},
ylabel={Upper bound},
axis background/.style={fill=white},
]
\addplot [color=black, line width=1.0pt]
  table[row sep=crcr]{%
0.7	0.0349214067972171\\
0.703030303030303	0.0348262450129694\\
0.706060606060606	0.0347327054186976\\
0.709090909090909	0.0346407411051554\\
0.712121212121212	0.034550307887449\\
0.715151515151515	0.034461363972159\\
0.718181818181818	0.0343738696908818\\
0.721212121212121	0.0342877872852096\\
0.724242424242424	0.0342030807252942\\
0.727272727272727	0.0341197155591512\\
0.73030303030303	0.0340376587841039\\
0.733333333333333	0.0339568787350656\\
0.736363636363636	0.0338773449874991\\
0.739393939393939	0.0337990282731121\\
0.742424242424242	0.0337219004050291\\
0.745454545454545	0.0336459342077638\\
0.748484848484848	0.0335711034608819\\
0.751515151515151	0.0334973828420905\\
0.754545454545454	0.0334247478778772\\
0.757575757575757	0.0333531749001927\\
0.76060606060606	0.0332826410028407\\
0.763636363636364	0.0332131240202174\\
0.766666666666667	0.0331446024181151\\
0.76969696969697	0.0330770554089543\\
0.772727272727273	0.0330104627727127\\
0.775757575757576	0.0329448049060081\\
0.778787878787879	0.0328800627766558\\
0.781818181818182	0.0328162179053685\\
0.784848484848485	0.0327532523375234\\
0.787878787878788	0.032691148623703\\
0.790909090909091	0.0326298898001935\\
0.793939393939394	0.0325694593690781\\
0.796969696969697	0.0325098412807987\\
0.8	0.0324510199164994\\
0.803030303030303	0.0323929800731203\\
0.806060606060606	0.0323357069470316\\
0.809090909090909	0.0322791861201868\\
0.812121212121212	0.0322234035468812\\
0.815151515151515	0.0321683455392322\\
0.818181818181818	0.0321140012454139\\
0.821212121212121	0.0320603501927118\\
0.824242424242424	0.0320073871646018\\
0.827272727272727	0.0319550973015471\\
0.83030303030303	0.0319034685346553\\
0.833333333333333	0.0318524890892576\\
0.836363636363636	0.0318021474718186\\
0.839393939393939	0.0317524324637204\\
0.842424242424242	0.0317033331122287\\
0.845454545454545	0.0316548387207988\\
0.848484848484849	0.0316069388431206\\
0.851515151515151	0.031559623273304\\
0.854545454545454	0.031512882042795\\
0.857575757575758	0.0314667054036814\\
0.86060606060606	0.0314210838378199\\
0.863636363636364	0.0313760080342815\\
0.866666666666667	0.0313314688932301\\
0.86969696969697	0.0312874575154865\\
0.872727272727273	0.0312439651985602\\
0.875757575757576	0.0312009834308228\\
0.878787878787879	0.0311585038836381\\
0.881818181818182	0.0311165184116473\\
0.884848484848485	0.0310750190434211\\
0.887878787878788	0.0310339979755222\\
0.890909090909091	0.0309934475731136\\
0.893939393939394	0.0309533603612222\\
0.896969696969697	0.0309137290201674\\
0.9	0.0308745463857711\\
0.903030303030303	0.0308358054400803\\
0.906060606060606	0.0307974993107796\\
0.909090909090909	0.030759621266\\
0.912121212121212	0.0307221647109337\\
0.915151515151515	0.0306851231846252\\
0.918181818181818	0.0306484903561443\\
0.921212121212121	0.0306122600216835\\
0.924242424242424	0.0305764261008076\\
0.927272727272727	0.0305409826336963\\
0.93030303030303	0.0305059237782173\\
0.933333333333333	0.0304712438066984\\
0.936363636363636	0.030436937103222\\
0.939393939393939	0.0304029981610387\\
0.942424242424242	0.0303694215798736\\
0.945454545454545	0.0303362020633519\\
0.948484848484848	0.0303033344140052\\
0.951515151515152	0.0302708135390383\\
0.954545454545455	0.0302386344338634\\
0.957575757575758	0.030206792196299\\
0.960606060606061	0.0301752820113253\\
0.963636363636364	0.0301440991547694\\
0.966666666666667	0.0301132389903151\\
0.96969696969697	0.0300826969678829\\
0.972727272727273	0.0300524686196946\\
0.975757575757576	0.0300225495614459\\
0.978787878787879	0.0299929354870304\\
0.981818181818182	0.0299636221703284\\
0.984848484848485	0.0299346054600422\\
0.987878787878788	0.029905881279978\\
0.990909090909091	0.0298774456272176\\
0.993939393939394	0.0298492945694303\\
0.996969696969697	0.0298214242446866\\
1	0.0297938308590001\\
1.0030303030303	0.0297665106850105\\
1.00606060606061	0.0297394600606373\\
1.00909090909091	0.0297126753877705\\
1.01212121212121	0.0296861531304293\\
1.01515151515152	0.0296598898138342\\
1.01818181818182	0.0296338820230433\\
1.02121212121212	0.0296081264022163\\
1.02424242424242	0.0295826196513558\\
1.02727272727273	0.0295573585273872\\
1.03030303030303	0.0295323398427347\\
1.03333333333333	0.0295075604620016\\
1.03636363636364	0.0294830173035666\\
1.03939393939394	0.0294587073369864\\
1.04242424242424	0.0294346275818067\\
1.04545454545455	0.0294107751079245\\
1.04848484848485	0.029387147032504\\
1.05151515151515	0.0293637405208926\\
1.05454545454545	0.029340552784684\\
1.05757575757576	0.0293175810808997\\
1.06060606060606	0.0292948227103552\\
1.06363636363636	0.029272275019183\\
1.06666666666667	0.0292499353952799\\
1.06969696969697	0.0292278012679587\\
1.07272727272727	0.0292058701187577\\
1.07575757575758	0.02918413942915\\
1.07878787878788	0.0291626067795373\\
1.08181818181818	0.0291412697498715\\
1.08484848484848	0.0291201259675448\\
1.08787878787879	0.0290991730972579\\
1.09090909090909	0.0290784088404046\\
1.09393939393939	0.0290578309338329\\
1.0969696969697	0.0290374371507417\\
1.1	0.0290172252961367\\
1.1030303030303	0.0289971932116324\\
1.10606060606061	0.0289773387705153\\
1.10909090909091	0.0289576598787303\\
1.11212121212121	0.0289381544744703\\
1.11515151515152	0.0289188205268798\\
1.11818181818182	0.0288996560364053\\
1.12121212121212	0.0288806590320751\\
1.12424242424242	0.0288618275749871\\
1.12727272727273	0.0288431597530176\\
1.13030303030303	0.0288246536838826\\
1.13333333333333	0.028806307512552\\
1.13636363636364	0.0287881194118887\\
1.13939393939394	0.0287700875816865\\
1.14242424242424	0.0287522102477678\\
1.14545454545455	0.028734485662241\\
1.14848484848485	0.0287169121031557\\
1.15151515151515	0.0286994878730356\\
1.15454545454545	0.0286822112992251\\
1.15757575757576	0.0286650807333041\\
1.16060606060606	0.0286480945503795\\
1.16363636363636	0.0286312511490111\\
1.16666666666667	0.0286145489507109\\
1.16969696969697	0.028597986399623\\
1.17272727272727	0.0285815619619645\\
1.17575757575758	0.0285652741250938\\
1.17878787878788	0.028549121398671\\
1.18181818181818	0.0285331023125117\\
1.18484848484848	0.0285172154175797\\
1.18787878787879	0.0285014592848741\\
1.19090909090909	0.0284858325052605\\
1.19393939393939	0.0284703336893206\\
1.1969696969697	0.0284549614667356\\
1.2	0.0284397144860529\\
1.2030303030303	0.0284245914148207\\
1.20606060606061	0.028409590938298\\
1.20909090909091	0.0283947117597431\\
1.21212121212121	0.0283799526005413\\
1.21515151515152	0.0283653121992721\\
1.21818181818182	0.0283507893111113\\
1.22121212121212	0.0283363827125387\\
1.22424242424242	0.0283220911808577\\
1.22727272727273	0.0283079135324068\\
1.23030303030303	0.0282938485850416\\
1.23333333333333	0.0282798951748422\\
1.23636363636364	0.0282660521542841\\
1.23939393939394	0.0282523183907374\\
1.24242424242424	0.0282386927668173\\
1.24545454545455	0.0282251741792008\\
1.24848484848485	0.0282117615398085\\
1.25151515151515	0.0281984537743157\\
1.25454545454545	0.028185249823705\\
1.25757575757576	0.0281721486378391\\
1.26060606060606	0.0281591491874877\\
1.26363636363636	0.0281462508372751\\
1.26666666666667	0.0281334514245526\\
1.26969696969697	0.028120751112827\\
1.27272727272727	0.0281081485352949\\
1.27575757575758	0.0280956427236038\\
1.27878787878788	0.0280832327220658\\
1.28181818181818	0.0280709175870725\\
1.28484848484848	0.0280586963868287\\
1.28787878787879	0.0280465682011631\\
1.29090909090909	0.0280345321216639\\
1.29393939393939	0.0280225872512988\\
1.2969696969697	0.0280107327041746\\
1.3	0.0279989676054511\\
};

\end{axis}
\end{tikzpicture}%
\vspace*{-0.3in}
\caption{Upper bound on performance loss of Algorithm \ref{algorithm1} by comparing (\ref{optimalpricing_trip}) and (\ref{optimalpricing_relax}) under distinct values of $N_0$.} 
\label{sensitivity5_mz}
\end{minipage}
\begin{minipage}[b]{0.005\linewidth}
\hfill
\end{minipage}
\begin{minipage}[b]{0.32\linewidth}
\centering
%
%
\definecolor{mycolor1}{rgb}{0.00000,0.44700,0.74100}%
\definecolor{mycolor2}{rgb}{0.85000,0.32500,0.09800}%

\pgfplotsset{scaled y ticks=false}
\begin{tikzpicture}

\begin{axis}[%
width=1.694in,
height=1.03in,
at={(1.358in,0.0in)},
scale only axis,
unbounded coords=jump,
xmin=0.7,
xmax=1.3,
xtick={0.7,0.8,0.9, 1.0, 1.1,1.2, 1.3},
xticklabels={{0.7},{0.8},{0.9},{1.0}, {1.1},{1.2}, {1.3}},
xlabel style={font=\color{white!15!black}},
xlabel={$\text{scale of }\lambda{}_\text{0}$},
ymin=0.02,
ymax=0.06,
ytick={0.02,0.03,0.04, 0.05, 0.06},
yticklabels={{2\%},{3\%},{4\%}, {5\%}, {6\%}},
ylabel style={font=\color{white!15!black}},
ylabel={Upper bound},
axis background/.style={fill=white},
]
\addplot [color=black, line width=1.0pt]
  table[row sep=crcr]{%
0.7	0.0518985769873291\\
0.703030303030303	0.0515361241506478\\
0.706060606060606	0.0511779164948533\\
0.709090909090909	0.0508238898483923\\
0.712121212121212	0.0504739812254326\\
0.715151515151515	0.0501281289955538\\
0.718181818181818	0.049786272671615\\
0.721212121212121	0.0494483528864847\\
0.724242424242424	0.0491143115353787\\
0.727272727272727	0.0487840915568589\\
0.73030303030303	0.0484576370028952\\
0.733333333333333	0.0481348929933888\\
0.736363636363636	0.0478158056906519\\
0.739393939393939	0.0475003222762932\\
0.742424242424242	0.0471883909258455\\
0.745454545454545	0.0468799607872245\\
0.748484848484848	0.0465749819616606\\
0.751515151515151	0.0462734054602608\\
0.754545454545454	0.0459751832271326\\
0.757575757575757	0.0456802680742333\\
0.76060606060606	0.0453886136835959\\
0.763636363636364	0.0451001745833447\\
0.766666666666667	0.0448149061286521\\
0.76969696969697	0.0445327644833409\\
0.772727272727273	0.0442537066014003\\
0.775757575757576	0.0439776902104825\\
0.778787878787879	0.0437046737934827\\
0.781818181818182	0.0434346165733338\\
0.784848484848485	0.0431674784955413\\
0.787878787878788	0.042903220212981\\
0.790909090909091	0.0426418030710278\\
0.793939393939394	0.0423831890915392\\
0.796969696969697	0.0421273409591451\\
0.8	0.0418742220063096\\
0.803030303030303	0.0416237962003317\\
0.806060606060606	0.0413760281320969\\
0.809090909090909	0.0411308833573437\\
0.812121212121212	0.0408883265672747\\
0.815151515151515	0.0406483252389888\\
0.818181818181818	0.0404108459476479\\
0.821212121212121	0.0401758561933282\\
0.824242424242424	0.039943324024026\\
0.827272727272727	0.0397132180226298\\
0.83030303030303	0.0394855072961549\\
0.833333333333333	0.0392601614654747\\
0.836363636363636	0.0390371506541189\\
0.839393939393939	0.0388164454785126\\
0.842424242424242	0.0385980170374568\\
0.845454545454545	0.0383818369029995\\
0.848484848484849	0.038167877110155\\
0.851515151515151	0.0379561101479823\\
0.854545454545454	0.0377465089518613\\
0.857575757575758	0.0375390468905476\\
0.86060606060606	0.0373336977627271\\
0.863636363636364	0.0371304357841784\\
0.866666666666667	0.0369292355833209\\
0.86969696969697	0.0367300721888225\\
0.872727272727273	0.0365329210260762\\
0.875757575757576	0.0363377579064962\\
0.878787878787879	0.0361445590211953\\
0.881818181818182	0.0359533009328773\\
0.884848484848485	0.0357639605701938\\
0.887878787878788	0.0355765152189626\\
0.890909090909091	0.0353909425161493\\
0.893939393939394	0.0352072204432549\\
0.896969696969697	0.0350253273195942\\
0.9	0.0348452417952118\\
0.903030303030303	0.0346669428465706\\
0.906060606060606	0.034490409767576\\
0.909090909090909	0.0343156221662839\\
0.912121212121212	0.0341425599574177\\
0.915151515151515	0.033971203357\\
0.918181818181818	0.0338015328769028\\
0.921212121212121	0.0336335293196567\\
0.924242424242424	0.0334671737724491\\
0.927272727272727	0.0333024476028331\\
0.93030303030303	0.0331393324508554\\
0.933333333333333	0.0329778102293414\\
0.936363636363636	0.0328178631133318\\
0.939393939393939	0.0326594735386955\\
0.942424242424242	0.0325026241964802\\
0.945454545454545	0.0323472980280817\\
0.948484848484848	0.032193478220837\\
0.951515151515152	0.032041148215681\\
0.954545454545455	0.0318902916446248\\
0.957575757575758	0.0317408924417198\\
0.960606060606061	0.0315929347242499\\
0.963636363636364	0.031446402845806\\
0.966666666666667	0.0313012813812068\\
0.96969696969697	0.0311575551217002\\
0.972727272727273	0.0310152090726418\\
0.975757575757576	0.0308742284486536\\
0.978787878787879	0.0307345986702762\\
0.981818181818182	0.0305963053603953\\
0.984848484848485	0.0304593343408795\\
0.987878787878788	0.0303236716291152\\
0.990909090909091	0.0301893034342443\\
0.993939393939394	0.0300562161546331\\
0.996969696969697	0.0299243963742489\\
1	0.0297938308590001\\
1.0030303030303	0.0296645065546563\\
1.00606060606061	0.0295364105831875\\
1.00909090909091	0.0294095302397849\\
1.01212121212121	0.0292838529902558\\
1.01515151515152	0.0291593664678667\\
1.01818181818182	0.0290360584706998\\
1.02121212121212	0.0289139169591355\\
1.02424242424242	0.0287929300525084\\
1.02727272727273	0.0286730860269101\\
1.03030303030303	0.0285543733136722\\
1.03333333333333	0.0284367804938656\\
1.03636363636364	0.0283202962994338\\
1.03939393939394	0.0282049096078618\\
1.04242424242424	0.0280906094416275\\
1.04545454545455	0.0279773849650324\\
1.04848484848485	0.0278652254819535\\
1.05151515151515	0.0277541204337322\\
1.05454545454545	0.0276440593970283\\
1.05757575757576	0.027535032081663\\
1.06060606060606	0.0274270283278679\\
1.06363636363636	0.0273200381059333\\
1.06666666666667	0.027214051511051\\
1.06969696969697	0.0271090587653767\\
1.07272727272727	0.0270050502122186\\
1.07575757575758	0.0269020163165906\\
1.07878787878788	0.0267999476622524\\
1.08181818181818	0.0266988349506236\\
1.08484848484848	0.0265986689970697\\
1.08787878787879	0.0264994407319247\\
1.09090909090909	0.0264011411973387\\
1.09393939393939	0.02630376154515\\
1.0969696969697	0.0262072930328033\\
1.1	0.0261117270300264\\
1.1030303030303	0.0260170550070813\\
1.10606060606061	0.0259232685394722\\
1.10909090909091	0.0258303593041947\\
1.11212121212121	0.0257383190785959\\
1.11515151515152	0.0256471397691685\\
1.11818181818182	0.0255568132583223\\
1.12121212121212	0.0254673317065021\\
1.12424242424242	0.0253786872471032\\
1.12727272727273	0.0252908721368037\\
1.13030303030303	0.025203878724043\\
1.13333333333333	0.0251176994473949\\
1.13636363636364	0.0250323268345555\\
1.13939393939394	0.0249477535003635\\
1.14242424242424	0.0248639721463536\\
1.14545454545455	0.0247809755592933\\
1.14848484848485	0.0246987566090772\\
1.15151515151515	0.0246173082569908\\
1.15454545454545	0.0245366235123768\\
1.15757575757576	0.0244566955147083\\
1.16060606060606	0.0243775174487707\\
1.16363636363636	0.0242990825858361\\
1.16666666666667	0.0242213842744477\\
1.16969696969697	0.0241444159379418\\
1.17272727272727	0.0240681710744277\\
1.17575757575758	0.0239926432560378\\
1.17878787878788	0.0239178261269314\\
1.18181818181818	0.0238437134028968\\
1.18484848484848	0.0237702988697642\\
1.18787878787879	0.0236975763852095\\
1.19090909090909	0.0236255398678463\\
1.19393939393939	0.0235541833137259\\
1.1969696969697	0.0234850952055404\\
1.2	0.0234134863897122\\
1.2030303030303	0.0233441343318282\\
1.20606060606061	0.0232754388578411\\
1.20909090909091	0.0232073942830031\\
1.21212121212121	0.0231399949847438\\
1.21515151515152	0.0230732354004743\\
1.21818181818182	0.0230071100295327\\
1.22121212121212	0.0229416134300811\\
1.22424242424242	0.0228767402187045\\
1.22727272727273	0.0228124850703429\\
1.23030303030303	0.0227488427170484\\
1.23333333333333	0.0226858079465931\\
1.23636363636364	0.0226233756031672\\
1.23939393939394	0.0225615405855752\\
1.24242424242424	0.0225002978465054\\
1.24545454545455	0.022439642392356\\
1.24848484848485	0.0223795692818081\\
1.25151515151515	0.0223200736256483\\
1.25454545454545	0.0222611505859847\\
1.25757575757576	0.0222027953756459\\
1.26060606060606	0.0221450032573996\\
1.26363636363636	0.0220877695430065\\
1.26666666666667	0.0220310895931149\\
1.26969696969697	0.0219749588159763\\
1.27272727272727	0.0219193726675404\\
1.27575757575758	0.0218643266509854\\
1.27878787878788	0.0218098163144951\\
1.28181818181818	0.0217558372525469\\
1.28484848484848	0.0217023851044232\\
1.28787878787879	0.0216494555538182\\
1.29090909090909	0.0215970443279241\\
1.29393939393939	0.021545147197693\\
1.2969696969697	0.0214937599764162\\
1.3	0.0214428785196086\\
};

\end{axis}
\end{tikzpicture}%
\vspace*{-0.3in}
\caption{Upper bound on performance loss of Algorithm \ref{algorithm1} by comparing (\ref{optimalpricing_trip}) and (\ref{optimalpricing_relax}) under distinct values of $\lambda_0$.} 
\label{sensitivity6_mz}
\end{minipage}
\end{figure*}

\subsection*{\bf{B: Sensitivity Analysis for Impacts of Cordon Price in Two-Zone Model}}
\label{AppendixC}
\begin{figure*}[h]%
\begin{minipage}[b]{0.32\linewidth}
\centering
%
%
\definecolor{mycolor1}{rgb}{0.00000,0.44700,0.74100}%
\definecolor{mycolor2}{rgb}{0.85000,0.32500,0.09800}%
\begin{tikzpicture}

\begin{axis}[%
width=1.694in,
height=1.03in,
at={(1.358in,0.0in)},
scale only axis,
xmin=0,
xmax=3,
xlabel style={font=\color{white!15!black}},
xlabel={One-directional cordon},
ymin=210,
ymax=320,
ylabel style={font=\color{white!15!black}},
ylabel={Idle driver in $\mathcal{S}$},
axis background/.style={fill=white},
legend style={at={(0.0,0.58)}, anchor=south west, legend cell align=left, align=left, font=\small, draw=white!15!black, nodes={scale=0.8, transform shape} }
]
\addplot [color=mycolor1, line width=1.0pt]
  table[row sep=crcr]{%
0	283.01838727922\\
0.157894736842105	283.276094168437\\
0.315789473684211	283.535987015319\\
0.473684210526316	283.794151282487\\
0.631578947368421	284.05087997306\\
0.789473684210526	284.305357256429\\
0.947368421052632	284.5\\
1.10526315789474	284.808829502229\\
1.26315789473684	285.057342954035\\
1.42105263157895	285.303754239804\\
1.57894736842105	285.548092047997\\
1.73684210526316	285.789850768655\\
1.89473684210526	286.029448477067\\
2.05263157894737	286.266389504615\\
2.21052631578947	286.500840850734\\
2.36842105263158	286.732931301646\\
2.52631578947368	286.962209091863\\
2.68421052631579	287.188941119612\\
2.84210526315789	287.412959650567\\
3	287.634304616335\\
};
\addlegendentry{$\text{1.3}\alpha$}

\addplot [color=mycolor2, line width=1.0pt]
  table[row sep=crcr]{%
0	254.860670713077\\
0.157894736842105	255.287332075303\\
0.315789473684211	255.713604872895\\
0.473684210526316	256.139943691003\\
0.631578947368421	256.56606071728\\
0.789473684210526	256.991298911363\\
0.947368421052632	257.41646323123\\
1.10526315789474	257.840917301465\\
1.26315789473684	258.264365083502\\
1.42105263157895	258.687151166618\\
1.57894736842105	259.108578455799\\
1.73684210526316	259.529028648608\\
1.89473684210526	259.948024199677\\
2.05263157894737	260.36542697812\\
2.21052631578947	260.781150900994\\
2.36842105263158	261.195573957945\\
2.52631578947368	261.607622292257\\
2.68421052631579	262.017912774861\\
2.84210526315789	262.426152837562\\
3	262.831828421841\\
};
\addlegendentry{$\alpha$}

\addplot [color=black, line width=1.0pt]
  table[row sep=crcr]{%
0	225.75\\
0.157894736842105	226.040310419175\\
0.315789473684211	226.388660063037\\
0.473684210526316	226.73751841374\\
0.631578947368421	227.087556224427\\
0.789473684210526	227.437837766854\\
0.947368421052632	227.788756977625\\
1.10526315789474	228.140357726323\\
1.26315789473684	228.492052889685\\
1.42105263157895	228.844283315899\\
1.57894736842105	229.197261056288\\
1.73684210526316	229.550100552301\\
1.89473684210526	229.903403176549\\
2.05263157894737	230.256609171317\\
2.21052631578947	230.610304060613\\
2.36842105263158	230.963324135865\\
2.52631578947368	231.31707577553\\
2.68421052631579	231.670188309779\\
2.84210526315789	232.023577226939\\
3	232.37634451009\\
};
\addlegendentry{$\text{0.7}\alpha$}

\end{axis}
\end{tikzpicture}%
\vspace*{-0.3in}
\caption{Number of idle drivers in southwest corner\protect\footnotemark under the one-directional cordon charge for distinct $\alpha$} 
\label{figure16_2z}
\end{minipage}
\begin{minipage}[b]{0.005\linewidth}
\hfill
\end{minipage}
\begin{minipage}[b]{0.32\linewidth}
\centering
%
%
\definecolor{mycolor1}{rgb}{0.00000,0.44700,0.74100}%
\definecolor{mycolor2}{rgb}{0.85000,0.32500,0.09800}%
\begin{tikzpicture}

\begin{axis}[%
width=1.694in,
height=1.03in,
at={(1.358in,0.0in)},
scale only axis,
xmin=0,
xmax=3,
xlabel style={font=\color{white!15!black}},
xlabel={One-directional cordon},
ymin=210,
ymax=320,
ylabel style={font=\color{white!15!black}},
ylabel={Idle driver in $\mathcal{S}$},
axis background/.style={fill=white},
legend style={at={(0.0,0.58)}, anchor=south west, legend cell align=left, align=left, font=\small, draw=white!15!black, nodes={scale=0.8, transform shape} }
]
\addplot [color=mycolor1, line width=1.0pt]
  table[row sep=crcr]{%
0	274.918541714329\\
0.157894736842105	275.338649968459\\
0.315789473684211	275.753005751167\\
0.473684210526316	276.165216159683\\
0.631578947368421	276.575457718729\\
0.789473684210526	276.983303528768\\
0.947368421052632	277.389827454803\\
1.10526315789474	277.791792122922\\
1.26315789473684	278.192528760771\\
1.42105263157895	278.589676330516\\
1.57894736842105	278.984619892861\\
1.73684210526316	279.376707393038\\
1.89473684210526	279.765268753863\\
2.05263157894737	280.151183260436\\
2.21052631578947	280.533936559014\\
2.36842105263158	280.913264697834\\
2.52631578947368	281.289293204459\\
2.68421052631579	281.661775364278\\
2.84210526315789	282.031286501301\\
3	282.396090928874\\
};
\addlegendentry{$1.3\epsilon$}

\addplot [color=mycolor2, line width=1.0pt]
  table[row sep=crcr]{%
0	254.860670713077\\
0.157894736842105	255.287332075303\\
0.315789473684211	255.713604872895\\
0.473684210526316	256.139943691003\\
0.631578947368421	256.56606071728\\
0.789473684210526	256.991298911363\\
0.947368421052632	257.41646323123\\
1.10526315789474	257.840917301465\\
1.26315789473684	258.264365083502\\
1.42105263157895	258.687151166618\\
1.57894736842105	259.108578455799\\
1.73684210526316	259.529028648608\\
1.89473684210526	259.948024199677\\
2.05263157894737	260.36542697812\\
2.21052631578947	260.781150900994\\
2.36842105263158	261.195573957945\\
2.52631578947368	261.607622292257\\
2.68421052631579	262.017912774861\\
2.84210526315789	262.426152837562\\
3	262.831828421841\\
};
\addlegendentry{$\epsilon$}

\addplot [color=black, line width=1.0pt]
  table[row sep=crcr]{%
0	222.488392242521\\
0.157894736842105	222.909472730903\\
0.315789473684211	223.33430023458\\
0.473684210526316	223.762097340311\\
0.631578947368421	224.192977607138\\
0.789473684210526	224.626329047671\\
0.947368421052632	225.062625906635\\
1.10526315789474	225.501692861874\\
1.26315789473684	225.943232046771\\
1.42105263157895	226.387303793899\\
1.57894736842105	226.834439035769\\
1.73684210526316	227.281907114513\\
1.89473684210526	227.732365644026\\
2.05263157894737	228.184744500994\\
2.21052631578947	228.638959226903\\
2.36842105263158	229.094731222086\\
2.52631578947368	229.551890608187\\
2.68421052631579	230.010411239416\\
2.84210526315789	230.471038155217\\
3	230.931519974525\\
};
\addlegendentry{$0.7\epsilon$}

\end{axis}
\end{tikzpicture}%
\vspace*{-0.3in}
\caption{Number of idle drivers in southwest corner  under the one-directional cordon charge for distinct $\epsilon$.}
\label{figure17_2z}
\end{minipage}
\begin{minipage}[b]{0.005\linewidth}
\hfill
\end{minipage}
\begin{minipage}[b]{0.32\linewidth}
\centering
%
%
\definecolor{mycolor1}{rgb}{0.00000,0.44700,0.74100}%
\definecolor{mycolor2}{rgb}{0.85000,0.32500,0.09800}%
\begin{tikzpicture}

\begin{axis}[%
width=1.694in,
height=1.03in,
at={(1.358in,0.0in)},
scale only axis,
xmin=0,
xmax=3,
xlabel style={font=\color{white!15!black}},
xlabel={One-directional cordon},
ymin=220,
ymax=320,
ylabel style={font=\color{white!15!black}},
ylabel={Idle driver in $\mathcal{S}$},
axis background/.style={fill=white},
legend style={at={(0.0,0.58)}, anchor=south west, legend cell align=left, align=left, font=\small, draw=white!15!black, nodes={scale=0.8, transform shape} }
]
\addplot [color=mycolor1, line width=1.0pt]
  table[row sep=crcr]{%
0	266.7518664862\\
0.157894736842105	267.108358570445\\
0.315789473684211	267.464735409935\\
0.473684210526316	267.820892471435\\
0.631578947368421	268.176729892005\\
0.789473684210526	268.532680221634\\
0.947368421052632	268.887575968076\\
1.10526315789474	269.241726468552\\
1.26315789473684	269.595383301011\\
1.42105263157895	269.946035741317\\
1.57894736842105	270.298149194461\\
1.73684210526316	270.649430606957\\
1.89473684210526	270.998403673759\\
2.05263157894737	271.345919459946\\
2.21052631578947	271.691514638233\\
2.36842105263158	272.035386091202\\
2.52631578947368	272.377534410012\\
2.68421052631579	272.717713861782\\
2.84210526315789	273.055662988906\\
3	273.391280403449\\
};
\addlegendentry{$1.3\sigma$}

\addplot [color=mycolor2, line width=1.0pt]
  table[row sep=crcr]{%
0	254.860670713077\\
0.157894736842105	255.287332075303\\
0.315789473684211	255.713604872895\\
0.473684210526316	256.139943691003\\
0.631578947368421	256.56606071728\\
0.789473684210526	256.991298911363\\
0.947368421052632	257.41646323123\\
1.10526315789474	257.840917301465\\
1.26315789473684	258.264365083502\\
1.42105263157895	258.687151166618\\
1.57894736842105	259.108578455799\\
1.73684210526316	259.529028648608\\
1.89473684210526	259.948024199677\\
2.05263157894737	260.36542697812\\
2.21052631578947	260.781150900994\\
2.36842105263158	261.195573957945\\
2.52631578947368	261.607622292257\\
2.68421052631579	262.017912774861\\
2.84210526315789	262.426152837562\\
3	262.831828421841\\
};
\addlegendentry{$\sigma$}

\addplot [color=black, line width=1.0pt]
  table[row sep=crcr]{%
0	237.536275288134\\
0.157894736842105	238.049244007986\\
0.315789473684211	238.583135429467\\
0.473684210526316	239.106484134797\\
0.631578947368421	239.62977800582\\
0.789473684210526	240.152325070894\\
0.947368421052632	240.674727904171\\
1.10526315789474	241.196455127636\\
1.26315789473684	241.716420965791\\
1.42105263157895	242.237997681046\\
1.57894736842105	242.756847929403\\
1.73684210526316	243.274724362358\\
1.89473684210526	243.791500088165\\
2.05263157894737	244.306739601509\\
2.21052631578947	244.820665907895\\
2.36842105263158	245.334224945693\\
2.52631578947368	245.843153809003\\
2.68421052631579	246.351774778504\\
2.84210526315789	246.85821675475\\
3	247.362458693449\\
};
\addlegendentry{$0.7\sigma$}

\end{axis}
\end{tikzpicture}%
\vspace*{-0.3in}
\caption{Number of idle drivers in southwest corner  under the one-directional cordon charge for distinct $\sigma$.}
\label{figure18_2z}
\end{minipage}
\end{figure*}
\footnotetext{This includes zip code zones 94122, 94116, 94132, 94112, 94124.}

\begin{figure*}[h]
\begin{minipage}[b]{0.32\linewidth}
\centering
%
%
\definecolor{mycolor1}{rgb}{0.00000,0.44700,0.74100}%
\definecolor{mycolor2}{rgb}{0.85000,0.32500,0.09800}%
\begin{tikzpicture}

\begin{axis}[%
width=1.694in,
height=1.03in,
at={(1.358in,0.0in)},
scale only axis,
xmin=0,
xmax=3,
xlabel style={font=\color{white!15!black}},
xlabel={One-directional cordon},
ymin=220,
ymax=300,
ylabel style={font=\color{white!15!black}},
ylabel={Idle driver in $\mathcal{S}$},
axis background/.style={fill=white},
legend style={at={(0.0,0.58)}, anchor=south west, legend cell align=left, align=left, font=\small, draw=white!15!black, nodes={scale=0.8, transform shape} }
]
\addplot [color=mycolor1, line width=1.0pt]
  table[row sep=crcr]{%
0	239.058296584753\\
0.157894736842105	239.631665352799\\
0.315789473684211	240.206688745366\\
0.473684210526316	240.783252081988\\
0.631578947368421	241.361305577212\\
0.789473684210526	241.940654182659\\
0.947368421052632	242.520117185913\\
1.10526315789474	243.100561354035\\
1.26315789473684	243.681811060772\\
1.42105263157895	244.262639989764\\
1.57894736842105	244.843666464458\\
1.73684210526316	245.424187564449\\
1.89473684210526	246.004149447437\\
2.05263157894737	246.58095232401\\
2.21052631578947	247.161514451067\\
2.36842105263158	247.738661682281\\
2.52631578947368	248.31352251217\\
2.68421052631579	248.887076241676\\
2.84210526315789	249.458709990858\\
3	250.028158777026\\
};
\addlegendentry{$1.3\eta$}

\addplot [color=mycolor2, line width=1.0pt]
  table[row sep=crcr]{%
0	254.860670713077\\
0.157894736842105	255.287332075303\\
0.315789473684211	255.713604872895\\
0.473684210526316	256.139943691003\\
0.631578947368421	256.56606071728\\
0.789473684210526	256.991298911363\\
0.947368421052632	257.41646323123\\
1.10526315789474	257.840917301465\\
1.26315789473684	258.264365083502\\
1.42105263157895	258.687151166618\\
1.57894736842105	259.108578455799\\
1.73684210526316	259.529028648608\\
1.89473684210526	259.948024199677\\
2.05263157894737	260.36542697812\\
2.21052631578947	260.781150900994\\
2.36842105263158	261.195573957945\\
2.52631578947368	261.607622292257\\
2.68421052631579	262.017912774861\\
2.84210526315789	262.426152837562\\
3	262.831828421841\\
};
\addlegendentry{$\eta$}

\addplot [color=black, line width=1.0pt]
  table[row sep=crcr]{%
0	268.476381230608\\
0.157894736842105	268.737498179687\\
0.315789473684211	268.999965708178\\
0.473684210526316	269.260226023615\\
0.631578947368421	269.520099545714\\
0.789473684210526	269.778642531359\\
0.947368421052632	270.040072805537\\
1.10526315789474	270.302116619137\\
1.26315789473684	270.549665795697\\
1.42105263157895	270.804751101435\\
1.57894736842105	271.058296939841\\
1.73684210526316	271.311096114055\\
1.89473684210526	271.562210195424\\
2.05263157894737	271.812996949562\\
2.21052631578947	272.061696152747\\
2.36842105263158	272.309102199975\\
2.52631578947368	272.555177736089\\
2.68421052631579	272.799881418449\\
2.84210526315789	273.042996344137\\
3	273.284792377884\\
};
\addlegendentry{$0.7\eta$}

\end{axis}
\end{tikzpicture}%
\vspace*{-0.3in}
\caption{Number of idle drivers in southwest corner under the one-directional cordon charge for distinct $\eta$.} 
\label{figure19_2z}
\end{minipage}
\begin{minipage}[b]{0.005\linewidth}
\hfill
\end{minipage}
\begin{minipage}[b]{0.32\linewidth}
\centering
%
%
\definecolor{mycolor1}{rgb}{0.00000,0.44700,0.74100}%
\definecolor{mycolor2}{rgb}{0.85000,0.32500,0.09800}%
\begin{tikzpicture}

\pgfplotsset{every axis y label/.style={
at={(-0.5,0.5)},
xshift=30pt,
rotate=90}}

\begin{axis}[%
width=1.694in,
height=1.03in,
at={(1.358in,0.0in)},
scale only axis,
xmin=0,
xmax=3,
xlabel style={font=\color{white!15!black}},
xlabel={One-directional cordon},
ymin=210,
ymax=320,
ytick={800,900, 1000, 1100},
yticklabels={{800},{900},{1K}, {1.1K}},
ylabel style={font=\color{white!15!black}},
ylabel={Idle driver in $\mathcal{S}$},
axis background/.style={fill=white},
legend style={at={(0.0,0.54)}, anchor=south west, legend cell align=left, align=left, font=\small, draw=white!15!black, nodes={scale=0.8, transform shape} }
]
\addplot [color=mycolor1, line width=1.0pt]
  table[row sep=crcr]{%
0	283.938106581046\\
0.157894736842105	284.302933533503\\
0.315789473684211	284.668084353545\\
0.473684210526316	285.033240521932\\
0.631578947368421	285.398264861325\\
0.789473684210526	285.762803520687\\
0.947368421052632	286.128640566076\\
1.10526315789474	286.491104363326\\
1.26315789473684	286.85415088794\\
1.42105263157895	287.216224542268\\
1.57894736842105	287.577475694769\\
1.73684210526316	287.937261421654\\
1.89473684210526	288.297457326315\\
2.05263157894737	288.652910784252\\
2.21052631578947	289.008810821091\\
2.36842105263158	289.36249354916\\
2.52631578947368	289.714448280547\\
2.68421052631579	290.06437919903\\
2.84210526315789	290.417715738071\\
3	290.757945342225\\
};
\addlegendentry{$1.3N_0$}

\addplot [color=mycolor2, line width=1.0pt]
  table[row sep=crcr]{%
0	254.860670713077\\
0.157894736842105	255.287332075303\\
0.315789473684211	255.713604872895\\
0.473684210526316	256.139943691003\\
0.631578947368421	256.56606071728\\
0.789473684210526	256.991298911363\\
0.947368421052632	257.41646323123\\
1.10526315789474	257.840917301465\\
1.26315789473684	258.264365083502\\
1.42105263157895	258.687151166618\\
1.57894736842105	259.108578455799\\
1.73684210526316	259.529028648608\\
1.89473684210526	259.948024199677\\
2.05263157894737	260.36542697812\\
2.21052631578947	260.781150900994\\
2.36842105263158	261.195573957945\\
2.52631578947368	261.607622292257\\
2.68421052631579	262.017912774861\\
2.84210526315789	262.426152837562\\
3	262.831828421841\\
};
\addlegendentry{$N_0$}

\addplot [color=black, line width=1.0pt]
  table[row sep=crcr]{%
0	227.935556194348\\
0.157894736842105	228.303839855208\\
0.315789473684211	228.671645674391\\
0.473684210526316	229.038296889925\\
0.631578947368421	229.40418050269\\
0.789473684210526	229.768842181773\\
0.947368421052632	230.133056992881\\
1.10526315789474	230.495617477258\\
1.26315789473684	230.857737451076\\
1.42105263157895	231.218339585307\\
1.57894736842105	231.577894966114\\
1.73684210526316	231.937185439309\\
1.89473684210526	232.29321709697\\
2.05263157894737	232.648999879913\\
2.21052631578947	233.003527605764\\
2.36842105263158	233.356354454158\\
2.52631578947368	233.708077334224\\
2.68421052631579	234.058195895262\\
2.84210526315789	234.406848352083\\
3	234.7539123417\\
};
\addlegendentry{$0.7N_0$}

\end{axis}
\end{tikzpicture}%
\vspace*{-0.3in}
\caption{Number of idle drivers in southwest corner  under the one-directional cordon charge for distinct $N_0$.}
\label{figure20_2z}
\end{minipage}
\begin{minipage}[b]{0.005\linewidth}
\hfill
\end{minipage}
\begin{minipage}[b]{0.32\linewidth}
\centering
%
%
\definecolor{mycolor1}{rgb}{0.00000,0.44700,0.74100}%
\definecolor{mycolor2}{rgb}{0.85000,0.32500,0.09800}%
\begin{tikzpicture}

\pgfplotsset{every axis y label/.style={
at={(-0.5,0.5)},
xshift=30pt,
rotate=90}}

\begin{axis}[%
width=1.694in,
height=1.03in,
at={(1.358in,0.0in)},
scale only axis,
xmin=0,
xmax=3,
xlabel style={font=\color{white!15!black}},
xlabel={One-directional cordon},
ymin=200,
ymax=320,
ytick={800,900,1000, 1100},
yticklabels={{800},{900},{1K}, {1.1K}},
ylabel style={font=\color{white!15!black}},
ylabel={Idle driver in $\mathcal{S}$},
axis background/.style={fill=white},
legend style={at={(0.0,0.54)}, anchor=south west, legend cell align=left, align=left, font=\small, draw=white!15!black, nodes={scale=0.8, transform shape} }
]
\addplot [color=mycolor1, line width=1.0pt]
  table[row sep=crcr]{%
0	273.86176346735\\
0.157894736842105	274.524357911786\\
0.315789473684211	275.190667890564\\
0.473684210526316	275.856248788242\\
0.631578947368421	276.5224157589\\
0.789473684210526	277.188897578687\\
0.947368421052632	277.856048467475\\
1.10526315789474	278.523122113084\\
1.26315789473684	279.190154411215\\
1.42105263157895	279.854674135402\\
1.57894736842105	280.523624459645\\
1.73684210526316	281.189697891938\\
1.89473684210526	281.85531262463\\
2.05263157894737	282.519602463856\\
2.21052631578947	283.183338528396\\
2.36842105263158	283.845752073385\\
2.52631578947368	284.506817577429\\
2.68421052631579	285.166712069901\\
2.84210526315789	285.824707881999\\
3	286.480792824009\\
};
\addlegendentry{$1.3\lambda_0$}

\addplot [color=mycolor2, line width=1.0pt]
  table[row sep=crcr]{%
0	254.860670713077\\
0.157894736842105	255.287332075303\\
0.315789473684211	255.713604872895\\
0.473684210526316	256.139943691003\\
0.631578947368421	256.56606071728\\
0.789473684210526	256.991298911363\\
0.947368421052632	257.41646323123\\
1.10526315789474	257.840917301465\\
1.26315789473684	258.264365083502\\
1.42105263157895	258.687151166618\\
1.57894736842105	259.108578455799\\
1.73684210526316	259.529028648608\\
1.89473684210526	259.948024199677\\
2.05263157894737	260.36542697812\\
2.21052631578947	260.781150900994\\
2.36842105263158	261.195573957945\\
2.52631578947368	261.607622292257\\
2.68421052631579	262.017912774861\\
2.84210526315789	262.426152837562\\
3	262.831828421841\\
};
\addlegendentry{$\lambda_0$}

\addplot [color=black, line width=1.0pt]
  table[row sep=crcr]{%
0	223.652314188533\\
0.157894736842105	223.824074206392\\
0.315789473684211	223.994967029872\\
0.473684210526316	224.164931665928\\
0.631578947368421	224.333777718327\\
0.789473684210526	224.501618083278\\
0.947368421052632	224.668184007815\\
1.10526315789474	224.833763820012\\
1.26315789473684	224.997488986743\\
1.42105263157895	225.159881893358\\
1.57894736842105	225.320645951691\\
1.73684210526316	225.479886065143\\
1.89473684210526	225.637141046871\\
2.05263157894737	225.79275012805\\
2.21052631578947	225.946253252941\\
2.36842105263158	226.097985377503\\
2.52631578947368	226.247340338505\\
2.68421052631579	226.394801499641\\
2.84210526315789	226.540153258227\\
3	226.682480774277\\
};
\addlegendentry{$0.7\lambda_0$}

\end{axis}
\end{tikzpicture}%
\vspace*{-0.3in}
\caption{Number of idle drivers in southwest corner  under the one-directional cordon charge for distinct $\lambda_0$.}
\label{figure21_2z}
\end{minipage}
\end{figure*}

\begin{figure*}[h]
\begin{minipage}[b]{0.32\linewidth}
\centering
%
%
\definecolor{mycolor1}{rgb}{0.00000,0.44700,0.74100}%
\definecolor{mycolor2}{rgb}{0.85000,0.32500,0.09800}%
\begin{tikzpicture}

\begin{axis}[%
width=1.694in,
height=1.03in,
at={(1.358in,0.0in)},
scale only axis,
xmin=0,
xmax=3,
xlabel style={font=\color{white!15!black}},
xlabel={One-directional cordon},
ymin=40,
ymax=100,
ylabel style={font=\color{white!15!black}},
ylabel={$\lambda_{22}$},
axis background/.style={fill=white},
legend style={at={(0.0,0.58)}, anchor=south west, legend cell align=left, align=left, font=\small, draw=white!15!black, nodes={scale=0.8, transform shape} }
]
\addplot [color=mycolor1, line width=1.0pt]
  table[row sep=crcr]{%
0	47.7701336527144\\
0.157894736842105	47.8443260879967\\
0.315789473684211	47.9184092870213\\
0.473684210526316	47.9922642527606\\
0.631578947368421	48.0658904818193\\
0.789473684210526	48.1392504036465\\
0.947368421052632	48.2\\
1.10526315789474	48.285143947462\\
1.26315789473684	48.3576294761065\\
1.42105263157895	48.4297891444665\\
1.57894736842105	48.5016083037615\\
1.73684210526316	48.5730466256538\\
1.89473684210526	48.6441003527884\\
2.05263157894737	48.7147339809787\\
2.21052631578947	48.7849383424219\\
2.36842105263158	48.8546815927737\\
2.52631578947368	48.9239347276811\\
2.68421052631579	48.9926794070369\\
2.84210526315789	49.060895776711\\
3	49.1285436199016\\
};
\addlegendentry{$\text{1.3}\alpha$}

\addplot [color=mycolor2, line width=1.0pt]
  table[row sep=crcr]{%
0	58.5373681759301\\
0.157894736842105	58.6389711202651\\
0.315789473684211	58.7406193745077\\
0.473684210526316	58.8423080536616\\
0.631578947368421	58.9440376961458\\
0.789473684210526	59.0457698592082\\
0.947368421052632	59.1475370016456\\
1.10526315789474	59.249305161649\\
1.26315789473684	59.3510597647947\\
1.42105263157895	59.4528139552078\\
1.57894736842105	59.5545173106246\\
1.73684210526316	59.6561963248957\\
1.89473684210526	59.7578171699758\\
2.05263157894737	59.8593561661616\\
2.21052631578947	59.9608137760196\\
2.36842105263158	60.0621915307725\\
2.52631578947368	60.1634306819543\\
2.68421052631579	60.2645378782795\\
2.84210526315789	60.3655037916252\\
3	60.4662903184067\\
};
\addlegendentry{$\alpha$}

\addplot [color=black, line width=1.0pt]
  table[row sep=crcr]{%
0	70.9\\
0.157894736842105	71.0152459367859\\
0.315789473684211	71.1258980073501\\
0.473684210526316	71.2363262240283\\
0.631578947368421	71.3465515076248\\
0.789473684210526	71.4565397816694\\
0.947368421052632	71.5663009633301\\
1.10526315789474	71.6758288849303\\
1.26315789473684	71.7851085926678\\
1.42105263157895	71.8941506676245\\
1.57894736842105	72.002950669299\\
1.73684210526316	72.1114825367896\\
1.89473684210526	72.2197606335621\\
2.05263157894737	72.3277548890995\\
2.21052631578947	72.4354811071012\\
2.36842105263158	72.5429001270952\\
2.52631578947368	72.6500470131758\\
2.68421052631579	72.7568761245218\\
2.84210526315789	72.8633999097482\\
3	72.9695993324812\\
};
\addlegendentry{$\text{0.7}\alpha$}

\end{axis}
\end{tikzpicture}%
\vspace*{-0.3in}
\caption{$\lambda_{22}$  under the one-directional cordon charge for distinct $\alpha$} 
\label{figure22_2z}
\end{minipage}
\begin{minipage}[b]{0.005\linewidth}
\hfill
\end{minipage}
\begin{minipage}[b]{0.32\linewidth}
\centering
%
%
\definecolor{mycolor1}{rgb}{0.00000,0.44700,0.74100}%
\definecolor{mycolor2}{rgb}{0.85000,0.32500,0.09800}%
\begin{tikzpicture}

\begin{axis}[%
width=1.694in,
height=1.03in,
at={(1.358in,0.0in)},
scale only axis,
xmin=0,
xmax=3,
xlabel style={font=\color{white!15!black}},
xlabel={One-directional cordon},
ymin=55,
ymax=65,
ylabel style={font=\color{white!15!black}},
ylabel={$\lambda_{22}$},
axis background/.style={fill=white},
legend style={at={(0.0,0.58)}, anchor=south west, legend cell align=left, align=left, font=\small, draw=white!15!black, nodes={scale=0.8, transform shape} }
]
\addplot [color=mycolor1, line width=1.0pt]
  table[row sep=crcr]{%
0	59.5757524115305\\
0.157894736842105	59.6909319427337\\
0.315789473684211	59.8054562536329\\
0.473684210526316	59.9193606525722\\
0.631578947368421	60.032648626944\\
0.789473684210526	60.145265895237\\
0.947368421052632	60.2572425476224\\
1.10526315789474	60.368439891913\\
1.26315789473684	60.4789814153839\\
1.42105263157895	60.5886770377207\\
1.57894736842105	60.6976572441375\\
1.73684210526316	60.8058117492441\\
1.89473684210526	60.9131020994747\\
2.05263157894737	61.0195371731969\\
2.21052631578947	61.1250661858338\\
2.36842105263158	61.2296431752078\\
2.52631578947368	61.3332638103665\\
2.68421052631579	61.4358679392681\\
2.84210526315789	61.5374588117634\\
3	61.6379115878983\\
};
\addlegendentry{$1.3\epsilon$}

\addplot [color=mycolor2, line width=1.0pt]
  table[row sep=crcr]{%
0	58.5373681759301\\
0.157894736842105	58.6389711202651\\
0.315789473684211	58.7406193745077\\
0.473684210526316	58.8423080536616\\
0.631578947368421	58.9440376961458\\
0.789473684210526	59.0457698592082\\
0.947368421052632	59.1475370016456\\
1.10526315789474	59.249305161649\\
1.26315789473684	59.3510597647947\\
1.42105263157895	59.4528139552078\\
1.57894736842105	59.5545173106246\\
1.73684210526316	59.6561963248957\\
1.89473684210526	59.7578171699758\\
2.05263157894737	59.8593561661616\\
2.21052631578947	59.9608137760196\\
2.36842105263158	60.0621915307725\\
2.52631578947368	60.1634306819543\\
2.68421052631579	60.2645378782795\\
2.84210526315789	60.3655037916252\\
3	60.4662903184067\\
};
\addlegendentry{$\epsilon$}

\addplot [color=black, line width=1.0pt]
  table[row sep=crcr]{%
0	56.4962923577386\\
0.157894736842105	56.5726393638096\\
0.315789473684211	56.6494630364947\\
0.473684210526316	56.7267344408667\\
0.631578947368421	56.8044691851841\\
0.789473684210526	56.8826387569202\\
0.947368421052632	56.9612890154014\\
1.10526315789474	57.0403863095106\\
1.26315789473684	57.1199301841025\\
1.42105263157895	57.199930161351\\
1.57894736842105	57.2803988268882\\
1.73684210526316	57.3612538313521\\
1.89473684210526	57.4425759266837\\
2.05263157894737	57.5243280475513\\
2.21052631578947	57.6065129857848\\
2.36842105263158	57.6891030029584\\
2.52631578947368	57.7721171268039\\
2.68421052631579	57.8555320689673\\
2.84210526315789	57.9393471878494\\
3	58.0235461918707\\
};
\addlegendentry{$0.7\epsilon$}

\end{axis}
\end{tikzpicture}%
\vspace*{-0.3in}
\caption{$\lambda_{22}$  under the one-directional cordon charge for distinct $\epsilon$.}
\label{figure23_2z}
\end{minipage}
\begin{minipage}[b]{0.005\linewidth}
\hfill
\end{minipage}
\begin{minipage}[b]{0.32\linewidth}
\centering
%
%
\definecolor{mycolor1}{rgb}{0.00000,0.44700,0.74100}%
\definecolor{mycolor2}{rgb}{0.85000,0.32500,0.09800}%
\begin{tikzpicture}

\begin{axis}[%
width=1.694in,
height=1.03in,
at={(1.358in,0.0in)},
scale only axis,
xmin=0,
xmax=3,
xlabel style={font=\color{white!15!black}},
xlabel={One-directional cordon},
ymin=54,
ymax=68,
ylabel style={font=\color{white!15!black}},
ylabel={$\lambda_{22}$},
axis background/.style={fill=white},
legend style={at={(0.0,0.58)}, anchor=south west, legend cell align=left, align=left, font=\small, draw=white!15!black, nodes={scale=0.8, transform shape} }
]
\addplot [color=mycolor1, line width=1.0pt]
  table[row sep=crcr]{%
0	60.3743650100686\\
0.157894736842105	60.4672522431064\\
0.315789473684211	60.5601683698846\\
0.473684210526316	60.6531075919951\\
0.631578947368421	60.7460683324713\\
0.789473684210526	60.8390588713245\\
0.947368421052632	60.9320327669717\\
1.10526315789474	61.0250058367351\\
1.26315789473684	61.1179663808624\\
1.42105263157895	61.2108689791505\\
1.57894736842105	61.303792616126\\
1.73684210526316	61.3966517209502\\
1.89473684210526	61.4894426442983\\
2.05263157894737	61.5821554843294\\
2.21052631578947	61.6747765136805\\
2.36842105263158	61.7673015332691\\
2.52631578947368	61.8596949497069\\
2.68421052631579	61.9519664498718\\
2.84210526315789	62.0440734825385\\
3	62.136020001892\\
};
\addlegendentry{$1.3\sigma$}

\addplot [color=mycolor2, line width=1.0pt]
  table[row sep=crcr]{%
0	58.5373681759301\\
0.157894736842105	58.6389711202651\\
0.315789473684211	58.7406193745077\\
0.473684210526316	58.8423080536616\\
0.631578947368421	58.9440376961458\\
0.789473684210526	59.0457698592082\\
0.947368421052632	59.1475370016456\\
1.10526315789474	59.249305161649\\
1.26315789473684	59.3510597647947\\
1.42105263157895	59.4528139552078\\
1.57894736842105	59.5545173106246\\
1.73684210526316	59.6561963248957\\
1.89473684210526	59.7578171699758\\
2.05263157894737	59.8593561661616\\
2.21052631578947	59.9608137760196\\
2.36842105263158	60.0621915307725\\
2.52631578947368	60.1634306819543\\
2.68421052631579	60.2645378782795\\
2.84210526315789	60.3655037916252\\
3	60.4662903184067\\
};
\addlegendentry{$\sigma$}

\addplot [color=black, line width=1.0pt]
  table[row sep=crcr]{%
0	55.6985012775576\\
0.157894736842105	55.8128498326512\\
0.315789473684211	55.9284985235342\\
0.473684210526316	56.0436013838732\\
0.631578947368421	56.1587621419649\\
0.789473684210526	56.273955567869\\
0.947368421052632	56.3891886900806\\
1.10526315789474	56.5044543350791\\
1.26315789473684	56.6196720796121\\
1.42105263157895	56.7350005112988\\
1.57894736842105	56.8502515873873\\
1.73684210526316	56.9654694079464\\
1.89473684210526	57.0806537532063\\
2.05263157894737	57.1957750436736\\
2.21052631578947	57.3108317712749\\
2.36842105263158	57.4257911444999\\
2.52631578947368	57.5406357485526\\
2.68421052631579	57.6553638935993\\
2.84210526315789	57.7699429196812\\
3	57.8843582813011\\
};
\addlegendentry{$0.7\sigma$}

\end{axis}
\end{tikzpicture}%
\vspace*{-0.3in}
\caption{$\lambda_{22}$  under the one-directional cordon charge for distinct $\sigma$.}
\label{figure24_2z}
\end{minipage}
\end{figure*}

\begin{figure*}[h]
\begin{minipage}[b]{0.32\linewidth}
\centering
%
%
\definecolor{mycolor1}{rgb}{0.00000,0.44700,0.74100}%
\definecolor{mycolor2}{rgb}{0.85000,0.32500,0.09800}%
\begin{tikzpicture}

\begin{axis}[%
width=1.694in,
height=1.03in,
at={(1.358in,0.0in)},
scale only axis,
xmin=0,
xmax=3,
xlabel style={font=\color{white!15!black}},
xlabel={One-directional cordon},
ymin=56,
ymax=65,
ylabel style={font=\color{white!15!black}},
ylabel={$\lambda_{22}$},
axis background/.style={fill=white},
legend style={at={(0.0,0.58)}, anchor=south west, legend cell align=left, align=left, font=\small, draw=white!15!black, nodes={scale=0.8, transform shape} }
]
\addplot [color=mycolor1, line width=1.0pt]
  table[row sep=crcr]{%
0	57.3894222705599\\
0.157894736842105	57.4982458826503\\
0.315789473684211	57.6075033147778\\
0.473684210526316	57.7171935233619\\
0.631578947368421	57.8272981539965\\
0.789473684210526	57.9378319940555\\
0.947368421052632	58.0487253105156\\
1.10526315789474	58.1600139862312\\
1.26315789473684	58.2716771400301\\
1.42105263157895	58.3836647129758\\
1.57894736842105	58.4959866555411\\
1.73684210526316	58.6086160115807\\
1.89473684210526	58.7215299040452\\
2.05263157894737	58.8346295996478\\
2.21052631578947	58.9481340479186\\
2.36842105263158	59.0617761808643\\
2.52631578947368	59.1755928952943\\
2.68421052631579	59.2895849330949\\
2.84210526315789	59.4037143413057\\
3	59.5179552744702\\
};
\addlegendentry{$1.3\eta$}

\addplot [color=mycolor2, line width=1.0pt]
  table[row sep=crcr]{%
0	58.5373681759301\\
0.157894736842105	58.6389711202651\\
0.315789473684211	58.7406193745077\\
0.473684210526316	58.8423080536616\\
0.631578947368421	58.9440376961458\\
0.789473684210526	59.0457698592082\\
0.947368421052632	59.1475370016456\\
1.10526315789474	59.249305161649\\
1.26315789473684	59.3510597647947\\
1.42105263157895	59.4528139552078\\
1.57894736842105	59.5545173106246\\
1.73684210526316	59.6561963248957\\
1.89473684210526	59.7578171699758\\
2.05263157894737	59.8593561661616\\
2.21052631578947	59.9608137760196\\
2.36842105263158	60.0621915307725\\
2.52631578947368	60.1634306819543\\
2.68421052631579	60.2645378782795\\
2.84210526315789	60.3655037916252\\
3	60.4662903184067\\
};
\addlegendentry{$\eta$}

\addplot [color=black, line width=1.0pt]
  table[row sep=crcr]{%
0	59.5277158337478\\
0.157894736842105	59.6185037170312\\
0.315789473684211	59.7088741514703\\
0.473684210526316	59.7988349797998\\
0.631578947368421	59.8883984847617\\
0.789473684210526	59.9775174501579\\
0.947368421052632	60.066307894524\\
1.10526315789474	60.1548436308802\\
1.26315789473684	60.2422349849889\\
1.42105263157895	60.3295545231122\\
1.57894736842105	60.4163750693554\\
1.73684210526316	60.5027187917631\\
1.89473684210526	60.5885503634091\\
2.05263157894737	60.6738620338232\\
2.21052631578947	60.7586330542892\\
2.36842105263158	60.8428576895511\\
2.52631578947368	60.9265252192021\\
2.68421052631579	61.0096166320954\\
2.84210526315789	61.0921108821516\\
3	61.1740095561787\\
};
\addlegendentry{$0.7\eta$}

\end{axis}
\end{tikzpicture}%
\vspace*{-0.3in}
\caption{$\lambda_{22}$  under the one-directional cordon charge for distinct $\eta$.} 
\label{figure25_2z}
\end{minipage}
\begin{minipage}[b]{0.005\linewidth}
\hfill
\end{minipage}
\begin{minipage}[b]{0.32\linewidth}
\centering
%
%
\definecolor{mycolor1}{rgb}{0.00000,0.44700,0.74100}%
\definecolor{mycolor2}{rgb}{0.85000,0.32500,0.09800}%
\begin{tikzpicture}

\begin{axis}[%
width=1.694in,
height=1.03in,
at={(1.358in,0.0in)},
scale only axis,
xmin=0,
xmax=3,
xlabel style={font=\color{white!15!black}},
xlabel={One-directional cordon},
ymin=50,
ymax=75,
ylabel style={font=\color{white!15!black}},
ylabel={$\lambda_{22}$},
axis background/.style={fill=white},
legend style={at={(0.0,0.54)}, anchor=south west, legend cell align=left, align=left, font=\small, draw=white!15!black, nodes={scale=0.8, transform shape} }
]
\addplot [color=mycolor1, line width=1.0pt]
  table[row sep=crcr]{%
0	62.8769882711088\\
0.157894736842105	62.974132945353\\
0.315789473684211	63.0713548616787\\
0.473684210526316	63.1686359941223\\
0.631578947368421	63.2659703292981\\
0.789473684210526	63.3633529829937\\
0.947368421052632	63.4608323864239\\
1.10526315789474	63.5582476085255\\
1.26315789473684	63.6557324211173\\
1.42105263157895	63.7532278962102\\
1.57894736842105	63.8507359372559\\
1.73684210526316	63.948225374613\\
1.89473684210526	64.0457651766928\\
2.05263157894737	64.1431497187005\\
2.21052631578947	64.2405556796678\\
2.36842105263158	64.3378955114296\\
2.52631578947368	64.4351758818602\\
2.68421052631579	64.5323564302648\\
2.84210526315789	64.6295490539066\\
3	64.7264218730752\\
};
\addlegendentry{$1.3N_0$}

\addplot [color=mycolor2, line width=1.0pt]
  table[row sep=crcr]{%
0	58.5373681759301\\
0.157894736842105	58.6389711202651\\
0.315789473684211	58.7406193745077\\
0.473684210526316	58.8423080536616\\
0.631578947368421	58.9440376961458\\
0.789473684210526	59.0457698592082\\
0.947368421052632	59.1475370016456\\
1.10526315789474	59.249305161649\\
1.26315789473684	59.3510597647947\\
1.42105263157895	59.4528139552078\\
1.57894736842105	59.5545173106246\\
1.73684210526316	59.6561963248957\\
1.89473684210526	59.7578171699758\\
2.05263157894737	59.8593561661616\\
2.21052631578947	59.9608137760196\\
2.36842105263158	60.0621915307725\\
2.52631578947368	60.1634306819543\\
2.68421052631579	60.2645378782795\\
2.84210526315789	60.3655037916252\\
3	60.4662903184067\\
};
\addlegendentry{$N_0$}

\addplot [color=black, line width=1.0pt]
  table[row sep=crcr]{%
0	52.5496734971082\\
0.157894736842105	52.6494453467493\\
0.315789473684211	52.7487521201205\\
0.473684210526316	52.8475556094863\\
0.631578947368421	52.9458737964399\\
0.789473684210526	53.0436582025812\\
0.947368421052632	53.1409383386136\\
1.10526315789474	53.2376570549919\\
1.26315789473684	53.3338522859124\\
1.42105263157895	53.4294658999988\\
1.57894736842105	53.5244987264944\\
1.73684210526316	53.6190572458796\\
1.89473684210526	53.712774564665\\
2.05263157894737	53.8059815269348\\
2.21052631578947	53.8985480686172\\
2.36842105263158	53.9904582613373\\
2.52631578947368	54.0816946905259\\
2.68421052631579	54.1722324578648\\
2.84210526315789	54.2620740631538\\
3	54.3511861106542\\
};
\addlegendentry{$0.7N_0$}

\end{axis}
\end{tikzpicture}%
\vspace*{-0.3in}
\caption{$\lambda_{22}$ under the one-directional cordon charge for distinct $N_0$.}
\label{figure26_2z}
\end{minipage}
\begin{minipage}[b]{0.005\linewidth}
\hfill
\end{minipage}
\begin{minipage}[b]{0.32\linewidth}
\centering
%
%
\definecolor{mycolor1}{rgb}{0.00000,0.44700,0.74100}%
\definecolor{mycolor2}{rgb}{0.85000,0.32500,0.09800}%
\begin{tikzpicture}

\begin{axis}[%
width=1.694in,
height=1.03in,
at={(1.358in,0.0in)},
scale only axis,
xmin=0,
xmax=3,
xlabel style={font=\color{white!15!black}},
xlabel={One-directional cordon},
ymin=30,
ymax=100,
ylabel style={font=\color{white!15!black}},
ylabel={$\lambda_{22}$},
axis background/.style={fill=white},
legend style={at={(0.0,0.54)}, anchor=south west, legend cell align=left, align=left, font=\small, draw=white!15!black, nodes={scale=0.8, transform shape} }
]
\addplot [color=mycolor1, line width=1.0pt]
  table[row sep=crcr]{%
0	76.6358431616033\\
0.157894736842105	76.7832892842849\\
0.315789473684211	76.9310456430798\\
0.473684210526316	77.0788600736236\\
0.631578947368421	77.2268216357945\\
0.789473684210526	77.3749075510406\\
0.947368421052632	77.5231117334806\\
1.10526315789474	77.6714255841975\\
1.26315789473684	77.8198279713793\\
1.42105263157895	77.9682090843751\\
1.57894736842105	78.1168664928894\\
1.73684210526316	78.2654534524387\\
1.89473684210526	78.4140974742215\\
2.05263157894737	78.5627367453618\\
2.21052631578947	78.7113886743401\\
2.36842105263158	78.8600162285862\\
2.52631578947368	79.0085894922991\\
2.68421052631579	79.157114764172\\
2.84210526315789	79.3055423213095\\
3	79.4538564209637\\
};
\addlegendentry{$\text{1.3}\lambda{}_\text{0}$}

\addplot [color=mycolor2, line width=1.0pt]
  table[row sep=crcr]{%
0	58.5373681759301\\
0.157894736842105	58.6389711202651\\
0.315789473684211	58.7406193745077\\
0.473684210526316	58.8423080536616\\
0.631578947368421	58.9440376961458\\
0.789473684210526	59.0457698592082\\
0.947368421052632	59.1475370016456\\
1.10526315789474	59.249305161649\\
1.26315789473684	59.3510597647947\\
1.42105263157895	59.4528139552078\\
1.57894736842105	59.5545173106246\\
1.73684210526316	59.6561963248957\\
1.89473684210526	59.7578171699758\\
2.05263157894737	59.8593561661616\\
2.21052631578947	59.9608137760196\\
2.36842105263158	60.0621915307725\\
2.52631578947368	60.1634306819543\\
2.68421052631579	60.2645378782795\\
2.84210526315789	60.3655037916252\\
3	60.4662903184067\\
};
\addlegendentry{$\lambda{}_\text{0}$}

\addplot [color=black, line width=1.0pt]
  table[row sep=crcr]{%
0	39.5332918183653\\
0.157894736842105	39.593181287845\\
0.315789473684211	39.653005673831\\
0.473684210526316	39.7127655382985\\
0.631578947368421	39.7724436907141\\
0.789473684210526	39.832048287028\\
0.947368421052632	39.8915558987159\\
1.10526315789474	39.9509761586458\\
1.26315789473684	40.0102772827217\\
1.42105263157895	40.0694663690517\\
1.57894736842105	40.128525312915\\
1.73684210526316	40.187450015519\\
1.89473684210526	40.2462237201159\\
2.05263157894737	40.3048434550787\\
2.21052631578947	40.3632844633962\\
2.36842105263158	40.4215573540403\\
2.52631578947368	40.4796262351186\\
2.68421052631579	40.5374942414626\\
2.84210526315789	40.5951501884357\\
3	40.6525539564438\\
};
\addlegendentry{$\text{0.7}\lambda{}_\text{0}$}

\end{axis}
\end{tikzpicture}%
\vspace*{-0.3in}
\caption{$\lambda_{22}$ under the one-directional cordon charge for distinct $\lambda_0$.}
\label{figure27_2z}
\end{minipage}
\end{figure*}

\clearpage
\subsection*{\bf{C: Results of the multi-zone simulation under congestion charges}}
\label{AppendixD}
\begin{tabular}{ |p{1.2cm}||p{1.2cm}|p{1.2cm}|p{1cm}|p{1.25cm}|p{1.25cm}|| p{1.2cm}|p{1.2cm}|p{1cm}|p{1.25cm}|p{1.25cm}|    }
\hline
&\multicolumn{5}{|c|}{Before one-directional cordon price}& \multicolumn{5}{|c|}{After one-directional cordon price} \\ 
\hline
Zip Code& Idle vehicle & All vehicle\footnotemark & Ride fare & Ouflow   \mbox{to $\mathcal{C}$\footnotemark} 
(/min) & Ouflow to $\mathcal{R}$ (/min)& Idle vehicle & All vehicle & Ride fare  & Ouflow to $\mathcal{C}$ (/min) & Ouflow to $\mathcal{R}$  (/min) \\
\hline
94104\footnotemark              & 72.9 & 200.0 &\$2.00 &  8.14 & 1.70 & 66.6&    180.8  &  \$2.09  &   7.35  & 1.44 \\
94103              & 110.5 & 314.2 &\$2.05 &  9.55 & 3.01 &  100.5 & 285.1 &   \$2.13  &   8.73  &  2.66 \\
94109              & 114.0 & 337.6  & \$2.08 &8.77 & 4.28 & 100.4 &  299.6&    \$2.16 &    7.89  &  3.75 \\
94115              & 87.9 &  217.0 &\$2.07 &  3.33 & 5.23 &    83.6&  204.1 &   \$2.00 &   2.55  &  5.27\\
94118              & 115.4 & 237.0 &\$2.09 &  4.63 & 6.89 &   108.1&  227.4 &   \$2.01 &   3.99  &  6.98 \\
94123              & 89.9 &  184.9 &\$2.08 &  3.44 & 5.48 &    86.6 & 177.3  &  \$2.00  &  2.71  &  5.59 \\
94108\footnotemark  & 114.5 & 279.5 &\$2.04 &  10.30 & 3.41 &  105.7&   253.6&    \$2.14 &    9.37   & 2.97 \\
94121              & 47.0 &  92.7 &\$1.69 &   1.43 & 2.40 &    46.5 &  91.8 &   \$1.64   &1.26  & 2.50 \\
94102              & 117.9 & 308.9 &\$2.06 &  8.08 & 3.86 &   106.7 & 278.0  &  \$2.17  &  7.29  &  3.37 \\
94117              & 92.6 &  203.6 &\$1.83 &  2.24 & 5.73 &   90.1 &  198.6&    \$1.76 &   1.80  &  5.87 \\
94122              & 78.5 & 182.6 &\$1.82 &  1.79 & 6.65 &   79.3&   184.3 &   \$1.77&    1.58  &  6.86 \\
94114              & 92.7 &  195.9 &\$1.97 &  3.02 & 5.51 &    88.4 & 189.0 &   \$1.87   & 2.51 &   5.62 \\
94107              & 88.8 &  217.8 &\$1.87 &  6.96 & 3.24 &    81.8&  198.3&   \$1.94 &   6.36 &  2.91 \\
94110              & 110.1 & 249.6 &\$1.88 &  3.43 & 8.86 &   108.6 & 245.9 &   \$1.82 &   2.75 &  9.18 \\
94131              & 69.6 &  138.6 &\$1.68 &  1.49 & 2.99 &    68.3 &  139.1 &\$1.56 &   1.40 &  3.11 \\
94116              & 46.2 &  91.3 &\$1.57 &   0.94 & 2.15 &     46.5 &  93.2 &   \$1.49 &  0.90  &  2.27 \\
94124              & 45.5 &  81.4 &\$1.51 &   0.88 & 1.77 &     45.5 &  81.3 &   \$1.44  & 0.76 &   1.86 \\
94132              & 33.0 &  59.6 &\$1.61 &   0.06 & 1.76 &     36.2 &  65.1 &   \$1.56 &  0.06  &  1.98 \\
94112              & 51.6 &  95.6 &\$1.65 &   0.11 & 3.11 &     55.3 &  101.9 &   \$1.60  & 0.10  & 3.37 \\
\hline
\end{tabular}
\captionof{table}{Ride-sourcing market before and after the one-directional cordon price is imposed.}
\vspace{0.3cm}
\footnotetext[30]{All vehicles include idle vehicles, occupied vehicles, and vehicles on the way to pick up a designated passenger.}
\footnotetext[31]{Outflow to $\mathcal{C}$ denotes the passenger flow (per minute) with destinations in the congestion area $\mathcal{C}$.}
\footnotetext[32]{We use 94104 to represent the aggregated zone for zip code zones 94104, 94105, 94111.}
\footnotetext[33]{We use 94108 to represent the aggregated zone for zip code zones 94108 and 94133.}

\newpage
\begin{tabular}{ |p{1.2cm}||p{1.2cm}|p{1.2cm}|p{1cm}|p{1.25cm}|p{1.25cm}|| p{1.2cm}|p{1.2cm}|p{1cm}|p{1.25cm}|p{1.25cm}|    }
\hline
&\multicolumn{5}{|c|}{Before bi-directional cordon price}& \multicolumn{5}{|c|}{After bi-directional cordon price} \\ 
\hline
Zip Code& Idle vehicle & All vehicle & Ride fare & Ouflow to $\mathcal{C}$ (/min)& Ouflow to $\mathcal{R}$ (/min)& Idle vehicle & All vehicle & Ride fare  & Ouflow to $\mathcal{C}$ (/min)& Ouflow to $\mathcal{R}$ (/min) \\
\hline
94104              & 72.9 & 200.0 &\$2.00 &  8.14 & 1.70 & 70.3&    185.3  &  \$2.05  &   7.84  & 1.15 \\
94103              & 110.5 & 314.2 &\$2.05 &  9.55 & 3.01 &  105.3 & 285.6 &   \$2.06  &   9.34  &  2.17 \\
94109              & 114.0 & 337.6  & \$2.08 &8.77 & 4.28 & 102.0 &  293.3&    \$2.04 &    8.52  &  3.11 \\
94115              & 87.9 &  217.0 &\$2.07 &  3.33 & 5.23 &    72.1&  179.9 &   \$2.09 &   2.12  &  4.57\\
94118              & 115.4 & 237.0 &\$2.09 &  4.63 & 6.89 &   101.0&  211.9 &   \$2.07 &   3.59  &  6.60 \\
94123              & 89.9 &  184.9 &\$2.08 &  3.44 & 5.48 &    78.4 & 160.8  &  \$2.09  &  2.32  &  5.05 \\
94108              & 114.5 & 279.5 &\$2.04 &  10.30 & 3.41 &  108.8&   256.2&    \$2.04 &    10.13   & 2.49 \\
94121              & 47.0 &  92.7 &\$1.69 &   1.43 & 2.40 &    43.2 &  84.8 &   \$1.68   &1.11  & 2.28 \\
94102              & 117.9 & 308.9 &\$2.06 &  8.08 & 3.86 &   110.2 & 277.7  &  \$2.03  &  8.05  &  2.88 \\
94117              & 92.6 &  203.6 &\$1.83 &  2.24 & 5.73 &   84.4 &  183.2&    \$1.85 &   1.51  &  5.34 \\
94122              & 78.5 & 182.6 &\$1.82 &  1.79 & 6.65 &   75.9&   173.8 &   \$1.83&    1.35  &  6.51 \\
94114              & 92.7 &  195.9 &\$1.97 &  3.02 & 5.51 &    83.4 & 175.8 &   \$1.97   & 2.19 &   5.20 \\
94107              & 88.8 &  217.8 &\$1.87 &  6.96 & 3.24 &    88.1&  208.1&   \$1.84 &   7.25 &  2.54 \\
94110              & 110.1 & 249.6 &\$1.88 &  3.43 & 8.86 &   99.8 & 226.4 &   \$1.88 &   2.41 &  8.43 \\
94131              & 69.6 &  138.6 &\$1.68 &  1.49 & 2.99 &    66.6 &  133.2 &\$1.64 &   1.28 &  2.98 \\
94116              & 46.2 &  91.3 &\$1.57 &   0.94 & 2.15 &     45.1 &  88.9 &   \$1.54 &  0.81  &  2.16 \\
94124              & 45.5 &  81.4 &\$1.51 &   0.88 & 1.77 &     43.5 &  77.4 &   \$1.47  & 0.70 &   1.78 \\
94132              & 33.0 &  59.6 &\$1.61 &   0.06 & 1.76 &     35.0 &  62.1 &   \$1.60 &  0.05  &  1.85 \\
94112              & 51.6 &  95.6 &\$1.65 &   0.11 & 3.11 &     53.5 &   97.6 &   \$1.64  & 0.08  & 3.20 \\
\hline
\end{tabular}
\captionof{table}{Ride-sourcing market before and after the bi-directional cordon price is imposed.}
\vspace{0.5cm}

\begin{tabular}{ |p{1.2cm}||p{1.2cm}|p{1.2cm}|p{1cm}|p{1.25cm}|p{1.25cm}|| p{1.2cm}|p{1.2cm}|p{1cm}|p{1.25cm}|p{1.25cm}|    }
\hline
&\multicolumn{5}{|c|}{Before trip-based charge}& \multicolumn{5}{|c|}{After trip-based charge} \\ 
\hline
Zip Code& Idle vehicle & All vehicle & Ride fare & Ouflow to $\mathcal{C}$ (/min) & Ouflow to $\mathcal{R}$ (/min) & Idle vehicle & All vehicle & Ride fare  & Ouflow to $\mathcal{C}$ (/min) & Ouflow to $\mathcal{R}$ (/min)  \\
\hline
94104              & 72.9 & 200.0 &\$2.00 &  8.14 & 1.70 & 58.0&    155.8  &  \$1.90  &   5.47  & 1.26 \\
94103              & 110.5 & 314.2 &\$2.05 &  9.55 & 3.01 &  92.5 & 255.8 &   \$1.95  &   6.96  &  2.28 \\
94109              & 114.0 & 337.6  & \$2.08 &8.77 & 4.28 & 96.1 &  275.6&    \$1.98 &    6.45  &  3.22 \\
94115              & 87.9 &  217.0 &\$2.07 &  3.33 & 5.23 &    74.7&  179.5 &   \$1.96 &   2.46  &  3.84\\
94118              & 115.4 & 237.0 &\$2.09 &  4.63 & 6.89 &   102.8&  208.8 &   \$1.98 &   4.04  &  5.71 \\
94123              & 89.9 &  184.9 &\$2.08 &  3.44 & 5.48 &    78.9 & 157.7  &  \$1.96  &  2.67  &  4.23 \\
94108              & 114.5 & 279.5 &\$2.04 &  10.30 & 3.41 &  102.3&   235.3&    \$1.92 &    7.99   & 2.83 \\
94121              & 47.0 &  92.7 &\$1.69 &   1.43 & 2.40 &    44.9 &  84.3 &   \$1.60   &1.32  & 1.96 \\
94102              & 117.9 & 308.9 &\$2.06 &  8.08 & 3.86 &   104.4 & 260.2  &  \$1.93  &  6.27  &  3.05 \\
94117              & 92.6 &  203.6 &\$1.83 &  2.24 & 5.73 &   83.8 &  177.1&    \$1.69 &   1.88  &  4.50 \\
94122              & 78.5 & 182.6 &\$1.82 &  1.79 & 6.65 &   74.0&   166.1&   \$1.72&    1.67  &  5.65 \\
94114              & 92.7 &  195.9 &\$1.97 &  3.02 & 5.51 &    83.9 & 171.2 &   \$1.83   & 2.58 &   4.40 \\
94107              & 88.8 &  217.8 &\$1.87 &  6.96 & 3.24 &    80.4&  186.2&   \$1.77 &   5.53 &  2.65 \\
94110              & 110.1 & 249.6 &\$1.88 &  3.43 & 8.86 &   100.2 & 218.9 &   \$1.74 &   2.87 &  7.14 \\
94131              & 69.6 &  138.6 &\$1.68 &  1.49 & 2.99 &    65.3 &  127.7 &\$1.54 &   1.41 &  2.57 \\
94116              & 46.2 &  91.3 &\$1.57 &   0.94 & 2.15 &     45.3 &  86.1 &   \$1.47 &  0.93  &  1.86 \\
94124              & 45.5 &  81.4 &\$1.51 &   0.88 & 1.77 &     45.1 &  78.1 &   \$1.38  & 0.85 &   1.61 \\
94132              & 33.0 &  59.6 &\$1.61 &   0.06 & 1.76 &     32.1 &  56.4 &   \$1.45 &  0.08  &  1.58 \\
94112              & 51.6 &  95.6 &\$1.65 &   0.11 & 3.11 &     49.3 &  89.0 &   \$1.48  & 0.12  & 2.77 \\
\hline
\end{tabular}
\captionof{table}{Ride-sourcing market before and after the trip-based charge is imposed.}
\vspace{0.5cm}

\end{document}